\documentclass[12pt]{elsarticle}
\usepackage{url} 
\usepackage{lineno}
\usepackage{soul}
\usepackage{comment}
\usepackage{algpseudocode}
\usepackage{algorithm}
\usepackage{graphicx}
\usepackage{color}
\usepackage{amsmath,amssymb,amsfonts,bm}
\usepackage{lineno}
\usepackage{caption}
\usepackage{subcaption}
\usepackage{float}
\usepackage[normalem]{ulem}
\DeclareMathOperator*{\argmin}{arg\,min}

\usepackage{amsfonts}       
\renewcommand{\vec}[1]{\boldsymbol{#1}}

\usepackage{color}
\usepackage{amsmath}

\begin{document}


\title{Physics-Informed Neural Network Method for Parabolic Differential Equations with Sharply Perturbed Initial Conditions }

\author[UIUC]{Yifei Zong}
\author[UM]{QiZhi He}

\author[UIUC,PNNL]{Alexandre M. Tartakovsky\corref{mycorrespondingauthor}}
\cortext[mycorrespondingauthor]{Corresponding author}
\ead{amt1998@illinois.edu}

\address[UIUC]{Department of Civil and Environmental Engineering, University of Illinois Urbana-Champaign, Urbana, IL 61801}
\address[UM]{Department of Civil, Environmental, and Geo- Engineering, University of Minnesota, Minneapolis, MN 55455}
\address[PNNL]{Pacific Northwest National Laboratory, Richland, WA 99352.}%


\begin{abstract}
In this paper, we develop a physics-informed neural network (PINN) model for parabolic problems with a sharply perturbed initial condition. As an example of a parabolic problem, we consider the advection-dispersion equation (ADE) with a point (Gaussian) source initial condition. In the $d$-dimensional ADE, perturbations in the initial condition decay  with time $t$ as $t^{-d/2}$, which can cause a large approximation error in the PINN solution. Localized large gradients in the ADE solution make the (common in PINN) Latin hypercube sampling of the equation's residual highly inefficient. Finally, the PINN solution of parabolic equations is sensitive to the choice of weights in the loss function. We propose a normalized form of ADE where the initial perturbation of the solution does not decrease in amplitude and demonstrate that this normalization significantly reduces the PINN approximation error. We propose criteria for weights in the loss function that produce a more accurate PINN solution than those obtained with the weights selected via other methods. Finally, we proposed an adaptive sampling scheme that significantly reduces the PINN solution error for the same number of the sampling (residual) points. We demonstrate  the accuracy of the proposed PINN model for forward, inverse, and backward ADEs.   
\end{abstract}

\maketitle
\section{Introduction}\label{sec:intro}
There is an increased interest in using machine learning algorithms for solving partial differential equations (PDEs). Here, we are interested in methods where the solution of a PDE is found as a global optimization problem, where the state variables are approximated with deep neural networks (DNNs) \cite{raissi2019physics} or Karhunen-Lo{\`e}ve expansion  \cite{TARTAKOVSKY2021109904,yeung2021physics}. One example of such methods is the physics-informed neural network (PINN) method \cite{raissi2019physics,karniadakis2021physics}. 

PINN is typically introduced as a method where the DNNs are trained using data and a penalty term  in the form of the square of the PDE residuals evaluated at a set of \emph{residual} points ( \cite{karniadakis2021physics,cuomo,cai2022physics}). To facilitate this discussion, we consider the general PDE problem  
\begin{eqnarray}\label{eq:general PDE}
\mathcal{L}(u(\mathbf{x}, t), \mathbf{v}(\mathbf{x})) &=& 0, \mathbf{x} \in \Omega, t\in[0,T]\\
\mathcal{B}(u(\mathbf{x}, t), \mathbf{v}(\mathbf{x})) &=& f(\mathbf{x}, t), \mathbf{x} \in \partial\Omega\\
u(\mathbf{x}, t = t^*) &=& U(\mathbf{x}),
\label{eq:ITC}
\end{eqnarray}
where $\mathcal{L}$ is a differential operator, $\mathcal{B}$ is the  boundary condition operator, $\mathbf{x}$ is the spatial coordinate vector, $t$ is time, $\Omega$ is a bounded domain in $R^n$ with  the boundary  $\partial\Omega$, $u(\mathbf{x},t)$ is the state variable, $\mathbf{v}(\mathbf{x})$ is a vector or scalar space-dependent parameter, and $U(\mathbf{x})$ is a known function. 
Here, we are mainly concerned with the advection-dispersion equation (ADE), where  
\begin{equation}
\mathcal{L}(u(\mathbf{x}, t), \mathbf{v}(\mathbf{x})) = 
u_t + \mathbf{p}(\mathbf{x}) \cdot \nabla   u(\mathbf{x}) - \nabla \cdot [ \mathbf{D} \nabla u(\mathbf{x})],
\label{eq:ADE_int}
\end{equation}
$\mathbf{v}(\mathbf{x})$ is the divergence-free velocity field, which could be given or found by solving the Darcy flow equation, and $\mathbf{D}$ is the dispersion coefficient. For $t^*=0$, Eq \eqref{eq:ITC} provides an initial condition for the \emph{forward} ADE. With  $t^*=T>0$,  Eq \eqref{eq:ITC} gives a terminal condition for the \emph{backward} ADE.

In the PINN method, the state variable $u(\mathbf{x},t)$ is represented with a deep neural network $\hat{u}(\mathbf{x},t;\gamma)$, where  $\gamma$ is the assembly of weights and biases of the DNN. The function $\hat{u}$ should be a (nonlinear) non-polynomial function of $(\mathbf{x},t;\gamma)$, a condition that is required for a DNN with a sufficiently large  size (a large $\gamma$ vector) to accurately approximate any function and its derivatives \cite{pinkas1995reasoning,barron1994approximation,wang2018exponential}. Unlike the polynomial piece-wise representations in the finite-element (FE) methods, DNN $\hat{u}(\mathbf{x},t;\gamma)$ provides a global approximation of $u(\mathbf{x},t)$ in the space-time domain $\mathbf{z}=(\mathbf{x},t)$, and the solution of a PDE can be found by solving a global (in time and space) residual least-square minimization problem for DNN parameters $\gamma$: 
\begin{eqnarray} \label{LLR}
\gamma^* &=& \argmin_\gamma \Big[
\frac{1}{||\Omega|| T }\int_\Omega \int_0^T \mathcal{L}(\hat{u}(\mathbf{x},t;\gamma))^2 dt d\mathbf{x} \\ \nonumber
&+& \frac{1}{||\partial \Omega || T } \int_{\partial \Omega} \int_0^T [\mathcal{B}(\hat{u}(\mathbf{x},t;\gamma)) - f(x,t)]^2 dt d\mathbf{x} \\ \nonumber
&+& \frac{1}{||\Omega|| }\int_\Omega  [\hat{u}(\mathbf{x},t=t^*;\gamma) - U(\mathbf{x}) ]^2  d\mathbf{x}  \Big]
\end{eqnarray}
In general, the integrals in this minimization problem cannot be analytically computed due to the complexity of the DNN approximation. In the PINN method, these integrals are replaced  with summations as
\begin{eqnarray} \label{DLLR}
\gamma^* &=& \argmin_\gamma \Big[
\frac{1}{N_\Omega N_T } \sum_i \sum_j \mathcal{L}(\hat{u}(\mathbf{x}_i,t_j;\gamma))^2  \\ \nonumber
&+& \frac{1}{N_{\partial \Omega} N_T}\sum_i \sum_j [\mathcal{B}(\hat{u}(\mathbf{x}_i,t_j;\gamma)) - f(x_i,t_j)]^2 \\ \nonumber
&+& \frac{1}{N_\Omega} \sum_i [\hat{u}(\mathbf{x}_i,t=t^*;\gamma) - U(\mathbf{x}_i) ]^2 \Big],
\end{eqnarray}
where $N_\Omega$, $N_{\partial \Omega}$, and $N_T$ are the number of points on  $\Omega$,  $\partial \Omega$, and [0,T], respectively. 
The summations in Eq \eqref{DLLR} can be considered as approximations of the corresponding integrals in Eq \eqref{LLR} using the mid-point  rule under the assumptions that the integration domains are discretized with equal-size non-overlapping subdomains and  the terms in the summations are estimated in the middle of these subdomains. Alternatively,  terms in the summations can be estimated at points randomly placed in the corresponding domains, in which case each term in the summations must be multiplied by a weight corresponding to the size of the volume/area associated with each discretization point. However, in the PINN literature, these conditions are usually not followed. 

In the finite element least-square method (FELS), $u(\mathbf{x},t)$ in Eq \eqref{LLR} is approximated with piecewise polynomials \cite{nguyen1984space}. In the FELS method, for linear PDEs such as ADE, the integrals in Eq \eqref{LLR} can be computed exactly for each element in the FE discretization. In \cite{nguyen1984space}, the FE equations and stability analysis are presented for the forward ADE. It follows from this analysis that the FELS equations are unconditionally stable for the forward ADE and unconditionally unstable for the backward ADE. 
 It can be shown that other FE methods as well as finite difference and finite volume schemes are also unconditionally unstable for the backward ADE. 
In this work, we demonstrate that the solution of the least-square problem that is based on the global DNN approximation of $u$ is stable under certain conditions.

The fact that the FELS formulation of linear PDEs leads to linear least-square problems is a significant advantage of this method because it guarantees a unique numerical solution of the well-posed PDE problem. On the other hand, the PINN method results in a nonlinear least-square minimization problem that might have multiple local minima. There is no guarantee that a minimization algorithm would find the global minimum even if it exists. 
Another challenge with nonlinear least-square minimization problems, which involve a loss function with multiple terms as in Eq \eqref{DLLR}, is that their solution might strongly depend on the weights associated with each term.        
To this end, a dynamic estimation of the weights during training based on the magnitude of the gradient of each term was proposed in \citep{wang2021understanding}.  For ADEs, an empirical relationship between the weights and Peclet number was developed in \cite{he2021physics}. Here, we consider an alternative approach for determining weights based on the dimensionless analysis of the ADE. 

For some initial and boundary conditions, the range/amplitude  of the solutions of parabolic equations such as the ADE  can change over time by several orders of magnitude. 
 Global approximations of such solutions in the space-time domain with a finite-size DNN might be challenging. To illustrate this point, we use a DNN to approximate the function
\begin{equation}
u(x,t)=\frac1{\sqrt{2\pi(\varepsilon^2 +2Dt)}}\exp \Big( -\frac{(x-x_0 -vt)^2}{2(\varepsilon^2+2Dt)} \Big),
\label{eq:1Dsolution}
\end{equation}
which is the solution of the one-dimensional ADE with the Gaussian source initial condition $u(x,t=0)=\frac1{\sqrt{2\pi\varepsilon^2}}\exp \Big( -\frac{(x-x_0 )^2}{2\varepsilon^2} \Big)$  ($\varepsilon$ and $x_0$ are the spread and center of the source, respectively) and the homogeneous Dirichlet BC at $\pm \infty$. Other parameters in Eq \eqref{eq:1Dsolution} are the advection velocity $v$ and the dispersion coefficient $D$. According to Eq \eqref{eq:1Dsolution}, the maximum value of $u(x,t)$ decays as $u_{max}(t)=\frac1{\sqrt{2\pi(\varepsilon^2+2Dt)}}$, while the minimum of $u(x,t)$ is bounded by 0. Here, we set $\varepsilon = 0.025$, $x_0 = 0.25$, $D = 0.093$, and $v = 3.15$, and approximate $u(x,t)$ in the space-time domain $[0,4]\times[0,1]$. We set the DNN size to five hidden layers and 60 neurons per layer and train it using samples of $u(x,t)$ on the $201\times 101$ uniform mesh in the $x-t$ space. 

\begin{figure}[htb]
	\centering
	\includegraphics[angle=0,width=0.85\textwidth]{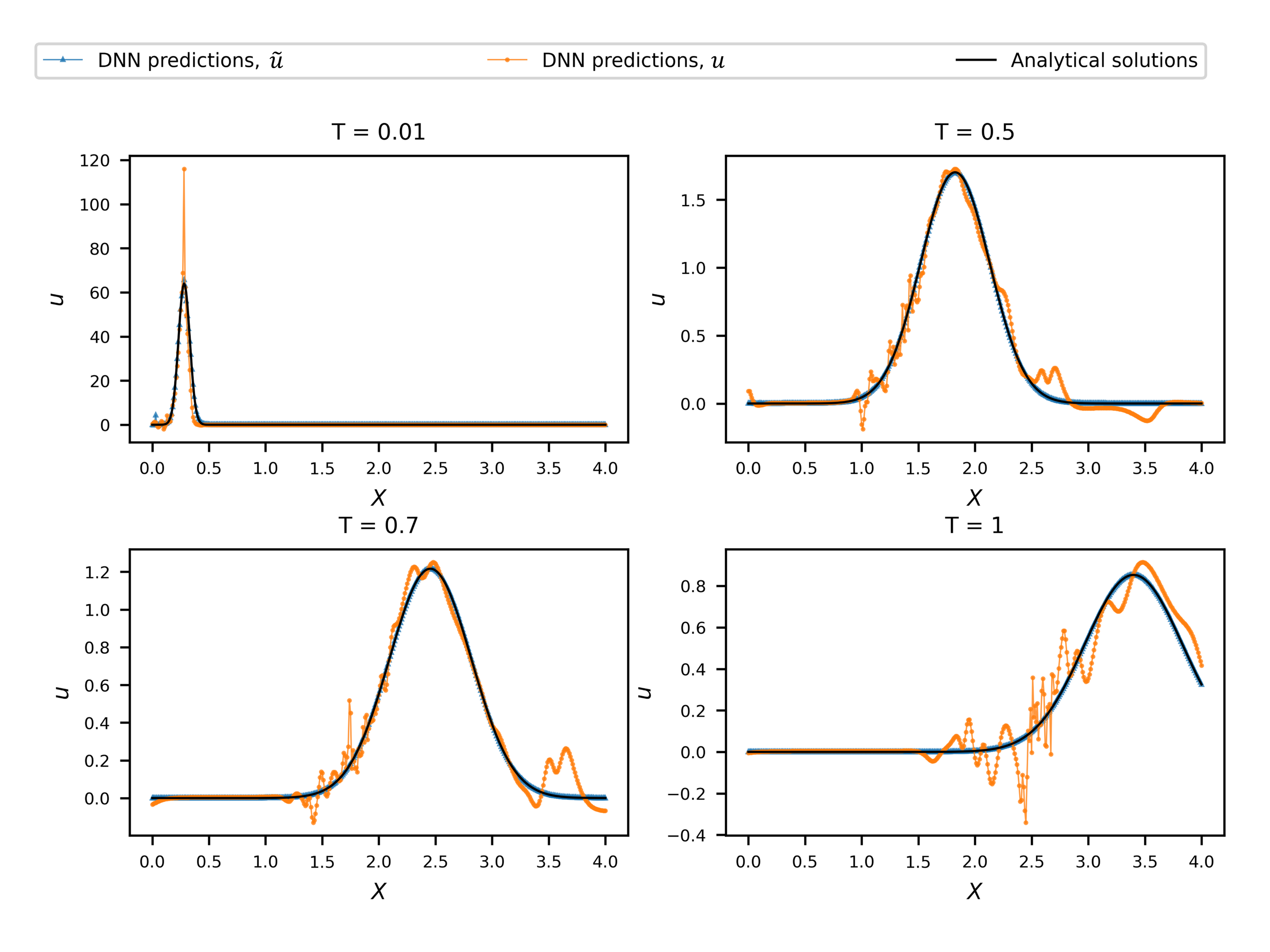}
	\caption{The reference function $u(x,t)$ and the DNN approximations 
	$u(x,t) \approx\hat{\tilde{u}}(x,t,\tilde{\gamma}) u_{max}(t)$ and $\hat{u}(x,t,\gamma) $ versus $x$
	at $t=0.01$, 0.5, 0.7, and 1. 
	}
	\label{fig:1Dcomparison}
\end{figure}

\begin{figure}[htb]
	\centering
	\includegraphics[angle=0,width=0.75\textwidth]{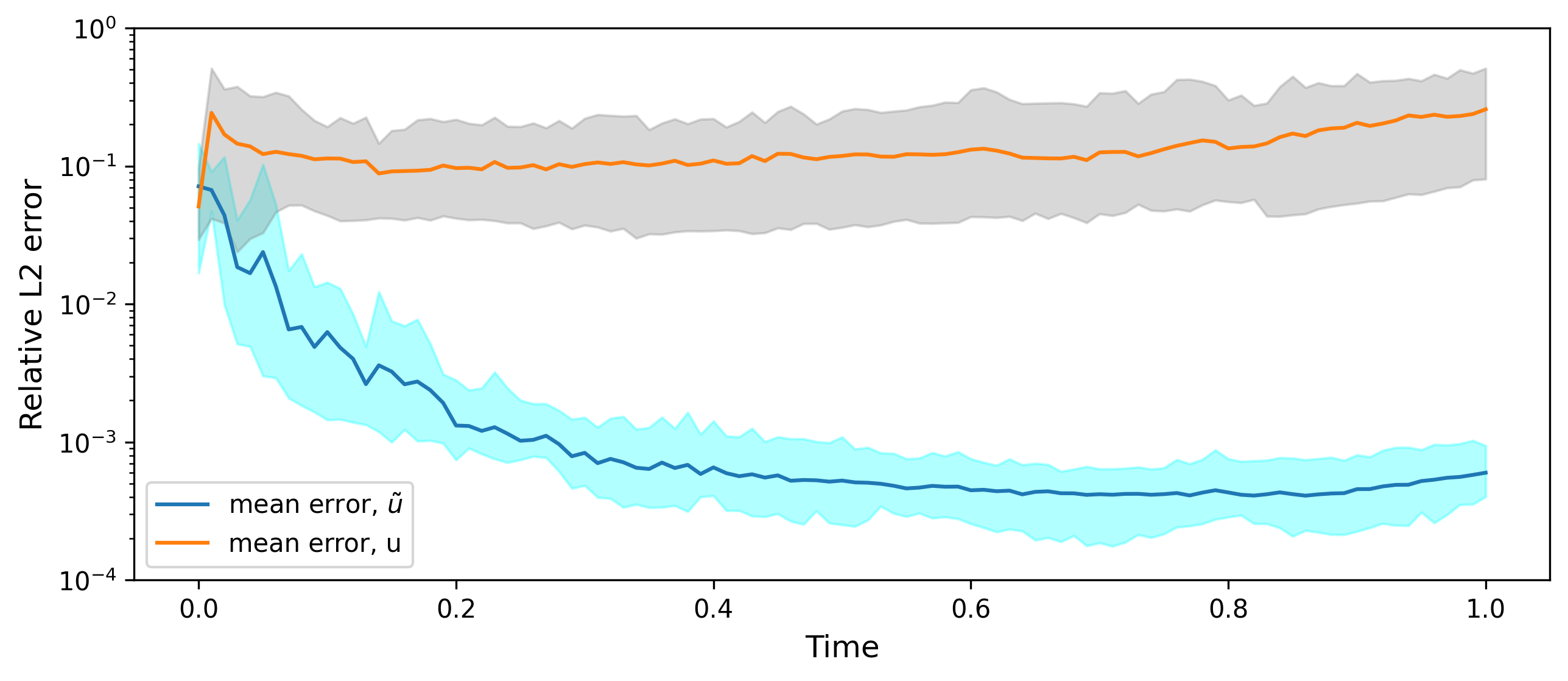}
	\caption{Average relative L2 error in the DNN approximations  of $u(x,t)$ obtained with the $u$ and normalized $\tilde{u}$ data obtained using six random initializations of the DNN weights. The upper and lower bounds of the shaded areas correspond to the maximum and minimum errors in the approximations obtained with the six DNN weight initializations.}
	\label{fig:1Derror_comparison}
\end{figure}

Figure \ref{fig:1Dcomparison} shows $u(x,t)$ and its DNN approximation $\hat{u}(x,t,\gamma)$, which is trained (i.e., the DNN parameters $\gamma$ are estimated) using $u$ samples. It can be seen that the DNN develops oscillations and, in general, is doing a poor job at approximating $u(x,t)$ at all times.  
Next, we normalize $u(x,t)$ as $\tilde{u}(x,t) = \frac{u(x,t)}{u_{max}(t)}$ such that the maximum and minimum of the normalized function are 1 and 0 at all times.  Then, we approximate ${u}(x,t)$ as ${u}(x,t)\approx \hat{\tilde{u}}(x,t,\tilde{\gamma}) u_{max}(t)$, where $\hat{\tilde{u}}(x,t,\tilde{\gamma})$  is a DNN (with the parameters $\tilde{\gamma}$) that is trained using the samples of $\tilde{u}(x,t)$. We set $\hat{\tilde{u}}(x,t,\tilde{\gamma})$ to have the same size as $\hat{{u}}(x,t,\gamma)$ and train it with the same number of samples. Figure \ref{fig:1Dcomparison} shows that $\hat{\tilde{u}}(x,t)$ provides a significantly more accurate approximation of $u$ than $\hat{u}(x,t)$. 
The DNN training is a nonlinear least-square minimization problem, whose solution (the values of the DNN's weights) can depend on the initialization of the minimization problem \cite{goodfellow2016deep, koturwar2017weight}. Figure \ref{fig:1Derror_comparison} shows the average relative L2 errors in the approximations of $u(x,t)$ obtained with the samples of $u(x,t)$ and $\tilde{u}(x,t)$ as functions of time found from six different random initializations of the DNN weights. This figure also presents the range of errors (the shaded areas) corresponding to the different DNN weights initializations.   
At $t=0$, $\tilde{u}(x,t) =  {u}(x,t)$, and the L2 errors in the two approximations are approximately the same. 
However, for the approximation $u(x,t) \approx\hat{\tilde{u}}(x,t,\tilde{\gamma}) u_{max}(t)$, the L2 error at $t=0$ is the largest and decreases with time by more than two orders of magnitude. This is because the gradient of $\tilde{u}(x,t)$ with respect to $(x,t)$ (the rate of change in space and time) is  largest at $t=0$ and decreases with time. This is not the case for the approximation based on $u(x,t) \approx \hat{u}(x,t,\gamma)$--its L2 error is smallest at $t=0$ and increases with time by one order of magnitude. Also,  Figure \ref{fig:1Derror_comparison} shows that the standard deviation of L2 error in the approximation $\hat{\tilde{u}}(x,t)u_{max}(t)$ is two orders of magnitude smaller than that in the $\hat{u}(x,t)$ approximation, i.e., a DNN approximation of  $\tilde{u}(x,t)$ is significantly less sensitive to initialization than that of ${u}(x,t)$. We note that the normalization of $u(x,t)$ with the constant $u_{max}(t=0)$, while making the normalized function to be in the [0,1] range at all times, does not significantly improve the DNN approximation of $u(x,t)$.  

The convergence of the PINN method with the number of collocation points was proven for elliptic and parabolic PDEs, given that the DNN is large enough to accurately represent the solution  \cite{shin2020convergence}. The results in Figures \ref{fig:1Dcomparison} and \ref{fig:1Derror_comparison} show that for certain initial conditions, a non-trivial normalization (i.e., normalization with a time-dependent factor rather than a constant) is required to allow an accurate representation with a DNN of a relatively small size.  

Later in this work, we demonstrate that the PINN method fails in obtaining an accurate solution of the two-dimensional ADE \eqref{eq:ADE_int} with a similar Gaussian source initial condition. Also, we demonstrate that the accuracy of the PINN method can be significantly improved by reformulating the ADE equation such that the maximum of the solution stays near 1 at all times.       

 Even when an accurate DNN representation is available, the number of the residual points required for obtaining an accurate solution can be prohibitively large. 
 The common practice is to place residual points uniformly or randomly using Latin hypercube sampling (LHS) \cite{raissi2019physics, tartakovsky2020physics,jiang2021solving}. The latter was found to require less residual points to achieve the same accuracy as in the case with uniformly placed residual points. However, for hyperbolic or advection-dominated parabolic problems, the LHS placement of residual points can be inadequate. In \cite{mao2020high-speed}, it was proposed to place the residual points adaptively according to the magnitude of the PDE residuals. Here, we propose an alternative strategy for identifying the locations of the residual points that is based on the rate of change of the solution.   
 
Despite the above-mentioned challenges,  multiple studies demonstrated that the PINN method with a sufficient number of collocation points can obtain accurate solutions of both linear and nonlinear PDEs. For example, PINN was used to solve a linear diffusion equation with a space-dependent diffusion coefficient describing saturated flow in porous media and a nonlinear diffusion equation describing flow in unsaturated porous media \cite{tartakovsky2020physics}.  
In \cite{he2020physics}, PINN was used to solve a joint inverse problem for a system of advection-diffusion and diffusion (flow) equations describing flow and transport in porous media.   In \cite{he2021physics}, some advantages of PINN over the FE method were demonstrated for a forward ADE with the direction of the advection velocity $\mathbf{v}$ being different from the grid orientation. 

We note that a comparative study of PINN with FE or other numerical methods is not the goal of this work. Our main arguments for using PINN are the backward ADE and inverse ADE problem where the unknown parameter is the hydraulic conductivity in the Darcy flow equation for the velocity $\mathbf{v}(x)$ in Eq \eqref{eq:ADE_int}. In the former case, most methods based on piece-wise continuous approximations fail. In the latter case, the standard inverse methods such as the maximum a posteriori (MAP) method suffer from the curse of dimensionality because its cost scales as $N^3$, where $N$ is the number of unknown parameters  \cite{yeung2021physics}.  Therefore, the main focus of this work is on the backward and inverse ADE problems. However, the ability of a method to accurately solve a forward ADE is a necessary condition for that method to be successful in solving inverse and backward ADEs. The proposed modifications of the PINN method improve its accuracy for the forward ADE and make it possible to apply PINN for the inverse and backward ADEs.

This paper is organized as follows.
In Section 2, we demonstrate how the PINN method can fail for the forward ADE with an instantaneous Gaussian source. We show that using either learning rate annealing \cite{wang2021understanding} or different combinations of fixed weights does not produce an accurate solution. 
In Section 3, we propose the normalized formulation of ADE, estimates of the weights in the PINN loss function, and the adaptive sampling scheme for residual points. 
In Section 4, we present forward and inverse PINN solutions of the normalized ADEs. 
We demonstrate that the assimilation of state measurements in the PINN model improves the estimations of both the state variable and parameters.
In Section 5, we discuss the PINN solutions of the backward normalized ADE equation. Conclusions are given in Section \ref{sec:conclusions}.

\section{PINN formulation for the advection-dispersion equation}\label{sec: problem setup}

\subsection{Problem statement}

We consider the two-dimensional non-reactive solute transport of a scalar field $u(\mathbf{x}, t)$ in a non-uniform flow field in a fully saturated heterogeneous porous media with space-dependent conductivity $K(\vec{x})$. 

The evolution of $u(\mathbf{x}, t)$ is governed by the ADE:
 \begin{equation}\label{eq:ADE}
 \begin{array}{ll}
 u_t +  \nabla \cdot [{\mathbf{v}} (\mathbf{x}) u(\mathbf{x})] & = \nabla \cdot [ \mathbf{D} \nabla u(\mathbf{x})], \quad \mathbf{x}\in\Omega,\quad t\in (0,T) \\
 u(\mathbf{x},t) & = g_d(\mathbf{x}), \quad \quad x_1 = 0 \\
\partial u(\mathbf{x},t) /  \partial x_1 & = g_n(\mathbf{x}), \quad \quad x_1 = L_1 \\
\partial u(\mathbf{x},t) /  \partial x_2 & = g_n(\mathbf{x}), \quad  \quad x_2 = 0 \quad{\rm{and}}\quad x_2 = L_2
\end{array}
\end{equation}

For the forward ADE, we assume the initial condition 
\begin{equation}
u(x,t=0) = \frac{M}{2\pi \varepsilon^2}\exp\left[ -\frac{(x_1 - x_1^*)^2+(x_2 - x_2^*)^2}{2\varepsilon^2}
\right], \quad \mathbf{x}\in \Omega
\label{eq:ic}
\end{equation}
which describes a Gaussian source centered at $(x_1^*,x_2^*)$ with the total mass $M$ and standard deviation $\varepsilon$. 
For the backward ADE, $u(\mathbf{x})$ is known at the terminal time $T>0$:
\begin{equation}\label{eq:terminal_cond}
u( \mathbf{x},t=T) = U( \mathbf{x}), \quad \mathbf{x}\in \Omega,
\end{equation}
where $U(x)$ is the (known) terminal condition. For both, the forward and inverse problems, we assume the homogeneous Direchlet and Neumann boundary conditions, i.e., $g_d(\mathbf{x})=g_d(\mathbf{x})=0$.

The dispersion tensor $\mathbf{D}$ is defined as (Bear, 1972):
\begin{equation}\label{eq:hydrodynamic dispersion}
\begin{split}
D_{xx} = D_w + a_L\frac{{v_1}^2}{\|v\|} + a_T\frac{{v_2}^2}{\|v\|}\\
D_{yy} = D_w + a_L\frac{{v_2}^2}{\|v\|} + a_T\frac{{v_1}^2}{\|v\|} \\
D_{xy} = D_{yx} = D_w + (a_L - a_T)\frac{v_1v_2}{\|v\|},
\end{split}
\end{equation}
where $\mathbf{v} (\mathbf{x})$ is the advection velocity given by the Darcy law: 
\begin{equation}\label{eq:velocity}
 \mathbf{v}(\mathbf{x})  = -\frac{K(\mathbf{x})}{\phi} \nabla h(\mathbf{x}),    
\end{equation}
where $\phi$ is the known porosity and $h$ is the hydraulic conductivity given by the Darcy flow equation: 
 \begin{equation}\label{eq:GWF}
 \begin{array}{ll}
 \begin{split}
 \nabla \cdot \left[ K(\mathbf{x}) \nabla h(\mathbf{x}) \right] & = 0, \quad \mathbf{x} \in \Omega \\
 h(\mathbf{x}) & = {H},\quad {x_1} = L_{1} \\
 -K(\mathbf{x}) \partial h(\mathbf{x}) /  \partial x_1 & = q, \quad x_1 = 0 \\
 -K(\mathbf{x}) \partial h(\mathbf{x}) /  \partial x_2 & = 0, \quad x_2 = 0 \:{\rm{and}}\: x_2 = L_2.
 \end{split}
 \end{array}
 \end{equation}
Here, $H$ and $q$ are the known head and flux at the Dirichlet and Neumann boundary conditions, respectively.  
We generate  $Y(\mathbf{x})=\ln K(\mathbf{x}) $  as a realization of a Gaussian random field with the constant mean and covariance function $cov(x, x^{\prime}) =\sigma^2_Y \exp(\frac{{\|x - x^{\prime}\|}^2}{\lambda^2 })$, where $\sigma^2_Y$ is the variance of and $\lambda$ is the correlation length. 

\begin{figure}[htb!]
	\centering
	\subfloat[] {\includegraphics[angle=0,width=0.45\textwidth]{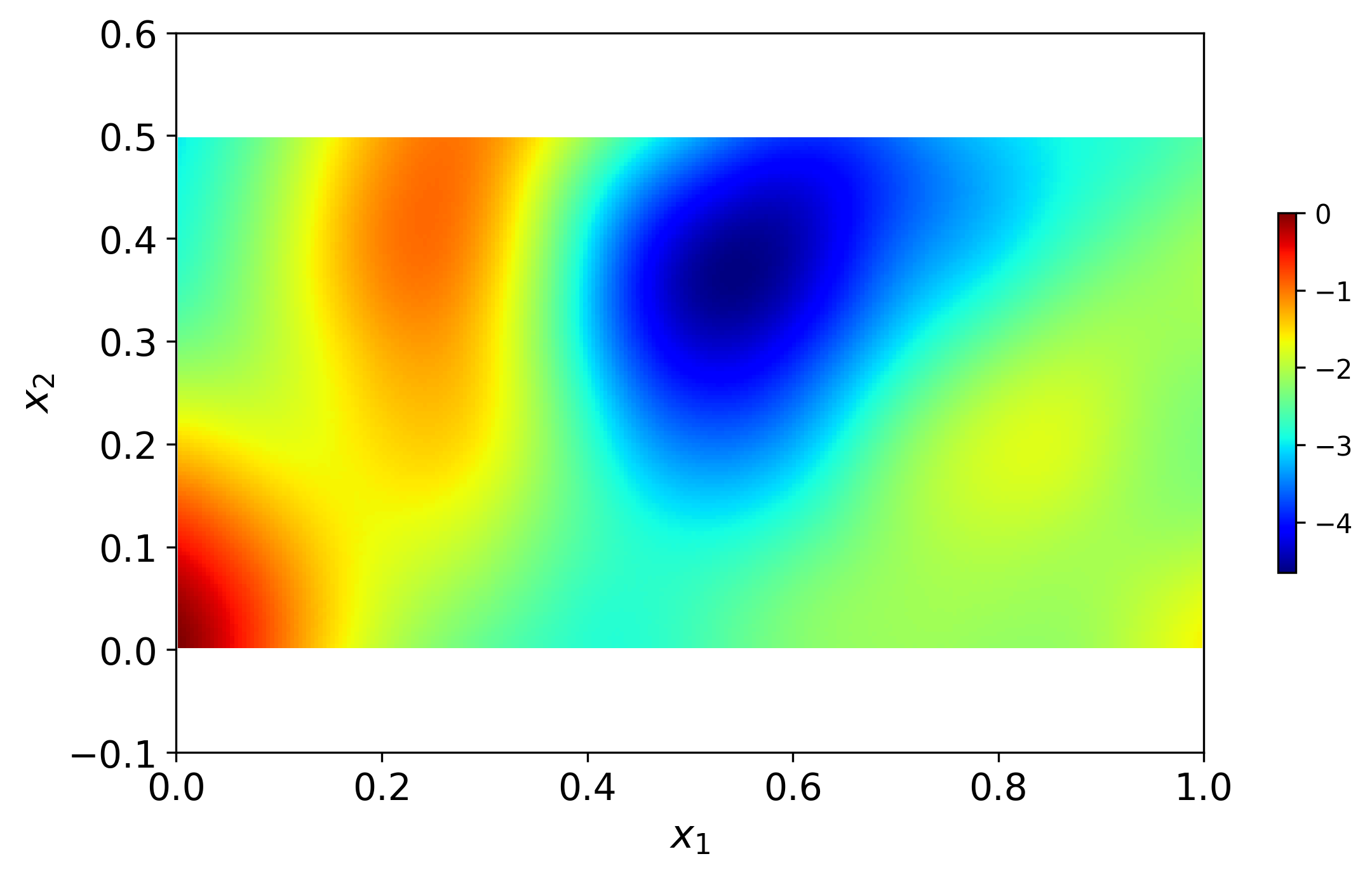}}
	\subfloat[] {\includegraphics[angle=0,width=0.45\textwidth]{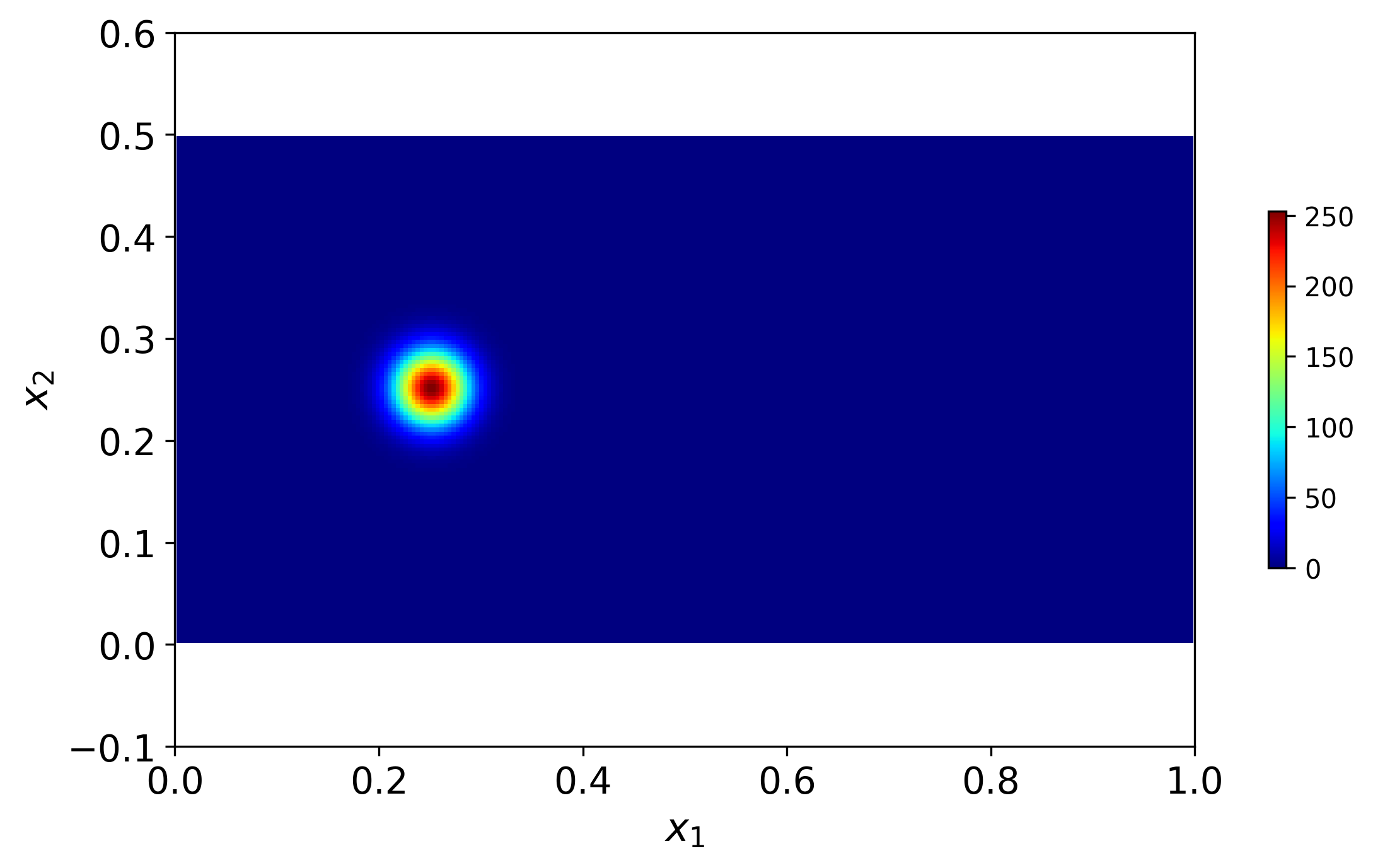}}
	\caption{\small (a) The reference field $Y(\mathbf{x})=\ln K(\mathbf{x}) $ with the correlation length $\lambda=0.5$. 
	(b) The initial condition for the solute transport problem in Eq (\ref{eq:ADE}).}
	\label{fig:K_field}
\end{figure}

In the considered examples, we set $x_1^* =  0.15$, $x_2^* = 0.25$, $M=1$, $\varepsilon = 0.025$, $\phi = 0.317$, $H=0$, $q=1$, $\lambda = 0.5$, and $\sigma_Y=0.9$. The conductivity field and initial condition for ADE are shown in Figure \ref{fig:K_field}.

To validate the PINN method, we solve the Darcy flow equation and forward ADE by the finite difference (FD) simulators MODFLOW \cite{harbaugh2005modflow} and MT3DMS \cite{zheng1999mt3dms}, respectively. The latter simulator uses the implicit Total Variation Diminishing (TVD) solver for ADE. 
The FD solution for $h(\mathbf{x})$ is obtained on a uniform $256\times128$ spatial grid and for $u(\mathbf{x},t)$ on the same grid with adaptive time stepping and is saved on a uniform $100\times256\times128$ time-space grid. 

\subsection{DNN representations in the  PINN method}\label{sec:PINN}
Here, our objective is to solve Eq \eqref{eq:ADE} with the PINN method. This equation contains the velocity $\mathbf{v}(\mathbf{x})$, which is given by Eq \eqref{eq:velocity}. In the forward and backward ADE problems, we assume that $K(\mathbf{x})$ and $h(\mathbf{x})$ in Eq \eqref{eq:velocity} are known. In the inverse ADE problem, we compute $\mathbf{v}(\mathbf{x})$ (and estimate $K(\mathbf{x})$ and $h(\mathbf{x})$) by solving a coupled inverse problem consisting of the ADE \eqref{eq:ADE} and the flow equations \eqref{eq:velocity} and \eqref{eq:GWF}. 

\subsubsection{DNN approximations of parameter fields and state variables}\label{sec:DNN_approximation}
To take advantage of automatic differentiation (AD) \cite{Baydin2015}, we approximate the state variables $u(\mathbf{x},t)$ and $h(\mathbf{x})$ and the space-dependent coefficient $K(\mathbf{x})$ with the three separate DNNs:

\begin{equation}\label{eq:NN def}
\begin{split}
K(\mathbf{x})\approx \hat{k}(\mathbf{x}; \alpha),  \\
h(\mathbf{x})\approx \hat{h}(\mathbf{x}; \beta),  \\
u(\mathbf{x}, t)\approx \hat{u}(\mathbf{x},t; \gamma), 
\end{split}
\end{equation}
where $\alpha$, $\beta$, and $\gamma$ are parameters (weights and biases) of the $\hat{k}$, $\hat{h}$, and $\hat{u}$ DNNs, respectively. Here, we use fully connected feed-forward DNNs, which have the form
\begin{equation}
    \mathcal{NN}(\mathbf{z}; \theta) := \sigma_{N_h+1}(\sigma_{N_h}(\sigma_{N_h-1}(...(\sigma_1(\mathbf{z}))))),
\end{equation}
where $\mathbf{z}$ is the input of the neural network (for $\hat{k}$ and $\hat{h}$, $\mathbf{z} = \mathbf{x}$; for $\hat{u}$, $\mathbf{z} = (\mathbf{x}, t)$), $\theta$ is the assembly of DNN parameters, $N_h$ is the number of the hidden layers, and 
\begin{equation}
    \sigma_{i}(\mathbf{z}_i) = \sigma \odot (\mathbf{W}_{i}\mathbf{z}_{i} + \mathbf{b}_{i})
\end{equation}
where $\sigma$ is the activation function, $\odot$ represents the element-wise operation, $\mathbf{z}_i = \sigma_{i-1}(\mathbf{z}_{i-1})$, and $\mathbf{z}_1 = \mathbf{z}$. Here, we use the hyperbolic tangent activation function $\sigma_i(\cdot)= \tanh(\cdot)$ for $i \in \{1, ..., N_h\}$, which satisfies the requirement of being at least twice continuously differentiable in $\mathbf{x}$, which is needed for solving the second-order PDE. 
The choice of the activation function in the output layer ($i = N_h + 1$) is problem specific. For the forward and inverse ADEs,  we use the identity activation function:
\begin{equation}\label{eq:identity_function}
\sigma_{N_h+1}(\mathbf{z}_{N_h+1}) = \mathbf{W}_{N_h+1}\mathbf{z}_{N_h+1} + \mathbf{b}_{N_h+1}.
\end{equation}
For the backward ADE, we employ a nonlinear positive (sigmoid) activation function to strictly enforce the positivity of the solution. For $i \in \{1, ..., N_{h + 1}\}$, $\mathbf{W}_{i} \in \mathbb{R}^{N_n^ {i + 1}\times{N_n^ {i}}}$ is the weight matrix connecting two adjacent layers and $\mathbf{b}_{i} \in \mathbb{R}^{N_n^ {i+1} }$ is the bias vector, where $N_n^ i$ is the number of neurons for layer $i$. 
There are other at least twice differential activation functions, including the sigmoid function, which varies from 0 to 1. We choose the tanh activation because it has better convergence properties \cite{lecun2012efficient}. 
In all DNNs,  we fix the number of hidden layers to $N_h=5$. The number of neurons in each hidden layer $N_n$ is kept the same. The effect of $N_n$ on the PINN solution accuracy is reported in Section \ref{sec:forward}.
 
\subsection{PINN formulation for forward and backward ADEs}
In groundwater models, the functional form of $K(\mathbf{x})$ is usually not known--instead, it is given at a set of points, e.g., on the FD mesh. The mesh values $K_n= K(\mathbf{x}_n^K)$ of $K$ are found by solving an inverse problem or generated as a realization of a correlated random field with known statistics. Here, we train $\hat{k}(\mathbf{x}; \alpha)$ using a set of conductivity measurements  $\{ K_n \}_{n=1}^{N_k}$, which are randomly sampled from the $K$ mesh values, by minimizing the loss function $L_k(\alpha)$,  
\begin{equation}\label{eq:k net}
\begin{split}
\alpha^* = \min_\alpha \left\{ L_{\hat{k}}(\alpha) =  \frac{1}{N_k} \sum_{n = 1}^{N_k} [\hat{k}({\mathbf{x}}_n; \alpha) - {K}_n]^2 \right\}.
\end{split}
\end{equation}
For the forward and backward ADEs, we assume that $h(\mathbf{x})$ is also known on the same mesh where $K(x)$ is specified, and we train $\hat{h}(\mathbf{x}; \beta)$ by minimizing the loss function $L_h(\beta)$
\begin{eqnarray}\label{eq:h net}
\beta^* = \min_\beta \left\{ L_{\hat{h}}(\beta) =  \frac{1}{N_h} \sum_{n = 1}^{N_h} [\hat{h}({\mathbf{x}}_n; \beta) - h_n]^2 \right\},
\end{eqnarray}
where $\{ h_n\}_{n=1}^{N_h}$ are sampled from the head values on the finite difference mesh. 
For the considered problem and DNN size, we find that the DNN approximation errors for both $K$ and $h$ are less than $0.5\%$, which is much smaller than the errors in the PINN solutions of ADEs presented below. Therefore, we disregard the DNN approximations of the known $K$ and $h$ fields as potential sources of errors in the PINN solutions. 

The $\hat{k}$ and $\hat{h}$ DNNs are used to compute a DNN approximation of $\mathbf{v}(x)$ based on the Darcy equation \eqref{eq:velocity}:
\begin{equation}\label{eq:DNN_velocity}
\hat{\mathbf{v}}(\mathbf{x};\alpha,\beta)  = -\frac{\hat{k}(\mathbf{x};\alpha)}{\phi} \nabla \hat{h}(\mathbf{x};\beta),    
\end{equation}
Here, $\alpha$ and $\beta$ are estimated in Eqs \eqref{eq:k net} and \eqref{eq:h net}, and the gradients are evaluated using AD. 

The PINN solution of the forward or backward ADE  \eqref{eq:ADE} is formulated as
\begin{eqnarray}\label{eq:u net}
 \gamma^* = \min_\gamma \Big\{ L_{\hat{u}} (\gamma) &=&  \frac{\lambda_{ic}}{N_{ic}} \sum_{n = 1}^{N_{ic}} \left[ \hat{u}({\mathbf{x}}_n, t = 0; \gamma) - u({\mathbf{x}}_n, t = 0)\right]^2  \\  \nonumber
&+& \frac{\lambda_{bcn}}{N_{bcn}^x \times N_{bcn}^t} \sum_{n = 1}^{N_{bcn}^x} \sum_{m = 1}^{N_{bcn}^t} \left[ 
\mathcal{B}_u (\mathbf{x}_n, t_m; \gamma)  \right]^2  \\ \nonumber
&+& \frac{\lambda_{bcd}}{N_{bcd}^x \times N_{bcd}^t} \sum_{n = 1}^{N_{bcd}^x} \sum_{m = 1}^{N_{bcd}^t}  [\hat{u}(\mathbf{x}_n, t_n; \gamma) - g_d(\mathbf{x}_n)]^2 \\ \nonumber
&+& \frac{\lambda_{res}}{N_{res}} \sum_{n = 1}^{N_{res}} [\mathcal{R}_{\hat{u}}(\mathbf{x}_n, t_n; \gamma) ]^2 \Big \}
\end{eqnarray}
where $\mathcal{R}_{\hat{u}}(\mathbf{x}, t; \gamma) =
\hat{u}_t (\mathbf{x}, t; \gamma) +  \nabla \cdot [\hat{\mathbf{v}} (\mathbf{x}) \hat{u}(\mathbf{x}, t; \gamma)]  - \nabla \cdot [ \mathbf{D} \nabla \hat{u}(\mathbf(\mathbf{x}, t; \gamma)]$ is the residual term, 
$\mathcal{B}_u (\mathbf{x}_n, t_m; \gamma) = \nabla \hat{u}(\mathbf{x},t_m; \gamma)
\Big \rvert_{\mathbf{x} = \mathbf{x}_n} 
\cdot \mathbf{n} - g_n(\mathbf{x}_n)$ is the Neumann boundary operator, 
 and $\mathbf{n}$ is the unit vector normal to boundaries where $\mathbf{x}_n$ is sampled.  
The coefficients $\lambda_{ic}$, $\lambda_{bcn}$, $\lambda_{bcd}$, and $\lambda_{res}$ are the weights corresponding to the initial, Neumann boundary, and  Dirichlet boundary conditions, and the ADE residual, respectively.
We randomly sample the initial concentration values at ${N_{ic}}$ spacial locations, where $N_{ic}$ equals $75\%$ of the initial nodal values on the finite volume mesh (i.e., $N_{ic} = 0.75 \times 256\times128 = 26,274$).
We enforce the boundary conditions at $N^t_{bcn} = N^t_{bcd} = N_t =100$ equal time intervals. 
At each time interval, we use the LHS scheme to sample 54 points on boundaries in the $x_2$ direction and 128 points on the boundaries in the $x_1$ direction. As result, the Dirichlet BC is enforced at $N_{bcd}^x=54$ spatial locations and the Neumann BCs are enforced at $N_{bcn}^x=128\times 2+ 54$ spatial locations. The ADE residual is enforced at $N_{res}$ points in the space-time domain. In the following sections, we study the accuracy of PINN solutions as functions of the number and distributions of the residual points. 
Algorithm \ref{alg:1} below encapsulates the complete training process. 

The DNN approximations of $u$ and $\mathbf{v}$ make the PINN ADE model ``analytically differentiable,'' and all partial derivatives  in the ADE  \eqref{eq:ADE} can be computed using AD, which is important for estimating the ADE residuals without numerically discretizing the derivatives. Also, in the PINN ADE model, the gradients with respect to the parameter vectors  $\alpha$ and $\beta$  can be computed with AD, which is needed for the back-propagation step in the DNN training.

\begin{algorithm}[htb]
 {Step 1: Initialize parameters $\alpha$, $\beta$, and $\gamma$ of three neural networks using the Xavier initialization scheme.} \\
 {Step 2: Train $\hat{k}(\mathbf{x}; \alpha)$ by minimizing the loss function (\ref{eq:k net})
 so that ${\alpha}^* = argmin L_k(\alpha)$}\\
 {Step 3: Train $\hat{h}(\mathbf{x}; \beta)$ by minimizing the loss function (\ref{eq:h net})
 so that ${\beta}^* = argmin L_h(\beta)$}\\
 {Step 4: 
 Find solution in the form $\hat{u}(\mathbf{x}, t; \gamma)$ by minimizing the loss function (\ref{eq:u net}):
 ${\gamma}^* = argmin L_u(\gamma)$}
\caption{PINN algorithm for solving forward and backward ADEs \eqref{eq:ADE}}
\label{alg:1}
\end{algorithm}

\subsection{Neural network training}\label{training}

Following  \cite{he2020physics}, we use a "two-step" approach to minimize the loss function in Eqs \eqref{eq:k net}, \eqref{eq:h net}, and \eqref{eq:u net}, where  the Adam optimizer \cite{kingma2017adam} is employed at the first step and the L-BFGS-B (quasi-Newton) optimizer \cite{liu1989limited} is used at the second step.   For training $\hat{k}(\mathbf{x}; \alpha)$ and  $\hat{h}(\mathbf{x}; \beta)$, the Adam optimizer is run for $30,000$ epochs with a fixed learning rate of $0.001$. In training $\hat{u}(\mathbf{x},t; \gamma)$, the Adam optimizer is run for $200,000$ epochs with a fixed learning rate of $0.0002$. The mini-batch size is set to $1,000$. The L-BFGS-B algorithm is terminated either when the magnitude of the loss function is less than a prescribed tolerance of $10^{-7}$ or the number of iterations (training epochs) exceeds $50,000$. 
The DNNs weights and biases are randomly initialized using the Xavier initialization \cite{glorot2010understanding}. 

All simulations in this work are performed on a PC with NVIDIA GeForce RTX 2070 Max-Q GPU. We find that the training time increases with the DNN size and the number of residual points. In the considered examples, the training time ranges from approximately 8 core hours to 12 core hours.

\subsection{Limitations of PINN for the ADE with Gaussian source initial conditions}\label{sec:PINN_limitations}

PINN was shown to be accurate for solving ADEs with certain initial and boundary conditions. For example, in \cite{he2021physics},  PINN was used to solve the two-dimensional time-dependent ADE in a non-uniform flow field 
with the zero initial condition, $u=f(x)$ was prescribed at the upstream boundary, and the homogeneous Neumann BC was assigned at the downstream boundary. The maximum and minimum of $f(x)$ are 1 and 0, correspondingly. For these BCs, 
the maximum and minimum of $u(\mathbf{x},t)$ are 1 and 0, respectively, at all times. An accurate PINN solution for this problem was obtained with the empirically determined weights $\lambda_{ic} = \lambda_{bcn} = \lambda_{bcd} = 10$ and $\lambda_{res}=1$. 

\begin{figure}[h]
	\centering
	\subfloat[] {\includegraphics[angle=0,width=0.47\textwidth]{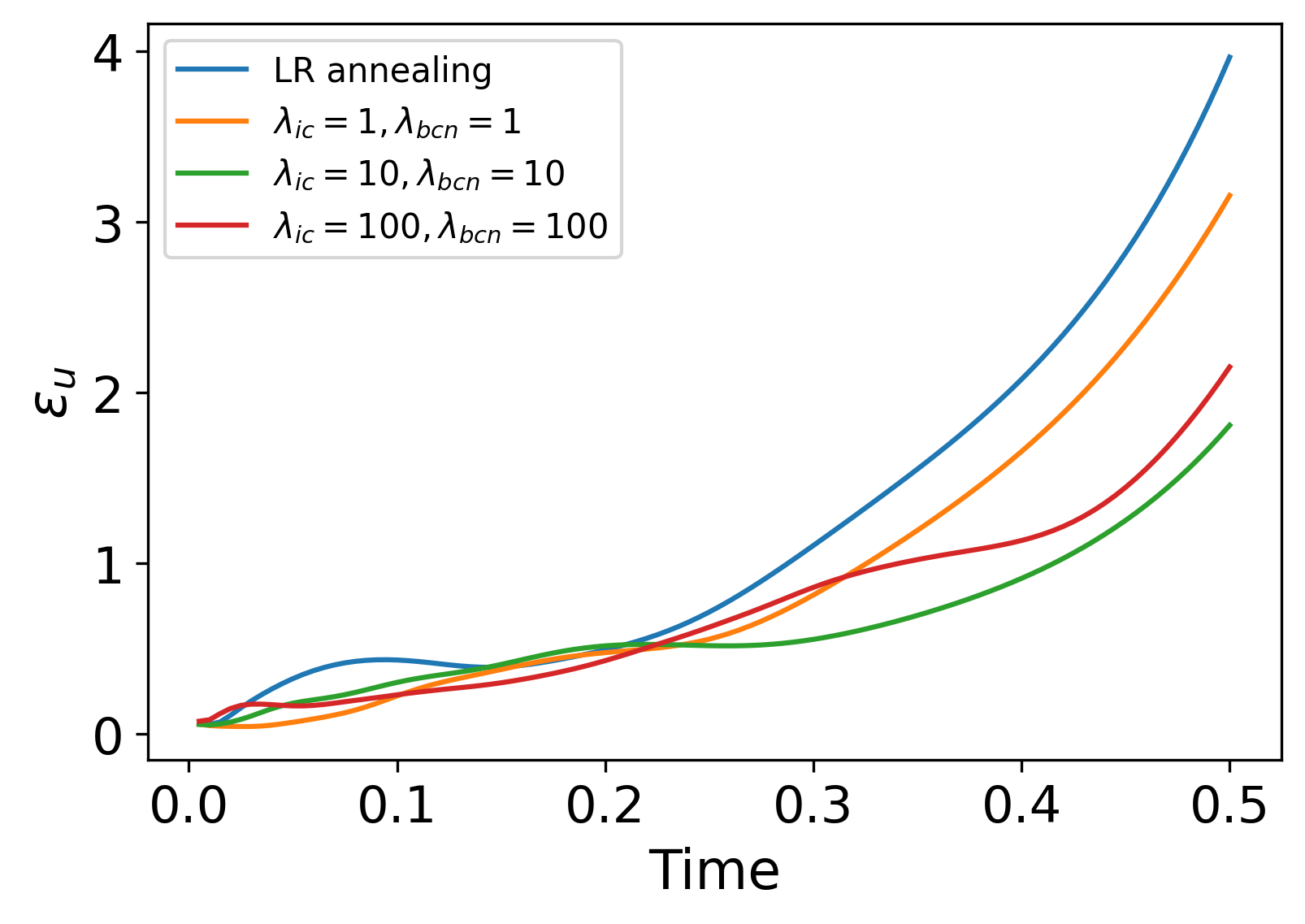}}
	\subfloat[] {\includegraphics[angle=0,width=0.5\textwidth]{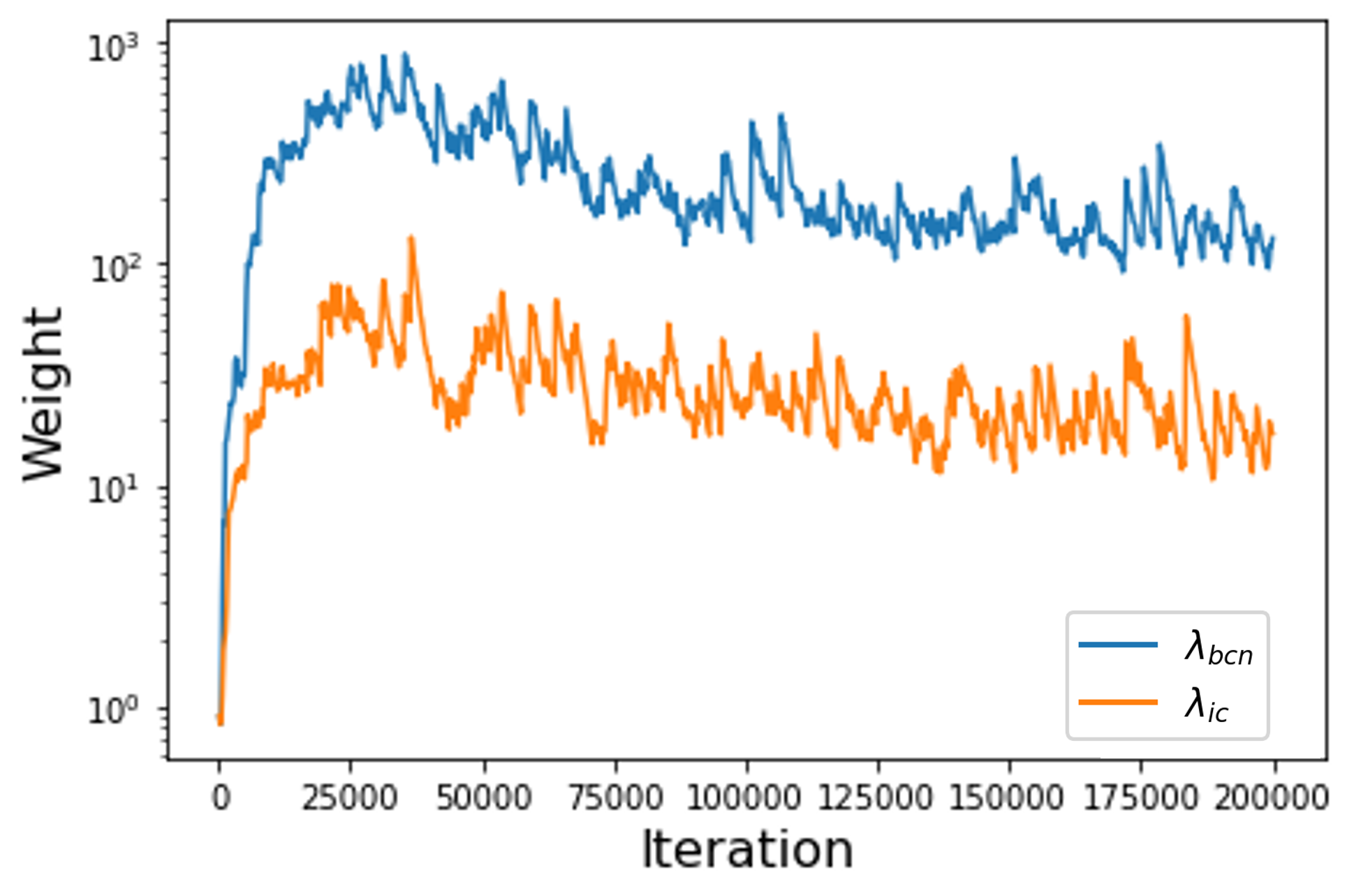}}
	\caption{ 
	(a) Relative $L_2$ error in the PINN solutions of the forward ADE \eqref{eq:ADE}
subject to the initial condition \eqref{eq:ic}
obtained with four different combinations of the weights in the loss function, including the weights selected according to the learning rate annealing method \protect\cite{wang2021understanding}. In all cases, $\lambda_{res}=1$ and $\lambda_{ic}=\lambda_{bd}$.
	(b) Evolution of the weights $\lambda_{ic}$ and $\lambda_{bn}$ (labeled as $\lambda_{bc}$) in the learning rate annealing algorithm. 
	}
	\label{fig:error_original_PINN}
\end{figure}

However, we find that applying the PINN method to the ADE \eqref{eq:ADE}  with a Gaussian source initial condition \eqref{eq:ic}
does not produce an accurate solution. The relative L2 error in this solution is shown in Figure \ref{fig:error_original_PINN}(a) for four different combinations of the weights in the loss function.   
Here and later in this work, we define the relative L2 error at time $t$ as 
\begin{equation}\label{eq:rl2_error}
\varepsilon_u(t) = \sqrt{\frac{\sum_{j=1}^{N_{tot}}\left[u(\mathbf{x}_j,t)-\hat{u}(\mathbf{x}_j,t)\right]^{2}}{\sum_{j=1}^{N_{tot}}\left[u(\mathbf{x}_j,t)\right]^{2}}} 
\end{equation}
where $N_{tot}$ is the mesh size in the reference FD solution.
For all considered  weight choices, the L2 errors increase with time, with the maximum errors ranging from 180\% to near 400\%. 
These solutions are obtained on the time domain $t\in(0,0.5]$, with $N_{res}=100,000$ residual points randomly chosen using LHS across the time-space domain and the DNN width set to $N_n = 60$. 
 The four combinations of weights are selected as follows. In all four cases, we set $\lambda_{res}=1$ and  $\lambda_{ic}=\lambda_{bcd}$. In the first case, the remaining weights $\lambda_{ic}$ and $\lambda_{bcn}$ are computed according to the learning rate annealing algorithm \cite{wang2021understanding}. This algorithm updates the weights during the training to balance the gradients of all terms in the loss function. The resulting values of $\lambda_{ic}$ and $\lambda_{bcn}$ are shown in Figure \ref{fig:error_original_PINN}(b).  In cases 2--4, we set the weights statically, i.e., the weights do not change during training. We use the weights  $\lambda_{ic} = \lambda_{bcd} = \lambda_{bcn} = 10$ from \cite{he2021physics} and two additional sets,  $\lambda_{ic} = \lambda_{bcd} = \lambda_{bcn} = 1$ and $100$. 
 The smallest L2 error corresponds to the weights selected according to  \cite{he2021physics}, and the largest errors are produced by the learning rate annealing algorithm.  In Figure \ref{fig:PINN_forward_original}, we plot four snapshots of the PINN solution obtained with $\lambda_{ic} = \lambda_{bcd} = \lambda_{bcn} = 10$. Significant point errors are observed at all times, with the PINN solution becoming negative at time $t=0.5$.  
 
The total error in a PINN solution can be decomposed in three parts: the \emph{optimization error}, the \emph{estimation error}, and the \emph{approximation error} \cite{shin2020convergence}. 
The approximation error comes from the finite-size of DNN and quantifies the distance between the class of functions that the given-size DNN can represent and the true PDE solution. The estimation error is due to 
approximating the integral loss function \eqref{LLR} with the \emph{discrete} loss function \eqref{DLLR}. In general, this error should decrease with the increasing number of residual points. The optimization error is the error in solving the nonlinear least-square problem. This error strongly depends on the weight selection in the PINN loss function. 

We attribute the high errors in the PINN solutions shown in Figure \ref{fig:error_original_PINN}(a) to several factors. First, the maximum value of the solution $u_{max}(t)$ exponentially decreases with time. This makes it difficult for a DNN with a finite size to approximate $u(\mathbf{x},t)$ in the space-time domain, as demonstrated in Figure \ref{fig:1Dcomparison} for a solution of the one-dimensional ADE. Second, the  LHS of the time-space domain with residual points in the PINN residual least-square formulation \eqref{DLLR} might not be optimal for the ADE problem where the large time derivatives exist at early times.
Finally, the selection of weights in the loss function in these solutions might not be optimal. 

We find that increasing the number of residual points to 200,000 and the DNN width to 120 does not significantly decrease the L2 error. 
Theoretically, further increasing the DNN size and number of residual points should reduce the PINN error. However, such a brute-force approach could be computationally intractable for many real-world applications. 
In the next section, we propose several improvements that reduce the L2 error to less than 6\% without increasing the DNN size and the number of residual points.

\section{Improved PINN method}\label{sec:methodology}
In this section, we propose a normalized ADE equation to reduce the approximation error in the PINN method. Next, we propose the adaptive sampling scheme for residual points to reduce the estimation error. Finally, we derive estimates of the weighting coefficients to reduce the minimization error.

\subsection{Normalized PINN}

A general rule for reducing the approximation and optimization errors in DNN models with the tanh and sigmoid activation functions is to normalize the DNN's training dataset such that inputs are in the [-0.5,0.5] range (this normalization is related to the fact that derivatives of the activation function $d \tanh(z)/dz$ are very small for $z$ outside of the [-0.5,0.5] range) and the outputs are in the [0,1] range \cite{kim1999normalization,lecun2012efficient}. Otherwise, the weights and biases of the DNN could be very large, which in turn might cause instabilities in the training process (large minimization errors). Usually, a dataset is normalized using a constant normalization factor. 
However, we demonstrated in Figure \ref{fig:1Dcomparison} that for $u(x,t)$ with exponentially decreasing $u_{max}(t) = \max_x u(x,t)$, a more accurate training is achieved if the DNN is trained with the samples of $u$ normalized with the time-dependent factor $u_{max}(t)$ rather than a constant factor, e.g., $\max_{x,t}u(x,t)$. 
Below, we propose a similar normalization of the state variable $u(\mathbf{x},t)$ in ADE \eqref{eq:ADE}. 

We normalize the space and time coordinates in Eq  \eqref{eq:ADE} as $\tilde{t} = \frac{t}{T}-0.5$, $\tilde{x}_1=\frac{x_1}{L_1} - 0.5$, and $\tilde{x}_2 = \frac{x_2}{L_2}-0.5$, where $T$ is the time horizon of the model.
Next, we introduce a dimensionless variable $\tilde{u}(\mathbf{x},t) = \frac{u(\mathbf{x},t)}{g(t)}$, where $g(t)$ is the maximum of the solution $\overline{u}(\mathbf{x},t)$ of the ``mean-field'' ADE
$$
\overline{u}_t +  \mathbf{V}(\mathbf{x}) \cdot \nabla \overline{u}(\mathbf{x})  = \nabla \cdot [ \mathbf{D} \nabla \overline{u}(\mathbf{x})] 
$$
defined on the two-dimensional infinite domain with the zero Dirichlet boundary conditions at plus/minus infinity and the initial condition \eqref{eq:ic}.
Here, $\mathbf{V}=[V,0]^T$ is the mean-field velocity that solves the flow equation  \eqref{eq:GWF} with $K(\mathbf{x})$ replaced by its geometric average. 
The analytical solution of this equation is 
\begin{equation}
\overline{u}(\mathbf{x},t) = \frac{M}
{2\pi \sqrt{\varepsilon^2 + 2 D_{x_1} t}\sqrt{\varepsilon^2 + 2 D_{x_2} t}}
\exp\left[ - \frac12 (
\frac{(x_1 -x_1^* -Vt)^2}{\varepsilon^2 +2 D_{x_1} t}
+\frac{(x_2 - x_2^*)^2}{\varepsilon^2 +2 D_{x_2} t }
)
\right].
\end{equation}
The maximum of $\overline{u}(\mathbf{x},t)$ over $\mathbf{x}$ is
\begin{equation}\label{eq:time-dependent normalizer}
g(t) = \frac{M}
{2\pi\sqrt{\varepsilon^2 + 2 D_{x_1} t}
 \sqrt{\varepsilon^2 + 2 D_{x_2} t}
}.
\end{equation}

The maximum of  $u(\mathbf{x},t)$ (the solution of Eq \eqref{eq:ADE}) over $\mathbf{x}$ is well approximated by $g(t)$ as long as $x_1^* + Vt < L_1$ (i.e., as long as the center of the plume stays within the domain). By introducing the additional dimensionless variables 
$\tilde{v}_{1} = v_{1} \frac{T}{L_{1}}$,  $\tilde{v}_{2} = v_{2} \frac{T}{L_{2}}$, $\tilde{D}_{11} =  D_{11} \frac{T}{L_{1}^2}$ , $\tilde{D}_{22} =  D_{22} \frac{T}{L_{2}^2}$, $\tilde{D}_{12} =  D_{12} \frac{T}{L_{1}L_{2}}$, and $\tilde{D}_{21} =  D_{21} \frac{T}{L_{2}L_{1}}$, we can rewrite ADE \eqref{eq:ADE} in a dimensionless form as:

\begin{equation}\label{eq:normalized_ADE}
\begin{split}
\frac{\partial \tilde{u}}{\partial \tilde{t}} + 
 \tilde{v}_{1} \frac{\partial \tilde{u}}{\partial \tilde{x}_1} +
 \tilde{v}_{2} \frac{\partial \tilde{u}}{\partial \tilde{x}_2}
 =
\frac{\partial}{\partial \tilde{x}_1}[{\tilde{D}_{11}}\frac{\partial \tilde{u}}{\partial \tilde{x}_1}] + 
\frac{\partial}{\partial \tilde{x}_2}[{\tilde{D}_{22}}\frac{\partial \tilde{u}}{\partial \tilde{x}_2}] \\
+ \frac{\partial}{\partial \tilde{x}_2}[{\tilde{D}_{12}}\frac{\partial \tilde{u}}{\partial \tilde{x}_1}] + 
 \frac{\partial}{\partial \tilde{x}_1}[{\tilde{D}_{21}}\frac{\partial \tilde{u}}{\partial \tilde{x}_2}]
 -  \frac{\tilde{u}}{g}\frac{dg}{dt}, 
\end{split}
\quad\quad \tilde{\mathbf{x}}\in \tilde{\Omega}
\end{equation}
where  $\tilde{\Omega}= (-0.5,0.5) \times(-0.5,0.5)$
.
The dimensionless form of the boundary conditions in Eq  \eqref{eq:ADE} is:
\begin{equation}\label{eq:BC_dimensionless}
 \begin{array}{ll}
 \begin{split}
 \tilde{u}(\tilde{\mathbf{x}},\tilde{t}) & = \tilde{g}_d(\tilde{\mathbf{x}},\tilde{t})=g_d(\tilde{\mathbf{x}})/g(\tilde{t}), \quad \quad \tilde{x}_1 = -0.5 \\
\partial \tilde{u}(\tilde{\mathbf{x}},\tilde{t}) /  \partial \tilde{x}_1 & = \tilde{g}_n(\tilde{\mathbf{x}},\tilde{t}) = 
\frac{L_1}{g(\tilde{t})}
g_n(\tilde{\mathbf{x}}) , \quad \quad \tilde{x}_1 = 0.5 \\
\partial \tilde{u}(\tilde{\mathbf{x}},\tilde{t}) /  \partial \tilde{x}_2 & = \tilde{g}_n(\tilde{\mathbf{x}},\tilde{t})=
\frac{L_2}{g(\tilde{t})}{g}_n(\tilde{\mathbf{x}})  
, \quad \quad \tilde{x}_2 = -0.5 \quad{\rm{and}}\quad \tilde{x}_2 = 0.5
 \end{split}
 \end{array} 
 \end{equation}
 and the initial condition
 \begin{equation}\label{eq:ic_normalized}
 \tilde{u}(\tilde{\mathbf{x}},\tilde{t}=0) = 
\exp\left[ 
-\frac{(\tilde{x}_1-\tilde{x}_1^*)^2 }{2\tilde{\varepsilon}^2_{1}} 
-\frac{(\tilde{x}_2-\tilde{x}_2^*)^2 }{2\tilde{\varepsilon}^2_{2}} 
\right], \quad \tilde{\mathbf{x}}\in \tilde{\Omega}    
 \end{equation}
 where $\tilde{\varepsilon}^2_{1} = \varepsilon^2_{1}/L^2_{1}$ and  $\tilde{\varepsilon}^2_{2} = \varepsilon^2_{2}/L^2_{2}$.  
 The dimensionless velocity $\tilde{\mathbf{v}} = [\tilde{v}_{1}, \tilde{v}_{2}]^T$ can be found by solving the dimensionless form of Eq  
 (\ref{eq:GWF}):
\begin{equation}\label{eq:dimensionless_GWF}
 \begin{array}{ll}
 \begin{split}
 \frac{1}{L_1^2}\frac{\partial}{\partial \tilde{x}_1}[{\tilde{K}(\mathbf{x})}\frac{\partial \tilde{h}(\mathbf{x})}{\partial \tilde{x}_1}] + 
 \frac{1}{L_2^2}\frac{\partial}{\partial \tilde{x}_2}[{\tilde{K}(\mathbf{x})}\frac{\partial \tilde{h}(\mathbf{x})}{\partial \tilde{x}_2}] = 0,& \quad {\tilde{\mathbf{x}}} \in \tilde{\Omega} \\
 \tilde{h}(\mathbf{x}) = {\frac{H_2}{h_{max}}}, &\quad {\tilde{x}_1} = 0.5 \\
 -\tilde{K}(\mathbf{x}) \frac{\partial \tilde{h}(\mathbf{x})} {\partial \tilde{x}_1}  = \frac{qL_{x_1}}{k_{max}h_{max}}, &\quad \tilde{x}_1 = -0.5 \\
 -\tilde{K}(\mathbf{x}) \frac{\partial \tilde{h}(\mathbf{x})} {\partial \tilde{x}_2}  = 0, &\quad \tilde{x}_2 = -0.5 \:{\rm{or}}\: \tilde{x}_2 = 0.5.
 \end{split}
 \end{array}
 \end{equation}
where $\tilde{h}(\mathbf{x}) = \frac{h(\mathbf{x})}{h_{max}}$, $\tilde{K}(\mathbf{x}) = \frac{K(\mathbf{x})}{K_{max}}$. Here, $h_{max}$ and $K_{max}$ are the maximum head and hydraulic conductivity, respectively. $K_{max}$ is found from the known $K(\mathbf{x})$ or estimated from $K$ measurements if $K(\mathbf{x})$ is not fully known. For the considered boundary conditions, $h_{max}$ must be at the left boundary located at $x_1=0$. If  $h(\mathbf{x})$ is unknown, then $h_{max}$ can be approximated by replacing $K(\mathbf{x})$ in Eq  \eqref{eq:GWF} with its geometric average $\overline{K}$, in which case Eq  \eqref{eq:GWF} has an analytical solution and $h_{max}=\frac{qL_{x_1}}{\overline{K}} +H_2$.

 Then, $\tilde{\mathbf{v}}$ can be found from the dimensionless form of Darcy's law \eqref{eq:velocity}
\begin{equation}\label{eq:dimensionless_flux}
\begin{split}
\tilde{v}_1 = -\frac{\tilde{K}}{{\phi}} \frac{T K_{max}h_{max}}{L^2_{1}} \frac{\partial \tilde{h}}{\partial \tilde{x}_1}    \\
\tilde{v}_2 = -\frac{\tilde{K}}{{\phi}} \frac{T K_{max}h_{max}}{L^2_{2}} \frac{\partial \tilde{h}}{\partial \tilde{x}_2}
\end{split}
\end{equation}

Finally, we approximate  $\tilde{K}(\mathbf{x})$, $\tilde{h}(\mathbf{x})$, and  $\tilde{u}(\mathbf{x},t)$ with the DNNs $\hat{K}(\mathbf{x}; \alpha)$, $\hat{h}(\mathbf{x}; \beta)$, and $\hat{u}(\mathbf{x},t; \gamma)$, respectively, and train these DNNs as described in Section \ref{sec:PINN}. 
 
\subsection{Criteria for the weights in the PINN loss function} 

The loss function of PINN for the normalized ADE takes the form 
\begin{eqnarray}\label{eq:loss_norm_ADE}
 \tilde{L}_{\hat{u}} (\gamma) &=&  \frac{\lambda_{ic}}{N_{ic}} \sum_{n = 1}^{N_{ic}} \left[ \hat{u}({\mathbf{x}}_n, t = 0; \gamma) - \tilde{u}({\mathbf{x}}_n, t = 0)\right]^2  \\  \nonumber
&+& \frac{\lambda_{bcn}}{N_{bcn}^x \times N_{bcn}^t} \sum_{n = 1}^{N_{bcn}^x} \sum_{m = 1}^{N_{bcn}^t} \left[  
\mathcal{B}_{\tilde{u}}(\mathbf{x}_n,t_m; \gamma) 
 \right]^2  \\ \nonumber
&+& \frac{\lambda_{bcd}}{N_{bcd}^x \times N_{bcd}^t} \sum_{n = 1}^{N_{bcd}^x} \sum_{m = 1}^{N_{bcd}^t}  [\hat{u}(\mathbf{x}_n, t_m; \gamma) - \tilde{g}_d(\mathbf{x}_n, t_m)]^2 \\ \nonumber
&+& \frac{\lambda_{res}}{N_{res}} \sum_{n = 1}^{N_{res}} [\tilde{\mathcal{R}}_{\hat{u}}(\mathbf{x}_n, t_n; \gamma) ]^2, 
\end{eqnarray}
where  $\hat{u}({\mathbf{x}}_n, t; \gamma)$ is a DNN approximation of $\tilde{u}({\mathbf{x}}, t)$ and  $\mathcal{B}_{\tilde{u}}(\mathbf{x},t; \gamma) =  \nabla \hat{u}(\mathbf{x},t; \gamma)  \cdot \mathbf{n} - \tilde{g}_n(\mathbf{x}, t)$. 
All terms in this loss function are dimensionless, but they might be of different orders of magnitude. We propose to select weights such that they balance different terms in the loss function. 
 
The loss terms in Eq \eqref{eq:loss_norm_ADE} can be normalized as follows:

\begin{eqnarray}\label{eq:loss_IC}
\lambda_{ic} \tilde{L}_{ic}(\gamma) = \frac{\lambda_{ic} \tilde{u}^2_{max}}{N_{ic}} \sum_{n = 1}^{N_{ic}} \left[\frac{\hat{u}({\mathbf{x}}_n, t = 0; \gamma)} {\tilde{u}_{max}}
- \frac{\tilde{u}({\mathbf{x}}_n, t = 0)}{\tilde{u}_{max}} \right]^2,
\end{eqnarray}
\begin{equation}\label{eq:loss_res}
\lambda_{res} \tilde{L}_{res}(\gamma) =  
\frac{\lambda_{res} f_{max}^2}{N_{res}} 
\sum_{n = 1}^{N_{res}} \left [ \frac{ \tilde{\mathcal{L}}_{\hat{u}} (\mathbf{x}_n, t_n; \gamma) }  {f_{max}} - \frac{f (\mathbf{x}_n, t_n;\gamma )}{f_{max}}  \right]^2 ,
\end{equation}
\begin{equation}\label{eq:loss_bcn}
\lambda_{bcn} \tilde{L}_{bcn}(\gamma) = 
\frac{\lambda_{bcn}  \left( \frac{ \tilde{u}_{max}}  {\Delta x} \right)^2 }{N_{bcn}^x \times N_{bcn}^t} \sum_{n = 1}^{N_{bcn}^x} \sum_{m = 1}^{N_{bcn}^t}  
\left[ \mathcal{B}_{\tilde{u}}(\mathbf{x}_n,t_m; \gamma) \frac{\Delta x}{ \tilde{u}_{max}} \right]^2
 \end{equation}
and
\begin{equation}\label{eq:loss_bcd}
\lambda_{bcd} \tilde{L}_{bcd}(\gamma)  =
\frac{\lambda_{bcd} \tilde{u}^2_{max} }{N_{bcd}^x \times N_{bcd}^t} \sum_{n = 1}^{N_{bcd}^x} \sum_{m = 1}^{N_{bcd}^t}  
\left[
\frac{ \hat{u}(\mathbf{x}_n, t_m; \gamma) } {\tilde{u}_{max} } -
\frac {\tilde{g}_d(\mathbf{x}_n, t_m)} { \tilde{u}_{max} } 
\right]^2 ,
\end{equation}
where $\tilde{u}_{max} =  \max_{\mathbf{x},t}  \tilde{u} (\mathbf{x},t) = 1$,
  $\tilde{\mathcal{R}}_{\hat{u}} = \tilde{\mathcal{L}}_{\hat{u}} - f$, and
\begin{eqnarray}\label{eq:normalized_L}
\tilde{\mathcal{L}}_{\hat{u}} &=& \frac{\partial \hat{u}}{\partial \tilde{t}} + 
 \tilde{v}_{1} \frac{\partial \hat{u}}{\partial \tilde{x}_1} +
 \tilde{v}_{2} \frac{\partial \hat{u}}{\partial \tilde{x}_2}
 -
\frac{\partial}{\partial \tilde{x}_1}[{\tilde{D}_{11}}\frac{\partial \hat{u}}{\partial \tilde{x}_1}] -
\frac{\partial}{\partial \tilde{x}_2}[{\tilde{D}_{22}}\frac{\partial \hat{u}}{\partial \tilde{x}_2}] \\ \nonumber
&-& \frac{\partial}{\partial \tilde{x}_2}[{\tilde{D}_{12}}\frac{\partial \hat{u}}{\partial \tilde{x}_1}] - 
 \frac{\partial}{\partial \tilde{x}_1}[{\tilde{D}_{21}}\frac{\partial \hat{u}}{\partial \tilde{x}_2}],
\end{eqnarray} 
\begin{equation}\label{eq:normalized_g}
f =  - \frac{\hat{u}T}{g}\frac{dg}{dt}, 
\end{equation}
and 
\begin{equation}
f_{max} = \max_{\mathbf{x},t} f(\mathbf{x},t)  = \max_{\mathbf{x},t} \left( -  \frac{\hat{u}(\mathbf{x},t;\gamma)T}{g(t)}\frac{dg(t)}{dt} \right)=  \max_{t}   \left( - \frac{T}{g}\frac{dg}{dt} \right).
\label{eq:fmax}
\end{equation}

In Eqs \eqref{eq:loss_IC}, \eqref{eq:loss_res}, and \eqref{eq:loss_bcd}, the terms in parentheses are normalized with the maximum of $\tilde{u}$ or $f$. In Eq \eqref{eq:loss_bcn}, we obtain the normalization factor  $\frac{ \tilde{u}_{max}}  {\Delta x}$ by approximating $ \tilde{\mathcal{B}}_{\hat{u}}$ with a finite-difference scheme. For example, for the Neumann boundary at $\tilde{x}_1=-0.5$,
\begin{eqnarray}\label{para_bc}
 \tilde{\mathcal{B}}_{\hat{u}} = \frac{\partial \hat{u}}{\partial \tilde{x}} \approx \frac{\hat{u}(\tilde{x}+\Delta \tilde{x},\tilde{t}) -\hat{u}(\tilde{x},\tilde{t}) }{\Delta \tilde{x}} 
\end{eqnarray}
where $\Delta \tilde{x}$ is a grid size, which should be much smaller then the domain size $\tilde{L}<1$.

Finally, from Eqs \eqref{eq:loss_IC}--\eqref{eq:loss_bcd}, the weights are estimated from the equation  
\begin{equation}
\lambda_{ic} \tilde{u}^2_{max} = \lambda_{res} f_{max}^2  = \lambda_{bcn}  \left( \frac{ \tilde{u}_{max}}  {\Delta x} \right)^2 = \lambda_{bcd} \tilde{u}^2_{max}, 
\end{equation}
where, without loss of generality, we set $\lambda_{res} = 1$. Recalling the expression for $f_{max}$ in Eq \eqref{eq:fmax} and that $\tilde{u}_{max} = 1$, the estimates of the remaining weights are obtained as:

\begin{equation}
  \lambda_{ic} =  \left[ f_{max} \right]^2  =  \left[  \frac{T}{\varepsilon^2_x} (D_x + D_y)        \right]^2,
\end{equation}
\begin{equation}
\lambda_{bcn}  =  \Delta \tilde{x}^2 {\lambda_{ic}},
\end{equation}
and
\begin{equation}
\lambda_{bcd}  =   \lambda_{ic}.
\end{equation}

For the considered problem, $-  \frac{T}{g}\frac{dg}{dt}  $ varies by several orders of magnitude over time, and we find that a better estimate of  $\lambda_{ic}$ is obtained in terms of the mean of  $-  \frac{T}{g}\frac{dg}{dt}  $  rather than its maximum: 
\begin{equation}
  \lambda_{ic} =  \left[ \frac{1}{T}\int_0^T (-1)  \frac{T}{g}\frac{dg}{dt} dt \right]^2  =  \left[   \frac{1}{2}\ln{(2D_xT + \varepsilon^2)}+ \frac{1}{2}\ln{(2D_yT + \varepsilon^2)} - \ln{{\varepsilon}^2}       \right]^2. 
  \label{eq:lambdaic}
\end{equation}

\subsection{Adaptive residual points}
In the PINN method, residual points $(\mathbf{x_{res}}, t_{res})$ are commonly placed in the space-time domain using the LHS method.
We find that for the considered problem, using only the LHS method is inefficient. Therefore, we propose an adaptive approach for placing  residual points to reduce the PINN error for a given number of the residual points. 
As with quadrature points, it is natural to place more points where the gradients of the loss function are larger. For example, it was shown in \cite{mao2020high-speed} that clustering residual points in the vicinity of a shock wave at any given time can improve the accuracy of PINN in modeling high-speed aerodynamic flows. However, for the problem considered here where the amplitude of the solution $g(t)$ rapidly decays in time, we find that it is more important to place the residual points adaptively in time than space. 
We normalize $g(t)$ in Eq \eqref{eq:time-dependent normalizer} such that an integral of the normalized $g(t)$ function equals 1:
\begin{equation}\label{eq:adaptive sampling}
    \tilde{g}(t) = \frac{g(t)}{\int_0^T g(t) dt }
\end{equation}
Then, we treat $\tilde{g}(t)$ as a probability density function for randomly placing 50\% of the residual points in the temporal domain. The remaining residual points are placed in the time domain using the LHS method. The latter is to ensure that there are some residual points in the time domain with very small $\tilde{g}(t)$. The spatial locations of the residual points are determined using the LHS method.

\section{Numerical results: normalized PINN for the forward and inverse ADE}\label{sec:forward}
In this section, we evaluate the accuracy of PINN with adaptive residual point sampling for solving the forward and inverse normalized ADE \eqref{eq:normalized_ADE}. We note that both the normalization of ADE and adaptive residual point sampling are needed to obtain an accurate solution of the considered ADE problems. 
Our results (not shown here) demonstrate that employing the normalization or adaptive sampling (but not both) leads to $\varepsilon_{u}(t)$ errors exceeding 100\%.   

In the remainder of this section, we  study the effect of the weights $\lambda_{ic}$, $\lambda_{bc}$, and $\lambda_{res}$, the neural network size, and the number of residual points on the accuracy of the PINN method. Finally, we show how the PINN method can be used for data assimilation, a process where the state measurements are used for solving the ADE with partially known $K(\mathbf{x})$ and $h(\mathbf{x})$.

In all simulations, we set $N_{res} = 100,000$ and choose residual points adaptively in time.
The space and time coordinates $(\mathbf{x}_n, t_m)$  in the initial and boundary condition terms are selected as in  Eq \eqref{eq:u net}. 
 The neural network structure has five hidden layers and 60 neurons per layer.

\subsection{Effect of weights in the loss function}\label{sec:weights}

 Here, we study the effect of the weights on the PINN accuracy. We set $\lambda_{res} = 1$. Then, for the ADE parameters used in this study ($\varepsilon = 0.025, D_x = 0.0929, D_y = 0.0645$), Eq \eqref{eq:lambdaic} yields $\lambda_{ic} = 4.83^2$, $\lambda_{bcd}=\lambda_{ic}$, and $\lambda_{bcn} = \kappa \lambda_{ic}$, where $\kappa <<1$.
Below, we compare the PINN solutions obtained with these values of weights and weights selected according to other methods.  

\begin{figure}[h]
	\centering
	\subfloat[] {\includegraphics[angle=0,width=0.45\textwidth]{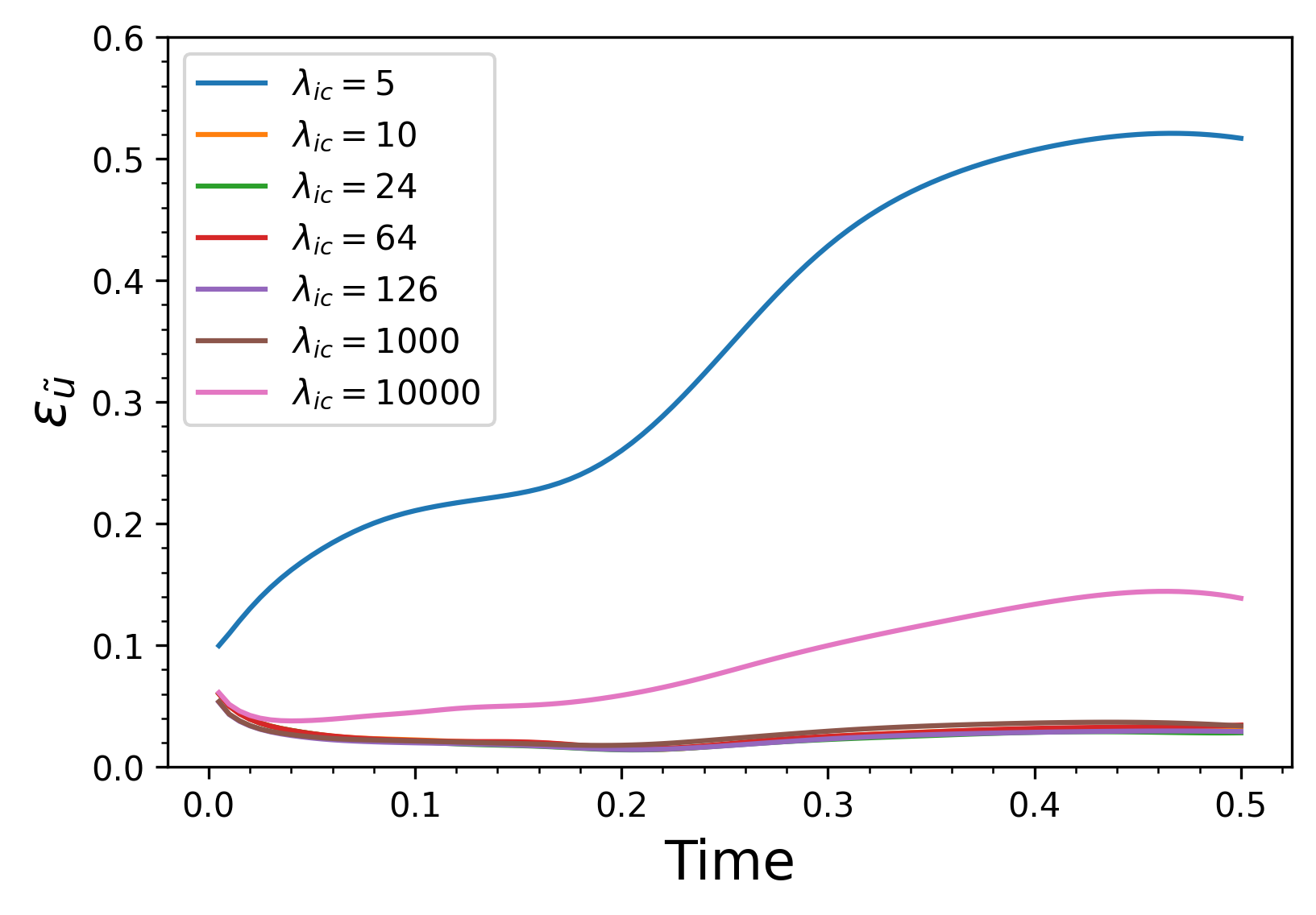}}
	\subfloat[] {\includegraphics[angle=0,width=0.45\textwidth]{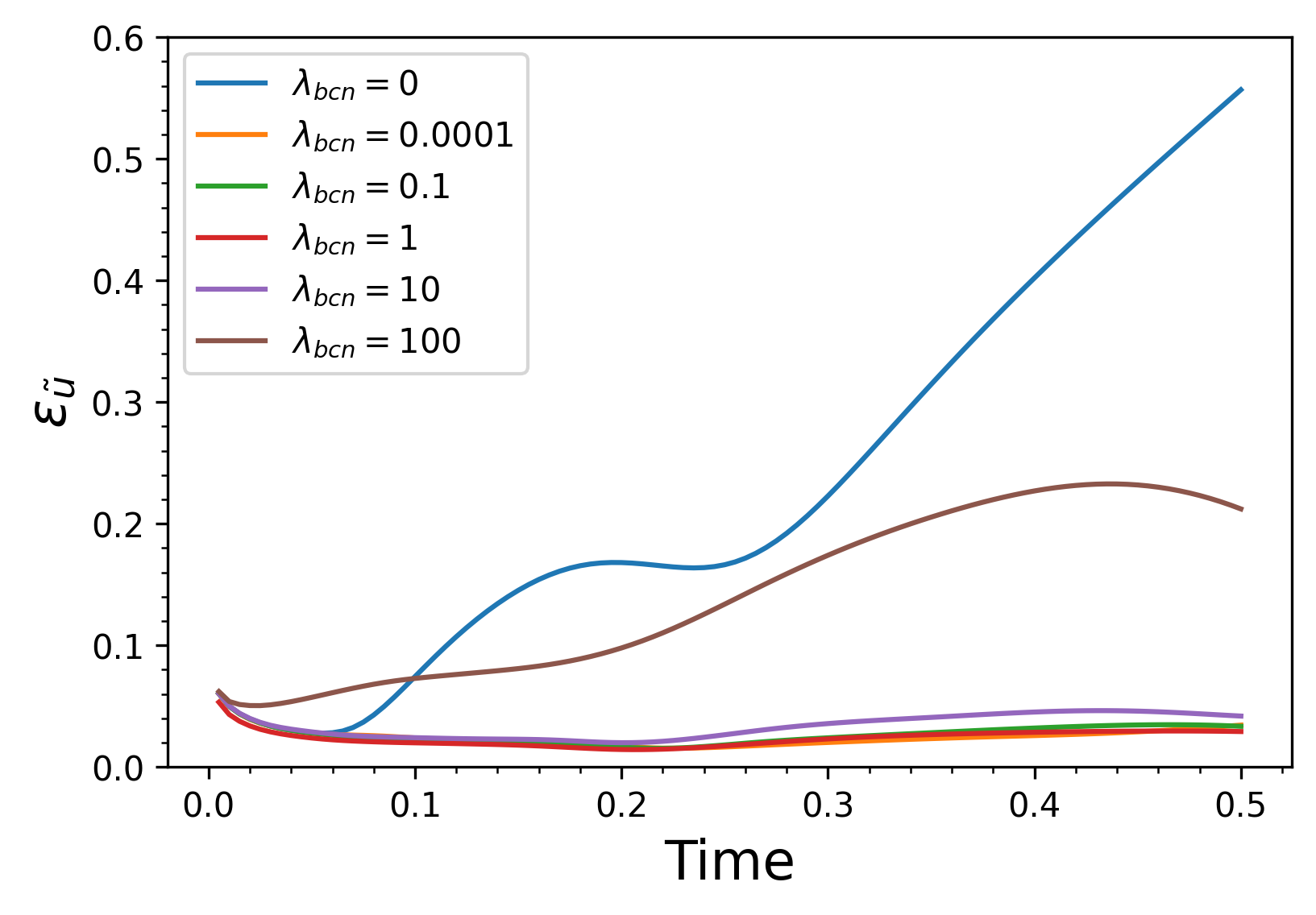}}
	\caption{(a) Relative $L_2$ errors $\varepsilon_u$ in the PINN solutions as functions of time with for different $\lambda_{ic}=\lambda_{bd}$ values. The other weights are $\lambda_{res}=\lambda_{bcn}=1$. 
	(b) Relative $L_2$ errors $\varepsilon_u$ in the PINN solutions as functions of time for different values of $\lambda_{bcn}$. The other weights are $\lambda_{res} = 1$ and $\lambda_{ic} = \lambda_{bd}=126$.
	}
	\label{fig:error_weight}
\end{figure}
Figure \ref{fig:error_weight}(a) shows the relative $L_2$ errors $\varepsilon_u(t)$ for $\lambda_{bcn}=1$ and  values of $\lambda_{ic} = \lambda_{bcd}$ in the range [5,10000]. In this figure, the largest error (exceeding 50\%) corresponds to $\lambda_{ic} = 5$. Not shown here,  the error exceeds 100\% for $\lambda_{ic} = 1$. 
For $\lambda_{ic}=10000$, the relative error is about 12\%.  
However, for $\lambda_{ic}$ in the range $[10, 1000]$, the PINN errors do not exceed 6\% and are practically independent of $\lambda_{ic}$ (vary by less than 1\%).  We find that the smallest error is achieved at  $\lambda_{ic}=126$. We plot the concentration field and the absolute point errors at selected time snapshots for this case in Fig \ref{fig:PINN_forward_normalized}.

The value estimated by the proposed criterion in Eq \eqref{eq:normalized_g} is $\lambda_{ic}= 24$, which is within the range of $\lambda_{ic}$ that produces an accurate solution. The error obtained with $\lambda_{ic}= 24$ is also shown in Figure \ref{fig:error_weight}(a) and is very close to the smallest error obtained with $\lambda_{ic}=126$. In the following examples, we set $\lambda_{ic}=126$. 

Figure \ref{fig:error_weight}(b) shows the relative $L_2$ errors $\varepsilon_u(t)$ for $\lambda_{bcn}$ varying in the range [0,100]. Other weights are set to $\lambda_{ic} =\lambda_{bcd} =126$ and $\lambda_{res} = 1$.  For $\lambda_{bcn} = 0$ and 100, the maximum relative errors are 56\% and 24\%, respectively.  The smallest errors are for $\lambda_{bcn}$ in the [0.0001,1] range. In this range, the errors do not exceed 6\% and are practically independent of $\lambda_{bcn}$.

These results agree with the $\lambda_{bcn}$ estimate $\lambda_{bcn} = \kappa \lambda_{ic}$ ($\kappa <<1$). 

\subsection{Effect of the neural network size and the number of residual points}\label{sec:size_residual}
\begin{figure}[h]
	\centering
	\subfloat[] {\includegraphics[angle=0,width=0.45\textwidth]{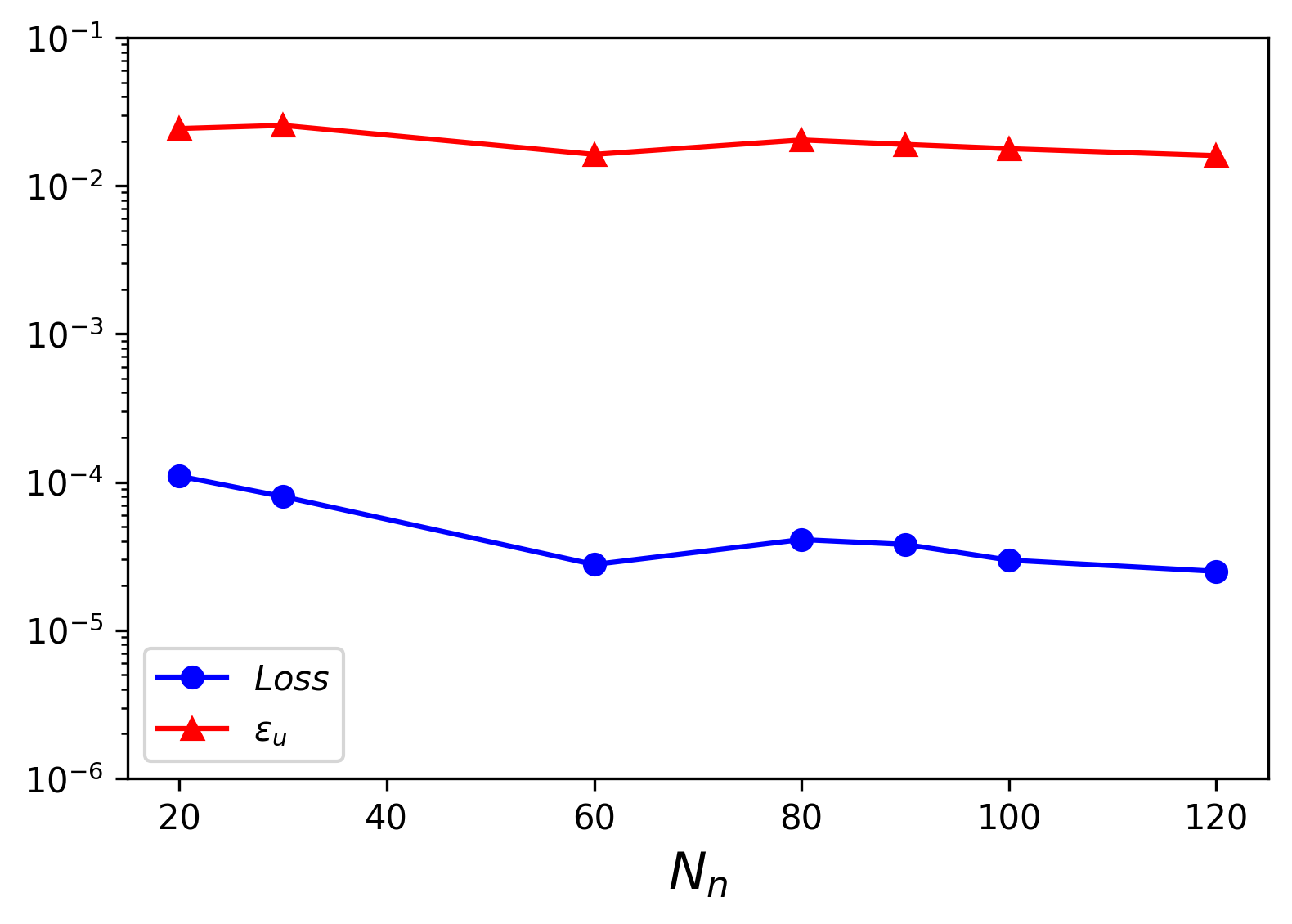}}
	\subfloat[] {\includegraphics[angle=0,width=0.45\textwidth]{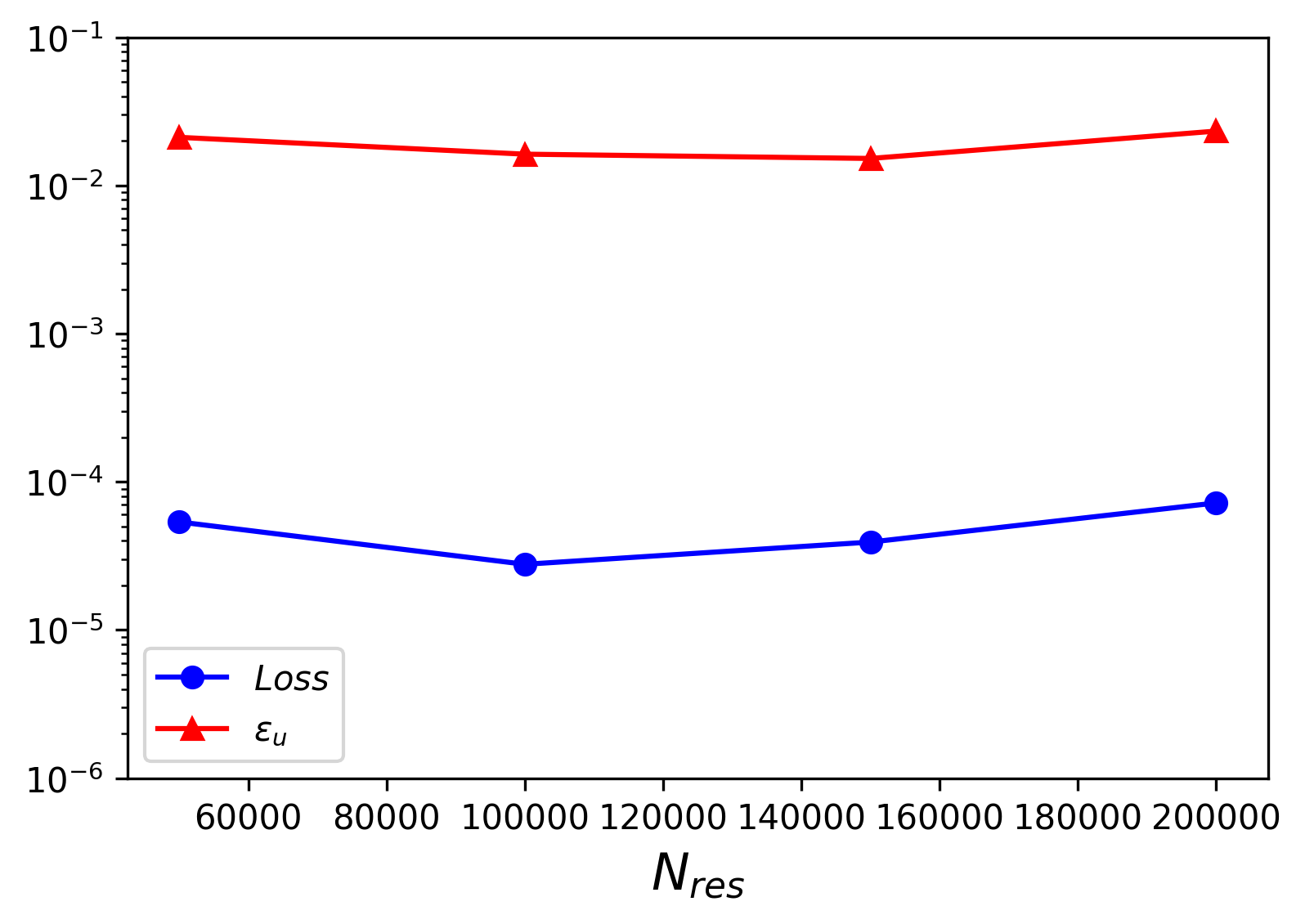}}
	\caption{(a) The loss function at the end of the training and the average-over-time relative L2 error as functions of $N_n$. $N_{res}=100,000$. (b) The  loss at the end of the training and the average-over-time relative L2 error as functions of $N_{res}$. $N_n =60$.}
	\label{fig:average_rl2_error_size_neuron}
\end{figure}

In this section, the dependence of PINN accuracy on the neural network size and the number of residual points is investigated. In all simulations, the weights are set to 
$\lambda_{ic} =\lambda_{bcd} =126$ and $\lambda_{res} = \lambda_{bcn} =1$. 
First, we study the PINN error as the function of the number of neurons per hidden layer $N_n$. The number of hidden layers is kept at 5, and 100,000 residual points are distributed adaptively in time. To facilitate the comparison between different cases, we show the averaged-over-time relative L2 error
\begin{equation}\label{eq:rl2_err_over_time}
\varepsilon_{\tilde{u}} = \sqrt{\frac{\sum_{i=1}^{N_t}\sum_{j=1}^{N_{tot}}\left[\tilde{u}(\mathbf{x}_j,t_i)-\hat{u}(\mathbf{x}_j,t_i)\right]^{2}}{\sum_{i=1}^{N_t}\sum_{j=1}^{N_{tot}}\tilde{u}(\mathbf{x}_j,t_i)^{2}}},
\end{equation}
where $N_t$ is the number of time snapshots in the simulation. Figure \ref{fig:average_rl2_error_size_neuron}(a) shows $\varepsilon_{\tilde{u}}$ and the magnitude of  the PINN loss as functions of $N_n$. For $N_n$ ranging from 20 to 120,  $\varepsilon_{\tilde{u}}$ and the loss are decreasing with  increasing $N_n$. However, the dependence on $N_n$ is small. We find that $\varepsilon_{\tilde{u}}$ is almost the same for $N_n=60$ and 120. 
The comparison of this figure with Figure  \ref{fig:error_original_PINN} reveals the importance of the proposed normalization of the ADE because the DNN size of five layers and $N_n=60$ is sufficient to produce an accurate PINN solution of the normalized ADE. However, this leads to a highly inaccurate solution when applied to a non-normalized ADE. 

Next, we study the effect of the number of residual points on the PINN error. Here, we fix the size of DNN at five hidden layers with $N_n= 60$ and consider the number of residual points in the range [50000, 200000]. 
The $\varepsilon_{\tilde{u}}$ error and PINN loss as functions of $N_{res}$ are shown in Figure \ref{fig:average_rl2_error_size_neuron}(b). For $N_{res}$ in the considered range, the error is practically independent of $N_{res}$ with the smallest errors obtained for $N_{res}=100,000$ and 160,000. The dependence of the loss on $N_{res}$ is weak with the smallest values of the loss terms obtained for $N_{res}=100,000$. A small increase in the value of the loss with increasing $N_{res}$ is expected because the DNN approximation of the solution is not exact.   

\subsection{Parameter and state estimation}
In Sections \ref{sec:weights} and \ref{sec:size_residual}, the PINN solutions were obtained without any measurements of $u$ except for the initial and boundary conditions. However, the parameter and state measurements can sometimes be available and can easily be included or ``assimilated'' in the PINN formulation of the PDE solution.  
In the following, we show that data assimilation can help accurately estimate $u$ even when the velocity field $\mathbf{v}$ is not fully known. Here, we assume that $N_K$, $N_h$, and $N_u$ measurements of the $K(\vec{x})$, $h(\vec{x})$, and  $u(\vec{x},t)$ fields are available. 
 Because of the incomplete knowledge of the $K$ field, $\mathbf{v}$ cannot be exactly found from the Darcy equation. 
However, the solutions of the Darcy and ADE equations can be found by formulating the joint PINN data assimilation problem as 
\begin{equation}\label{eq:loss_inverse}
\alpha^*, \beta^*, \gamma^* = \min_{\alpha,\beta,\gamma} \Big\{\tilde{L}(\alpha, \beta, \gamma) = \tilde{L}_{Da}(\alpha, \beta) + \tilde{L}_{AD}(\alpha, \beta, \gamma) + \tilde{L}_{data}(\alpha, \beta, \gamma) \Big\}
\end{equation}
where $\alpha$, $\beta$, and $\gamma$ are parameters of the $\hat{k}$, $\hat{h}$, and $\hat{u}$ DNNs, 
\begin{eqnarray}\label{eq:loss_Darcy_assimilation}
\tilde{L}_{Da}(\alpha, \beta)  &=& \frac{\lambda_{bcn}^{Da}}{N_{bcn}^{Da}} \sum_{n = 1}^{N_{bcn}^{Da}}  \left[ \frac{\partial \hat{h}(\mathbf{x}; \beta)}{\partial \mathbf{x}} \Big \rvert_{\mathbf{x} = \mathbf{x}_n} \cdot \mathbf{n} - \tilde{q}(\mathbf{x}_n) \right]^2  \\ \nonumber
&+& \frac{\lambda_{bcd}^{Da}}{N_{bcd}^{Da}} \sum_{n = 1}^{N_{bcd}^{Da}}  [\hat{h}(\mathbf{x}_n; \beta) - \tilde{h}(\mathbf{x}_n)]^2 + \frac{\lambda_{res}^{Da}}{N_{res}^{Da}} \sum_{n = 1}^{N_{res}^{Da}} [\tilde{\mathcal{R}}_{\hat{k},\hat{h}}(\mathbf{x}_n; \alpha, \beta) ]^2
\end{eqnarray}
is the loss term enforcing constraints due to the Darcy flow equation, and  
\begin{eqnarray}\label{eq:loss_data_assimilation}
\tilde{L}_{data}(\alpha, \beta, \gamma) &=& \frac{\lambda_{data}}{N_{x}N_t} \sum_{n = 1}^{N_{x}} \sum_{m = 1}^{N_{t}}  (\hat{u}(\mathbf{x}_n^u,t_m^u; \gamma) - u^*(\mathbf{x}_n^u, t_m))^2 \\
&+& \frac{\lambda_{data}}{N_{K}} \sum_{n = 1}^{N_{K}} (\hat{K}(\mathbf{x}_n^k,; \alpha) - K^*(\mathbf{x}_n^k))^2  
+\frac{\lambda_{data}}{N_{h}} \sum_{n = 1}^{N_{h}} (\hat{h}(\mathbf{x}_n^h; \beta) - h^*(\mathbf{x}_n^h))^2 \nonumber
\end{eqnarray}
is the loss term enforcing the measurements of $u$, $K$, and $h$ on the PDE solutions. The form of $\tilde{L}_{AD}(\alpha, \beta, \gamma)$ is identical to Eq \eqref{eq:loss_norm_ADE} except that $\alpha$ and $ \beta$ are now unknown parameters. 

In Eq \eqref{eq:loss_Darcy_assimilation}, $\tilde{\mathcal{R}}_{\hat{k},\hat{h}}(\mathbf{x}_n; \alpha, \beta) = \tilde{\nabla}\cdot (\hat{k}\tilde{\nabla} \hat{h})$ is the residual of the normalized Darcy equation \eqref{eq:dimensionless_GWF}.  The Dirichlet and Neumann BCs of the Darcy flow equation are enforced at the $N^{Da}_{bcd} = N^x_{bcn}$ and $N^{Da}_{bcn} = N^x_{bcd}$ number of points, respectively.

In Eq \eqref{eq:loss_data_assimilation}, $u^*(\mathbf{x}_n^u, t_m)$ is a normalized measurement of $u$ at the location $\mathbf{x}_n^u$ ($n=1,...,N_x$) and time  $t_m$ ($m=1,...,N_t$).  We assume that  $u^*$ is sampled uniformly in the time interval $(0,T]$ and that the number of $u^*$ measurements at each spatial location is $N_t = 100$, resulting in the total number of $u$ samples $N_u = N_x \times N_t$.  The sampling locations for $u^*$ are selected randomly in space.  Similarly,  $K^*(\mathbf{x}_n^k)$ and $h^*(\mathbf{x}_n^h)$ are the normalized measurements of $K$ and $h$ at the locations $\mathbf{x}_n^K$ and $\mathbf{x}_n^h$, respectively. 
In the following examples, we set $N_k=N_h=20$, and vary $N_{x}$ in the range [0,60]. The number of residual points in the Darcy and ADE residual losses are $N^{Da}_{res} = 10,000$ and $N^{AD}_{res} = 60,000$, respectively. The residual points for the Darcy equation are selected using the LHS scheme. The points for enforcing the initial and boundary conditions and the residuals in $\tilde{L}_{AD}$ are selected as in the forward ADE problem. The magnitude of the $\tilde{L}_{data}$ term is the same as that of the  $\tilde{L}_{ic}$ term in Eq \eqref{eq:loss_IC}. Therefore, we set 
 set $\lambda_{data}=\lambda_{ic} $. Following \cite{tartakovsky2020physics}, the weights in the loss term in Eq \eqref{eq:loss_Darcy_assimilation} are selected as $\lambda_{bcn}^{Da} = \lambda_{bcd}^{Da} =\lambda_{res}^{Da} =1$. 
Finally, the  $\hat{k}$, $\hat{h}$, and $\hat{u}$ DNNs are set to five hidden layers with 60 neurons per layer. 

\begin{figure}[h]
	\centering
	\subfloat[] {\includegraphics[angle=0,width=0.45\textwidth]{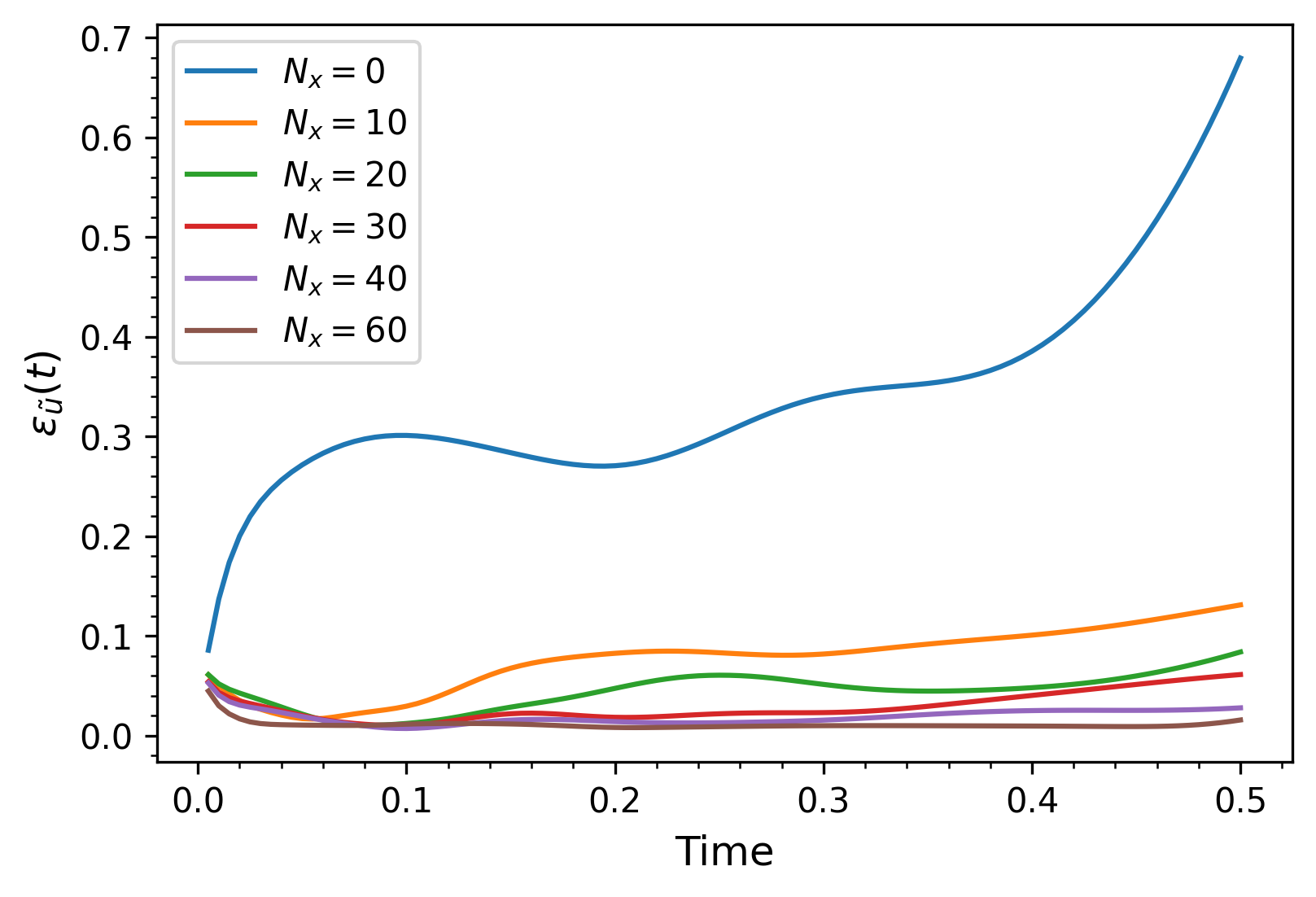}}
	\subfloat[] {\includegraphics[angle=0,width=0.45\textwidth]{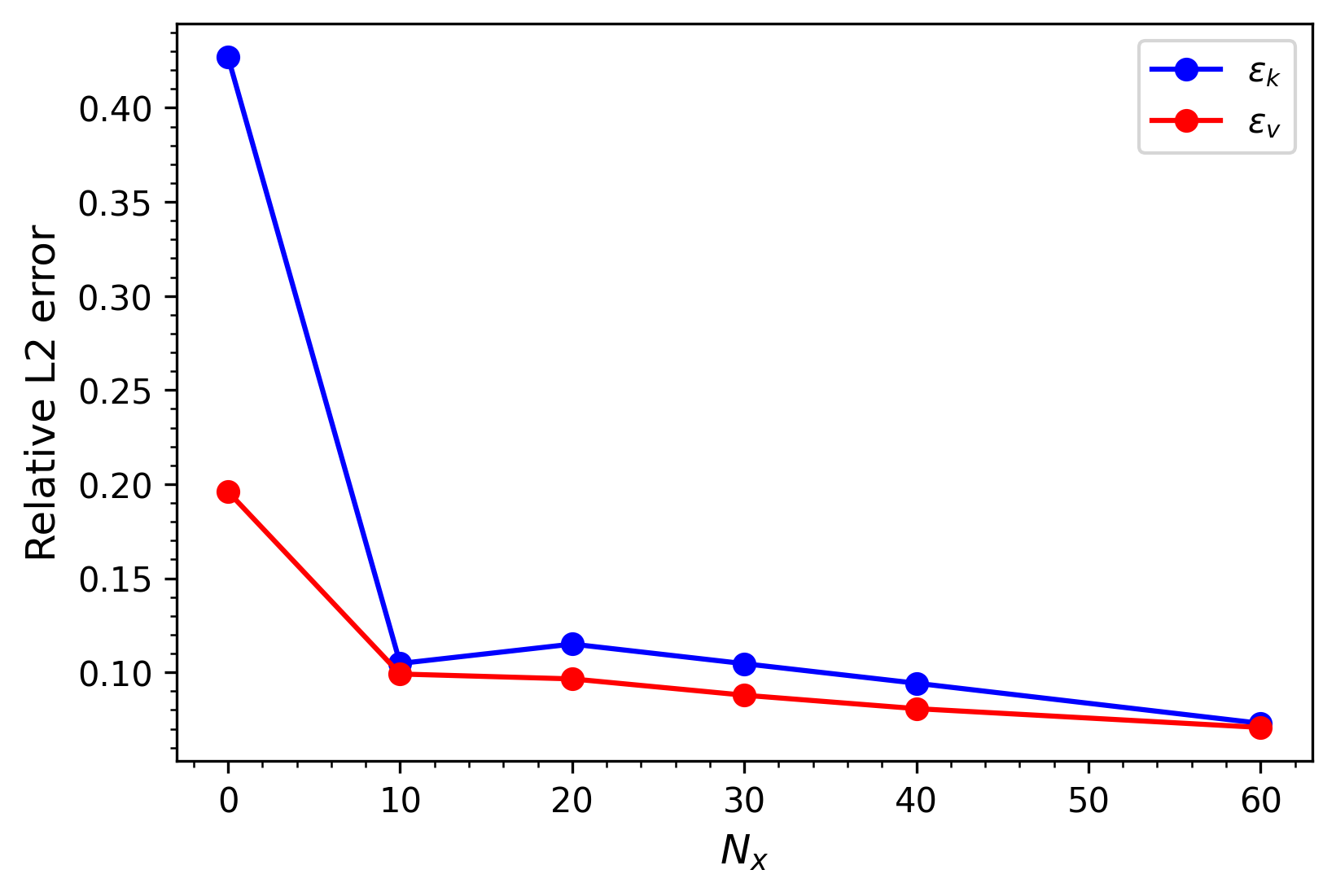}}
	\caption{(a) Relative $L_2$ errors in the PINN $\hat{u}$ solutions versus time for different $N_x$, the number of concentration measurements.  $N_k= N_h = 20$, $N^{Da}_{res}=10,000$, and $N^{AD}_{res}=60,000$. (b) Relative $L_2$ errors in the PINN estimates of conductivity ($\varepsilon_v$) and the magnitude of velocity ($\varepsilon_v$) as functions of $N_x$. 
	}
	\label{fig:error_PINN_inverse}
\end{figure}

Figure \ref{fig:error_PINN_inverse}(a) presents the relative PINN error in the estimated $u$ field as a function of time for different values of $N_x$. It can be seen that the error decreases with increasing $N_x$, with the the maximum error close to 70\% for $N_x=0$ and less than 5\% for $N_x=60$.  
Figure \ref{fig:error_PINN_inverse}(b) presents the relative L2 errors in the estimated hydraulic conductivity and the velocity magnitude as functions of $N_x$. 
The error in the estimated $K$ is reduced from over 43\% for $N_x=0$ to less than 8\% for $N_x=60$.
The error in the estimated velocity decreases from 20\% for $N_x=0$ to less that 8\% for $N_x=60$. 
 We note that the error in the estimated $u$ field  obtained with $N_x=60$ $u$ measurements is similar to the error in the estimated $u$ with the fully known velocity field in Section \ref{sec:size_residual}. These results show the ability of the PINN method to assimilate static and dynamic data for the state and parameter estimation. 

 \section{PINN for backward ADE}\label{sec:backward_ADE}
 \subsection{DNN with identity activation function in the last layer}\label{sec:identity_function}
\begin{figure}[h]
	\centering
	\includegraphics[angle=0,width=0.45\textwidth]{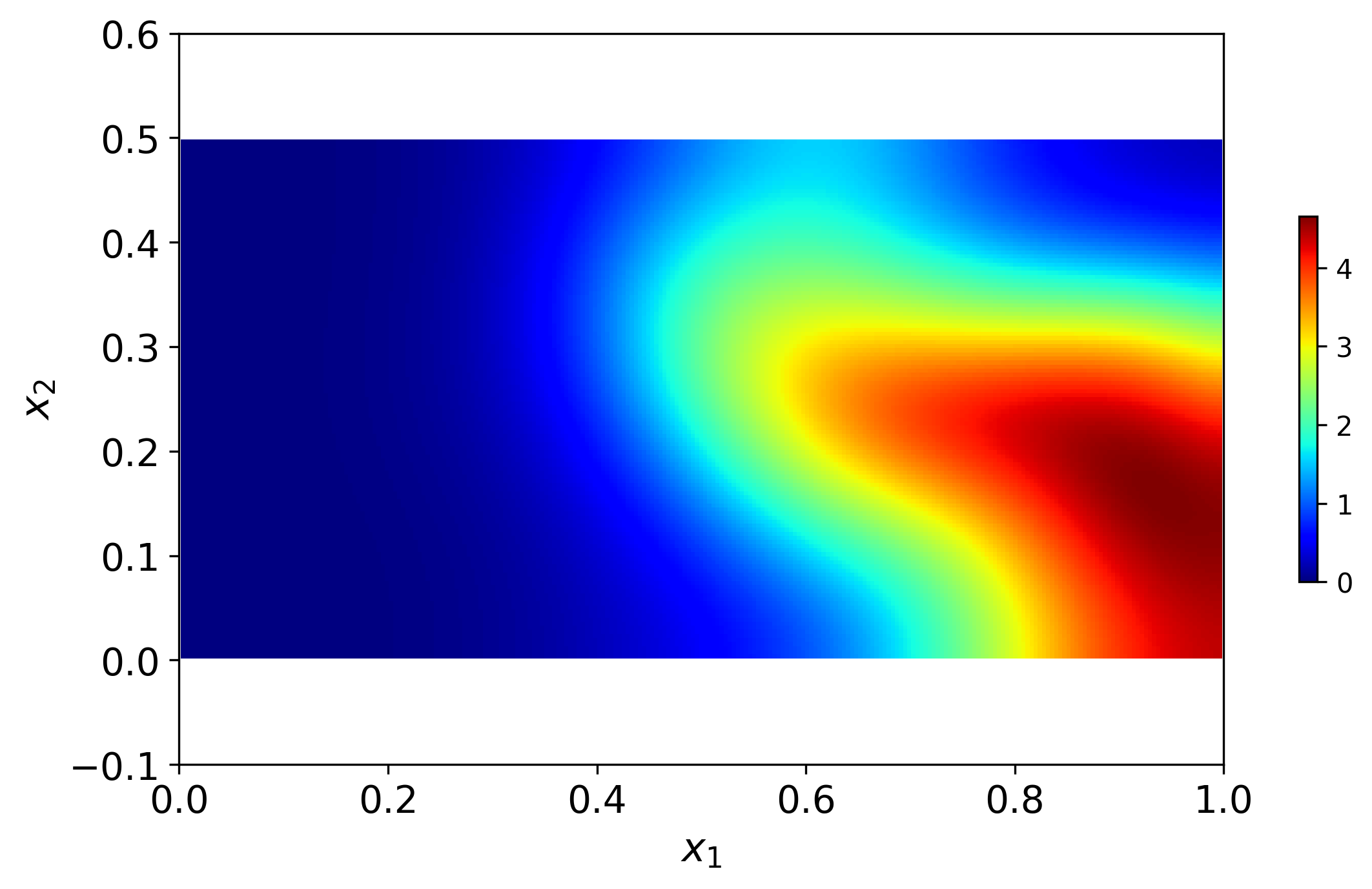}
	\caption{Terminal condition for the backward ADE at the terminal time $T = 0.2$.}
	\label{fig:terminal condition}
\end{figure}

 Here, we investigate the accuracy of PINN for solving backward ADEs,  which can arise in the source identification problems. We consider the backward ADE, which before normalization is given by Eq \eqref{eq:ADE} subject to the terminal condition \eqref{eq:terminal_cond}. The terminal time $T$ is set to 0.2. At $T=0.2$, approximately 70\% of the released mass remains in the computational domain $\Omega$. 
 The concentration field $U(\mathbf{x})$ at the terminal time is given by by the MT3DMS solution of the forward ADE at $t=T$ and is shown in Figure \ref{fig:terminal condition}.

 As for the forward problem, we normalize the backward ADE by introducing the dimensionless variable 
 $\tilde{u}(\mathbf{x},t) = \frac{u(\mathbf{x},t)}{g(t)}$ where $g(t)$ is defined by Eq
\eqref{eq:time-dependent normalizer}.
 This leads to the dimensionless form of the backward ADE  \eqref{eq:BC_dimensionless} subject to the dimensionless boundary conditions \eqref{eq:ic_normalized} and the terminal condition
 \begin{equation}\label{eq:terminal_cond_dim}
\tilde{u}( \tilde{\mathbf{x}},\tilde{t}=\tilde{T}) = \tilde{U}( \tilde{\mathbf{x}}).
\end{equation}
The proposed normalization assumes that the release time ($t=0$) and the released mass ($M=1$) are known. Also, this normalization assumes that the plume at the time of release has the Gaussian shape with the known standard deviation $\varepsilon$. The location of the plume at the time of release is treated as an unknown.

 The PINN solution of the dimensionless backward PDE is found by minimizing the loss function:
 \begin{equation}\label{eq:loss_backward_ADE}
     \tilde{L}_{\hat{u}}(\gamma) = \lambda_{tc}\times\tilde{L}_{tc}(\gamma) + \lambda_{bcn}\times\tilde{L}_{bcn}(\gamma) + \lambda_{bcd}\times\tilde{L}_{bcd}(\gamma) + \lambda_{res}\times\tilde{L}_{res}(\gamma)
 \end{equation}
Here, the weights $\tilde{L}_{bcn}$, $\tilde{L}_{bcd}$, and $\tilde{L}_{res}$ are the same as in the loss function for the normalized forward ADE.  The term $\tilde{L}_{tc}$ enforces the terminal condition \eqref{eq:terminal_cond_dim} and has the form:
  \begin{equation}
     \tilde{L}_{tc} = \frac{1}{N_{tc}} \sum_{n = 1}^{N_{tc}} (\hat{u}(\mathbf{x},t = T; \gamma) - \tilde{U}(\mathbf{x}))^2.
 \end{equation}
The analysis for selecting weights in the loss function \eqref{eq:loss_backward_ADE} is similar to one for the forward ADE in Section \ref{sec:weights}, yielding the same values for $\lambda_{bcn}$,  $\lambda_{bcd}$, and $\lambda_{res}$. The weight $\lambda_{tc}$ in the backward ADE equals to $\lambda_{ic}$ in the forward ADE. 

The terminal condition is enforced using 75\% of the terminal concentration data (given on the finite difference mesh) with the randomly selected measurements locations. 
To enforce the boundary conditions, we set $N_{bcn}^{t} = N_{bcd}^{t} = 40$. The spatial locations of the points on the boundaries are selected as in the forward problem. The number of residual points is $N_{res}=100,000$, and the locations of these points are chosen with the adaptive-in-time algorithm.  

We first test the same DNN architecture to approximate $\tilde{u}$ as in the forward ADE problem that has five hidden layers with 60 neurons per layer and an identity activation function in the last layer.
Snapshots of the PINN solutions are shown in Figure \ref{fig:PINN_backward_identity_0.2}. It can be seen that the solution develops negative concentration values and the magnitude of the negative values increases as time decreases. 

In the source identification problems, the important quantities of interest are the trajectory of the center of mass, the initial source location, and the  mass of the tracer as a function of time. The center of mass $(X_1^c(t), X_2^c(t))$ of the plume can be found as:
 \begin{equation}
 \begin{split}
     X_1^c(t)  = \frac{\int_{0}^{L_{x_2}}\int_{0}^{L_{x_1}}x_1\hat{u}(x_1, x_2, t)dx_1dx_2}{M(t)}\\
     X_2^c(t)  = \frac{\int_{0}^{L_{x_2}}\int_{0}^{L_{x_1}}x_2\hat{u}(x_1, x_2, t)dx_1dx_2}{M(t)},\\
\end{split}
\end{equation}
where $M(t) = \int_{0}^{L_{x_2}}\int_{0}^{L_{x_1}}\hat{u}(x_1, x_2, t)dx_1dx_2$ is the mass of the tracer time $t$.

 Figures \ref{fig:COM_PINN_backward}(a) and (b) display the position of the center of mass and the tracer mass as functions of time computed from the PINN and reference solutions. It can be seen that the PINN solution incorrectly estimates the source location at time $t=0$ and the tracer mass. 
In the following section, we propose a modification to the DNN network structure to improve the accuracy of PINN for backward problems by enforcing the positivity of the PINN solution.  

\subsection{Modified DNN with a positive wrapper function in the last layer}\label{sec:wrapper_function}
To impose the positivity of the PINN solution, we add a positive ``wrapper'' function $\sigma_w(z)$ to the  output node of the $\hat{u}$ DNN. The output of the modified DNN thus becomes $\hat{u}(\mathbf{x},t,\gamma) = \sigma_w (\sigma_{N_h +1} (\mathbf{x},t,\gamma))$. The added wrapper function maps the output of the identity activation function $\sigma_{N_h +1} (\mathbf{x},t,\gamma)$ to $\hat{u} \in [0,1]$.  
We use the sigmoid function as the wrapper function, $\sigma_w(z) =  \frac{1}{1+e^{-z}}$. Figure \ref{fig:PINN_backward_sigmoid_0.2} shows snapshots of the PINN solution with the sigmoid wrapper function and the point errors. It can be seen that the PINN solution is strictly positive. Also, the maximum point errors in the PINN solution with the modified structure are smaller than in the PINN solution with the original DNN structure. For example, at time t=0.1, the maximum errors in the PINN solution with the modified and original DNNs are 0.3 and 0.7, respectively. The errors are largest at time $t=0$, with the maximum errors for the modified and original DNNs being 120 and 200, respectively. As expected, errors in the PINN solution of the backward ADE are largest at $t=0$. 

\begin{figure}[h]
	\centering
	\subfloat[] {\includegraphics[angle=0,width=0.5\textwidth]{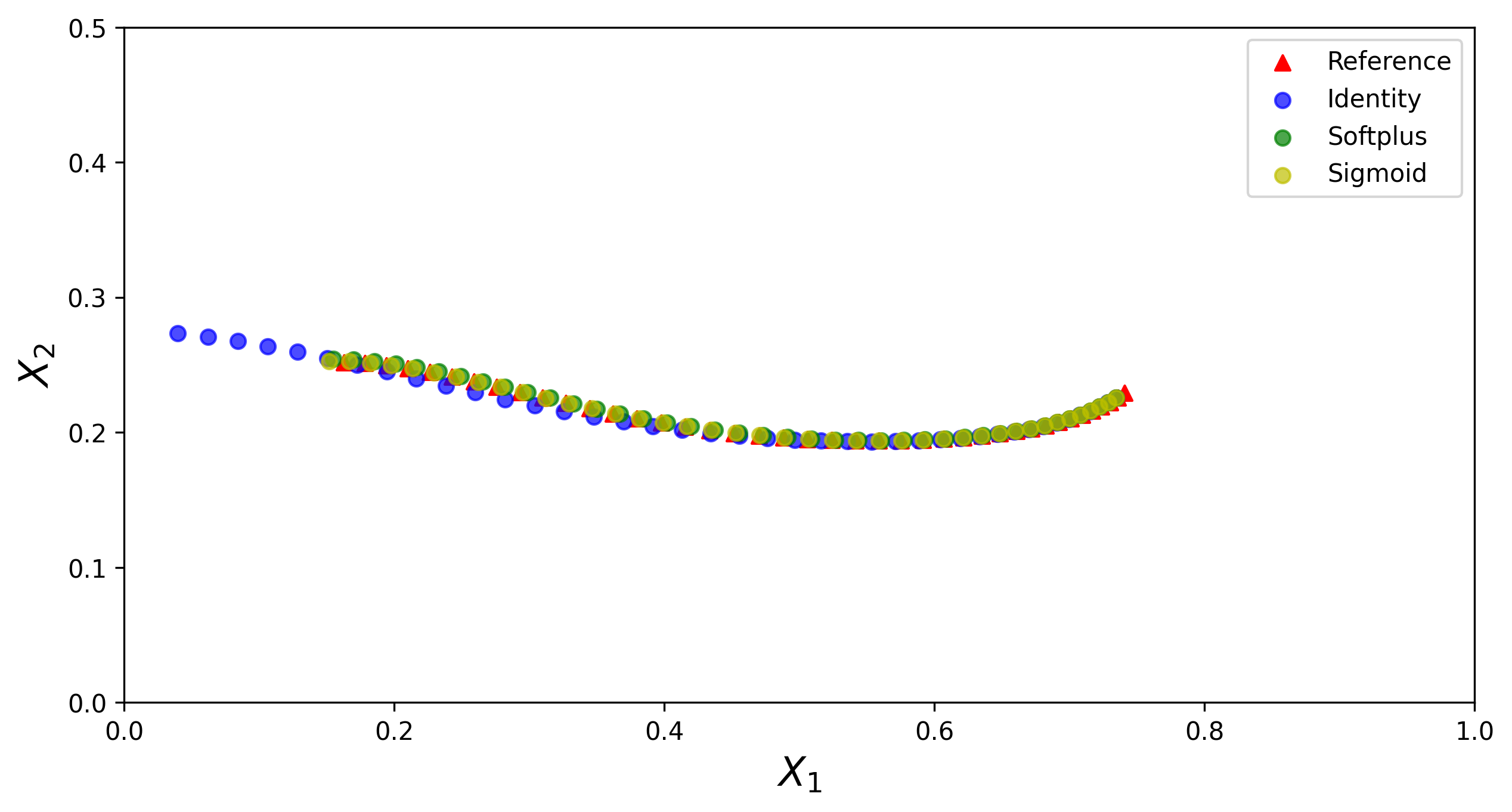}}
	\subfloat[] {\includegraphics[angle=0,width=0.5\textwidth]{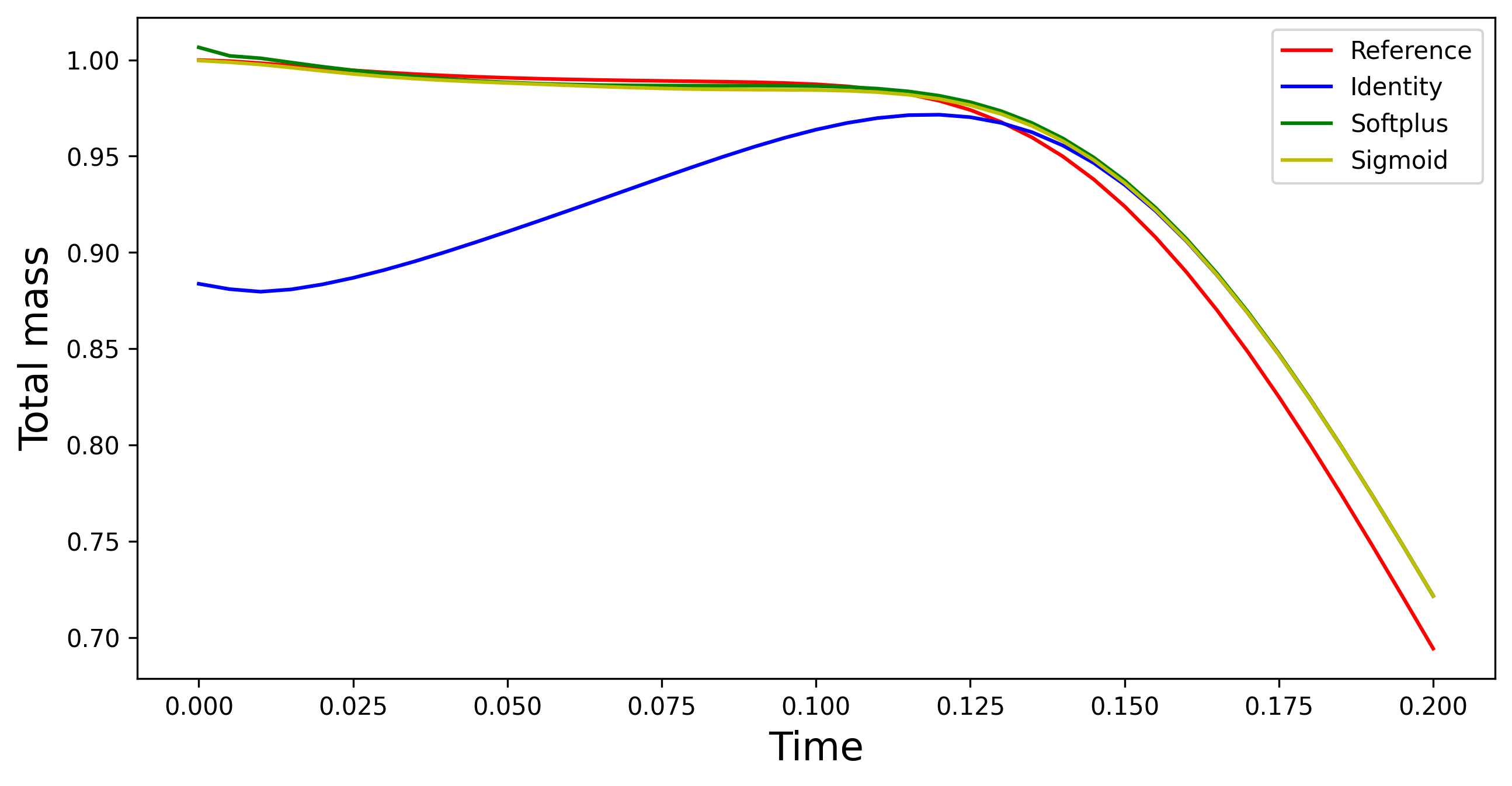}}
	\caption{(a) The location of the center of mass obtained from the reference solution and PINN solution with three different activation functions in the output layer. (b) Total mass in the domain computed from the reference solution and the PINN solutions with different activation functions in the last layer. No measurements of $u$ are used to obtain the PINN solution $\hat{u}$.}
	\label{fig:COM_PINN_backward}
\end{figure}

Figure \ref{fig:COM_PINN_backward}(a) compares the positions of the center of mass as functions of time obtained from the PINN solutions with the identity and sigmoid activation functions in the last layer as well as from the reference solution. The mass of the tracer, obtained from these solutions, is shown in Figure \ref{fig:COM_PINN_backward}(b). Positions of the center of mass predicted with the sigmoid  function are in significantly better agreement with the reference values than those predicted with the identity function. 
Also, the tracer mass in the PINN solution with the sigmoid wrapper function is in good agreement with the reference values, which is rather remarkable given that 30\% of the total mass had left the domain before the terminal time $T$. 

We note that while we assume that the total mass $M$ at the time of release and the release time are known, this information is only used to normalize the solution so its maximum at all times is approximately 1. Therefore, any reasonable guess about the time of release and $M$ would produce similarly accurate predictions. 
We also note that there is no benefit  of using the positive wrapper function for solving the forward ADE. We find that adding the wrapper function increases the error in the forward solution--for example, at late times the L2 error in the PINN solution for $u$ increases from 2\% (without the wrapper function) to 8\% (with the wrapper function).   As stated in Section \ref{sec:DNN_approximation}, the PINN method with the tanh activation functions (varying from -1 to 1) in the hidden layers   has  better convergence properties than that with  positive activation functions such as a sigmoid function. Adding the sigmoid activation function in the last layer might have a similar effect on the convergence as having it in all hidden layers. 

 In the considered problem, the normalized solution is in the [0,1] range. Therefore, the sigmoid function can be used as a wrapper. If a normalized solution is of the order of one but can exceed one, then the softplus function $softplus(z) = log(1+e^z)$ can be used as a wrapper instead of the sigmoid function. For the problem considered here, Figure \ref{fig:COM_PINN_backward} shows that the sigmoid and  softplus wrapper functions give similar results. 

\subsection{Data assimilation in the backward problem}
\begin{figure}[h]
	\centering
	\includegraphics[angle=0,width=0.5\textwidth]{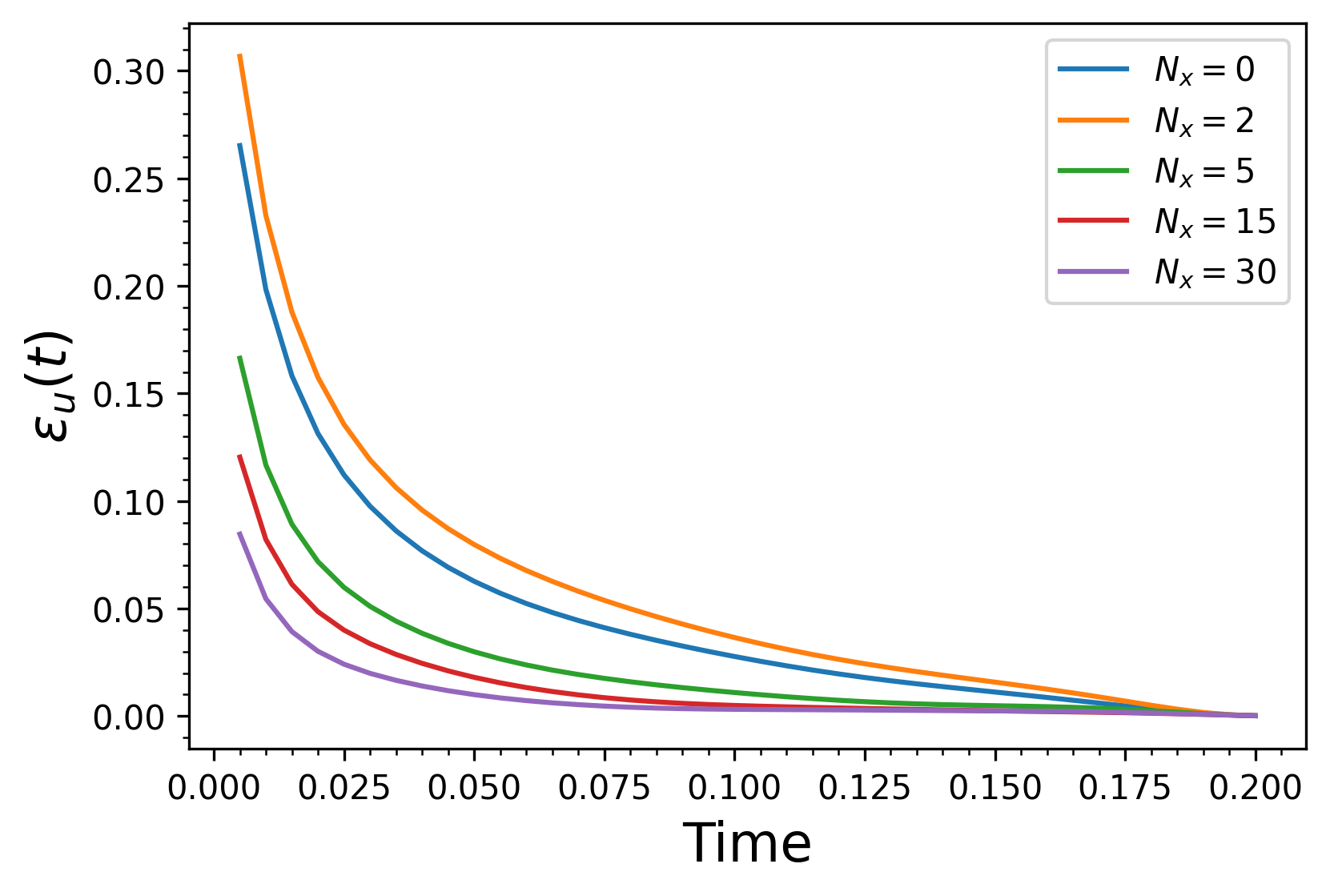}
	\caption{Relative $L_2$ errors in the PINN solutions $\hat{u}$ versus time as functions of the number of $u$ measurements, $N_x$.}
	\label{fig:error_Nc_PINN_backward}
\end{figure}

In many source identification problems, the measurements of $u$ could be available at multiple times and not only at the terminal time. 
The measurements of $u$ can be included in the PINN formulation of the backward ADE solution in the same way as for the inverse ADE, i.e., by adding the $\tilde{L}_{data}$ (with $N_k=N_h=0$) in the loss function \eqref{eq:loss_backward_ADE}.
We randomly select $N_x$ locations where $u$ is sampled over time and assume that 100 measurements of $u$ over time are available at each of the $N_x$ locations.  
Figure \ref{fig:error_Nc_PINN_backward} shows the relative L2 error as a function of time for $N_x$ ranging from 0 to 30. Our results show that at all times, the L2 error decreases with increasing $N_x$. At $t=0$ (when the errors are maximum), the errors are 30\% and 8\% for $N_x = 0$ and 30, respectively.
Snapshots of the backward PINN solution obtained with the measurements from 30 sensors ($N_x = 30$) and the corresponding point errors are shown in Figure \ref{fig:PINN_backward_sigmoid_0.2_Nx_30}.
The comparison with Figure \ref{fig:PINN_backward_sigmoid_0.2}, which depicts the solution obtained without $u$ measurements, demonstrates that $u$ measurements significantly improve solutions, especially at times close to 0 and near the (unknown) source location, with the maximum point errors decreasing from 1.2 to 0.2 at $t=0.05$ and from 120 to 60 at $t=0$.    

\section{Conclusion}\label{sec:conclusions}
In this study, we presented the PINN method for parabolic PDEs with strongly perturbed initial conditions. We demonstrated that an accurate PINN solution can be obtained after a PDE is normalized such that the amplitude of the perturbation in the initial condition does not decay in time. Such normalization of the solution allows an accurate approximation of the solution with a relatively small DNN, which is a necessary condition for obtaining a convergent PINN solution. We also developed criteria for the weights in the loss function and an adaptive-in-time sampling of the PDE's residuals that significantly improved the accuracy and efficiency of of the PINN method. 

As an example of the parabolic PDE, we considered a time-dependent ADE with a Gaussian source initial condition and a non-uniform advection velocity field described by the Darcy flow equation with the heterogeneous hydraulic conductivity. For the forward ADE problem (known velocity field), we demonstrated that an accurate PINN solution of the normalized ADE can be obtained with a relatively small DNN size and the number of residual points as long as the residual points are placed adaptively in time. For the same DNN size and the number of residual points, the PINN solution of the original ADE was found to have large errors and non-physical negative values, regardless of how the residuals are sampled.  Also, we found that there is a range of weight values that produce an accurate PINN solution of the normalized equation. The proposed criteria for weights resulted in the weights values within that range. 

For the inverse ADE problem (unknown conductivity and the velocity field), we demonstrated that the proposed normalization, the criteria for weights, and adaptive residual sampling strategy produced accurate estimates of the ADE solution as well as estimates of the velocity and conductivity fields in combination with sparse measurements of the ADE solution and conductivity. 

Finally, we demonstrated that the proposed methods allow for obtaining an accurate PINN solution of the backward ADE equation, including estimates of the center of mass and the mass of the tracer as functions of time within several percent of the reference values. With sparse measurements of the tracer concentration, the maximum L2 error in the backward solution (which occurs close to the location of the source and the time of release) was reduced from more than 30\% to less than 9\%.

\section{Acknowledgments}
This research was partially supported by the U.S. Department of Energy (DOE) Advanced Scientific Computing (ASCR) program. Pacific Northwest National Laboratory is operated by Battelle for the DOE under Contract DE-AC05-76RL01830.

\bibliographystyle{plain}
\bibliography{reference.bib}
\appendix
\section{Figures of PINN simulation results and errors}
\begin{figure}[h]
\begin{subfigure}[h]{0.32\textwidth}
\includegraphics[width=\linewidth]{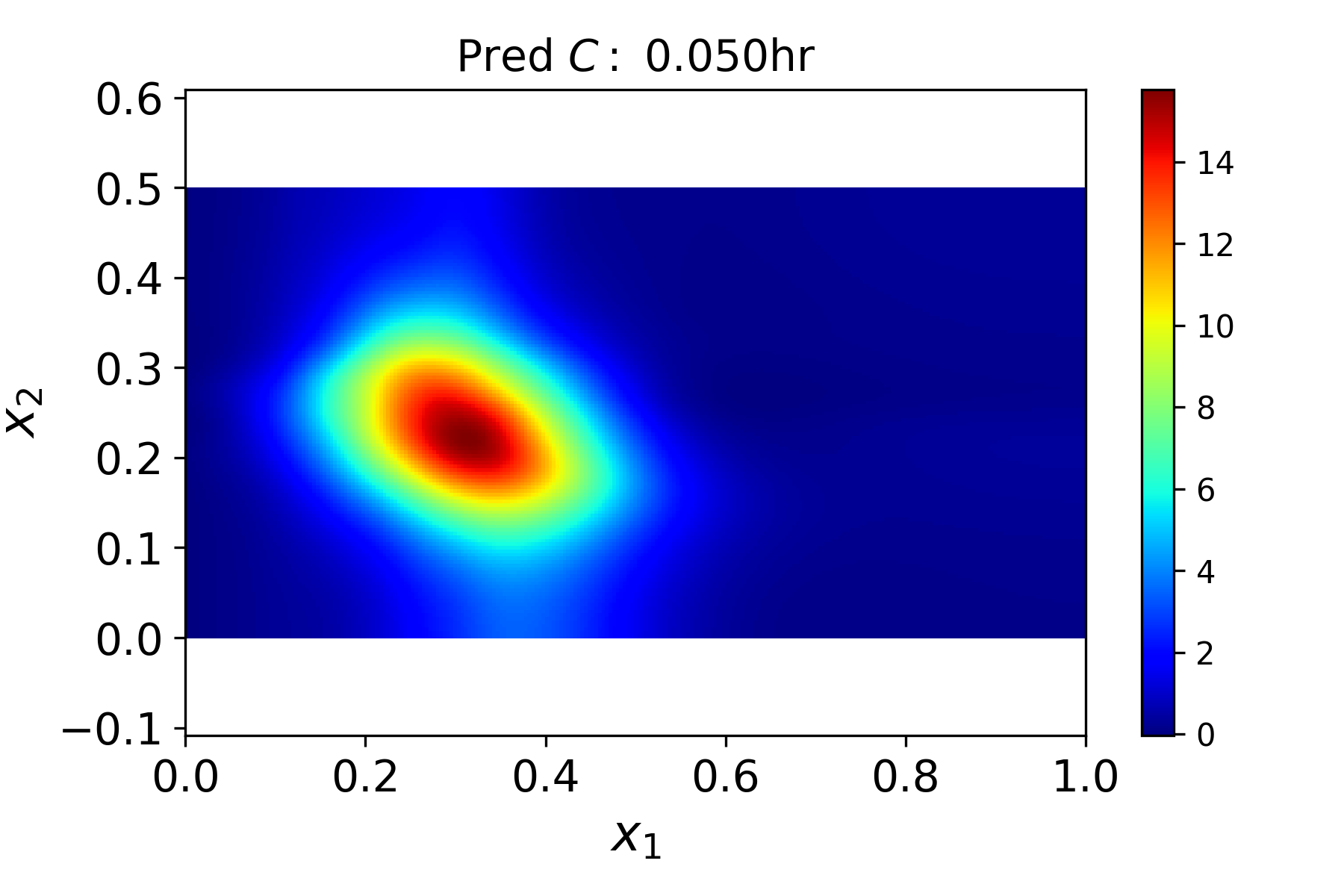}
\end{subfigure}
\begin{subfigure}[h]{0.32\textwidth}
\includegraphics[width=\linewidth]{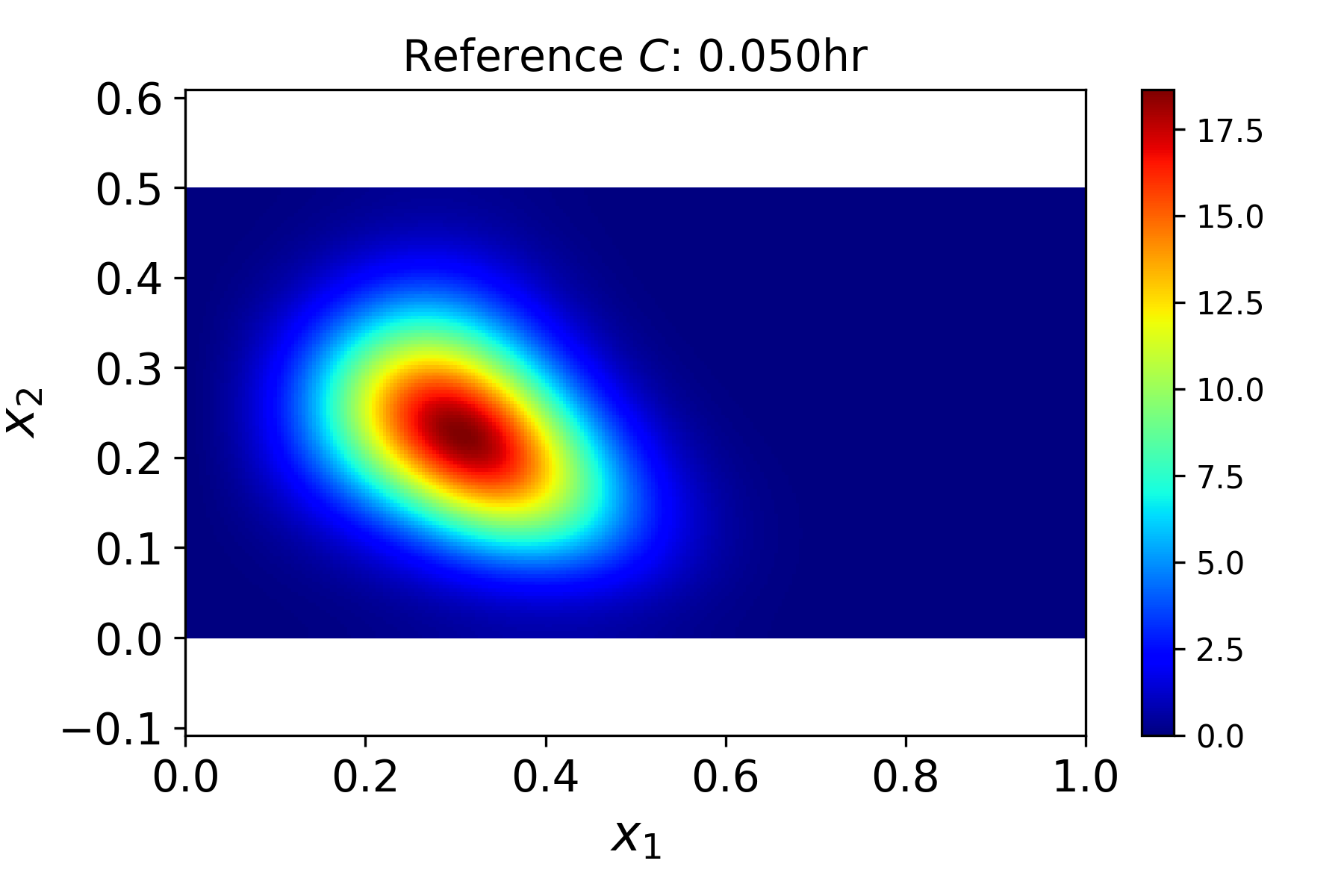}
\end{subfigure}
\begin{subfigure}[h]{0.32\textwidth}
  \includegraphics[width=\linewidth]{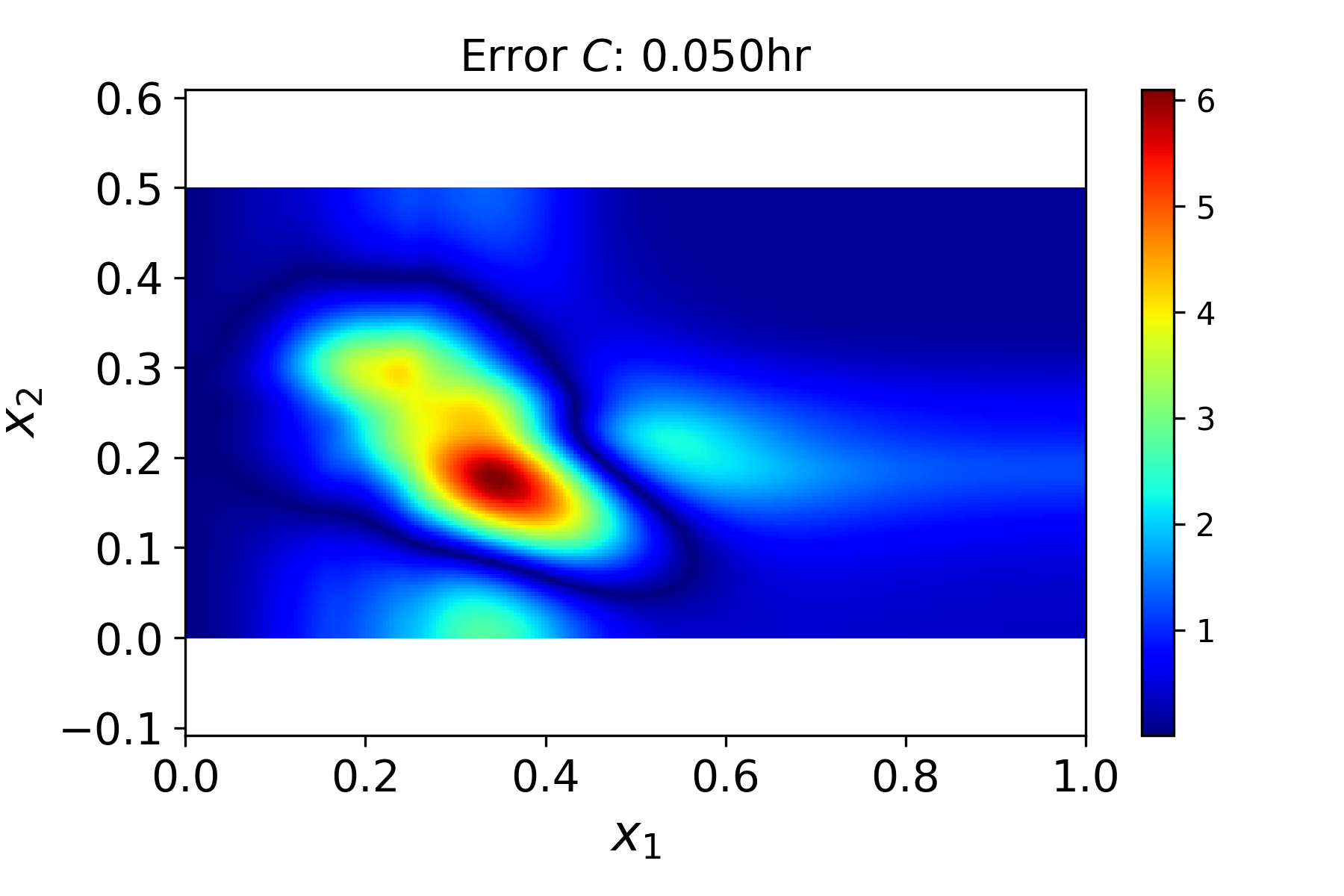}
\end{subfigure}

\begin{subfigure}[h]{0.32\textwidth}
\includegraphics[width=\linewidth]{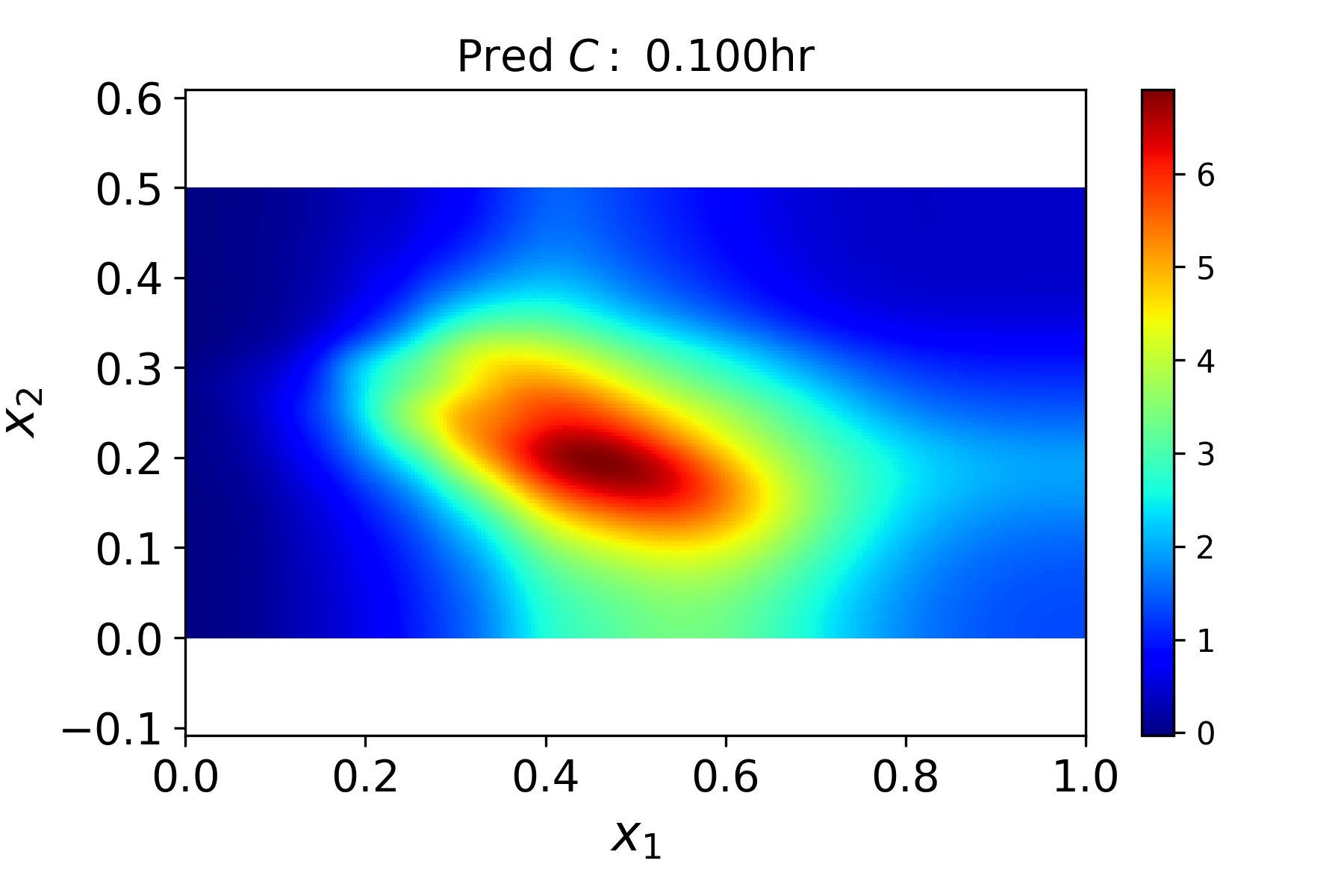}
\end{subfigure}
\begin{subfigure}[h]{0.32\textwidth}
\includegraphics[width=\linewidth]{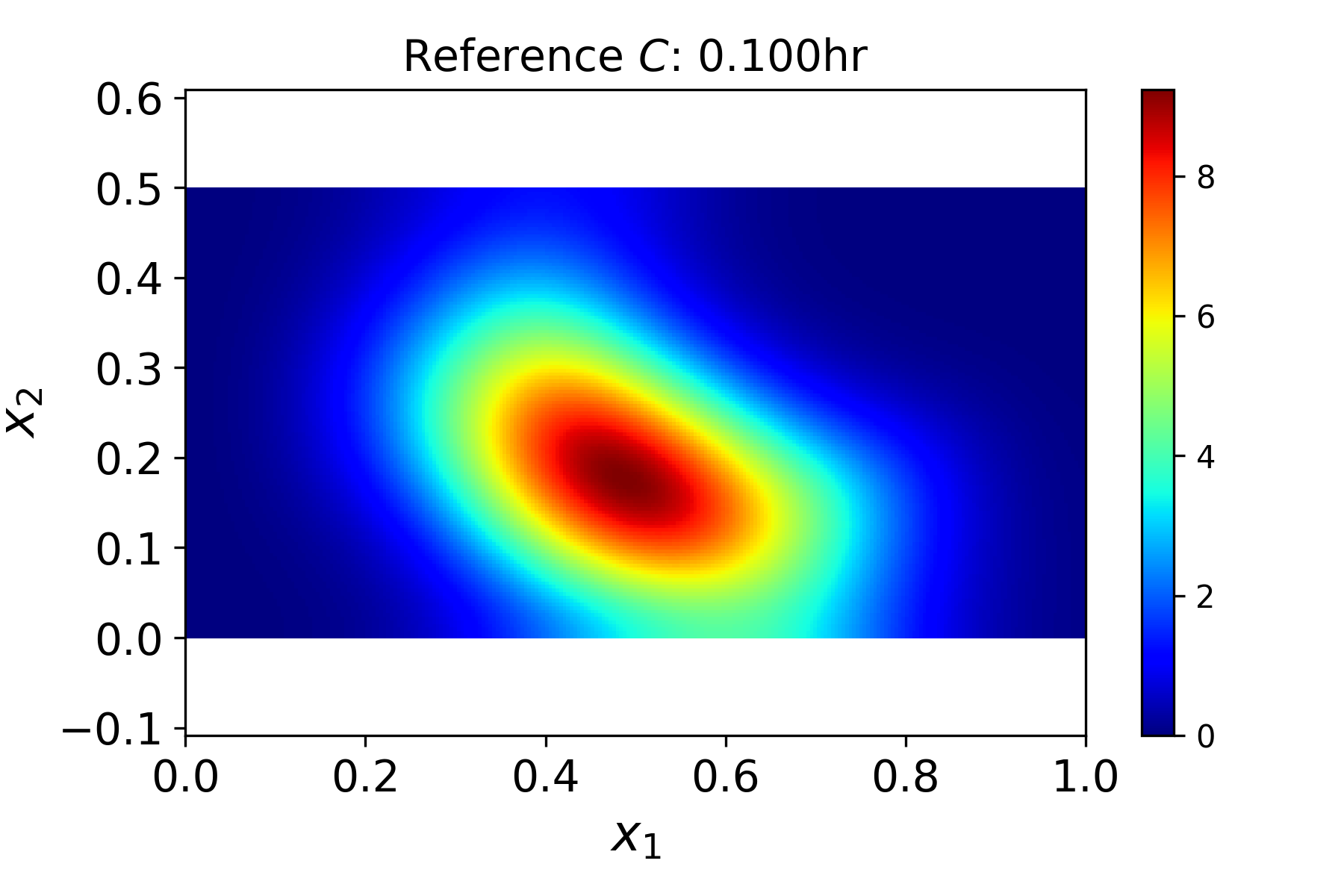}
\end{subfigure}
\begin{subfigure}[h]{0.32\textwidth}
  \includegraphics[width=\linewidth]{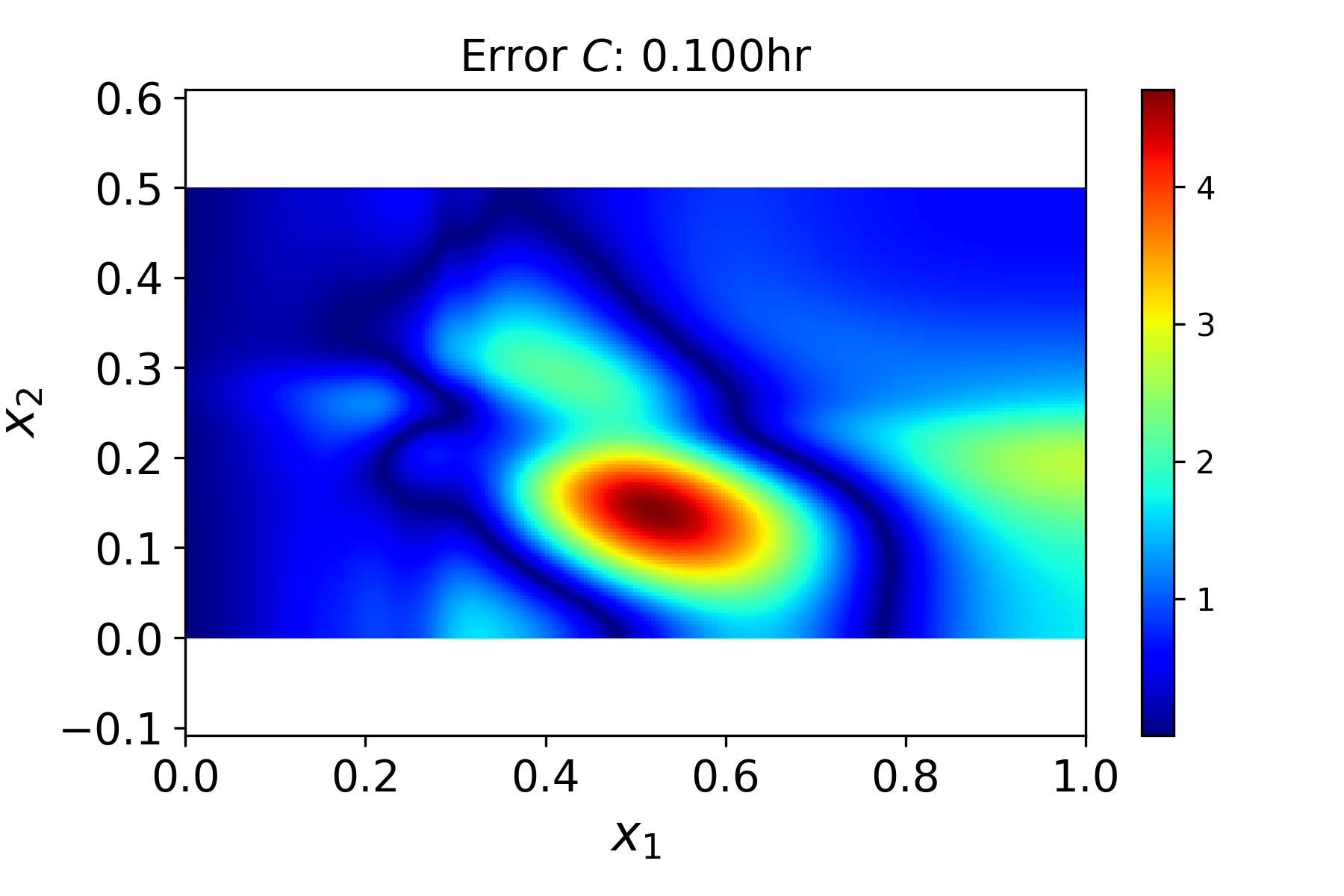}
\end{subfigure}

\begin{subfigure}[h]{0.32\textwidth}
\includegraphics[width=\linewidth]{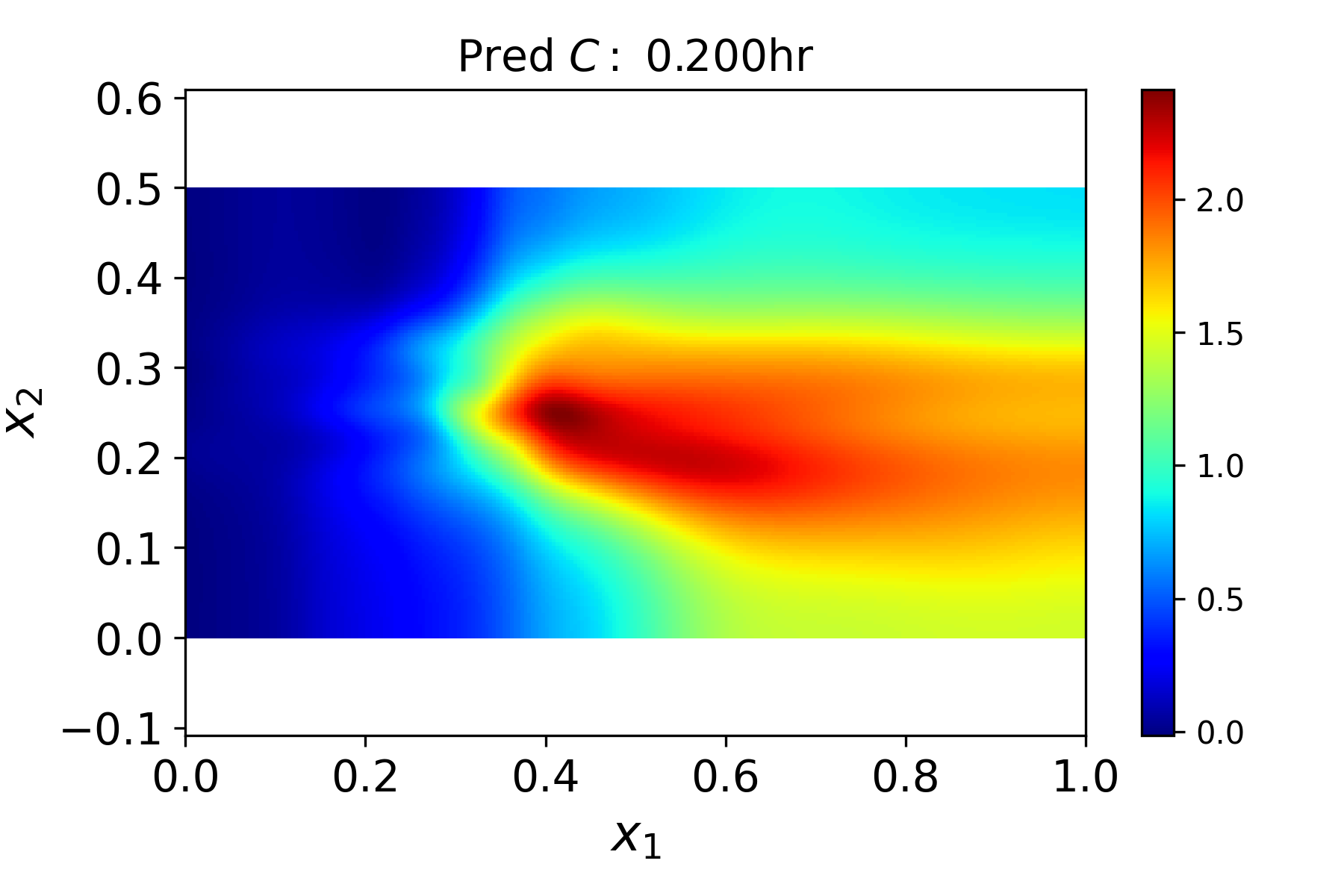}
\end{subfigure}
\begin{subfigure}[h]{0.32\textwidth}
\includegraphics[width=\linewidth]{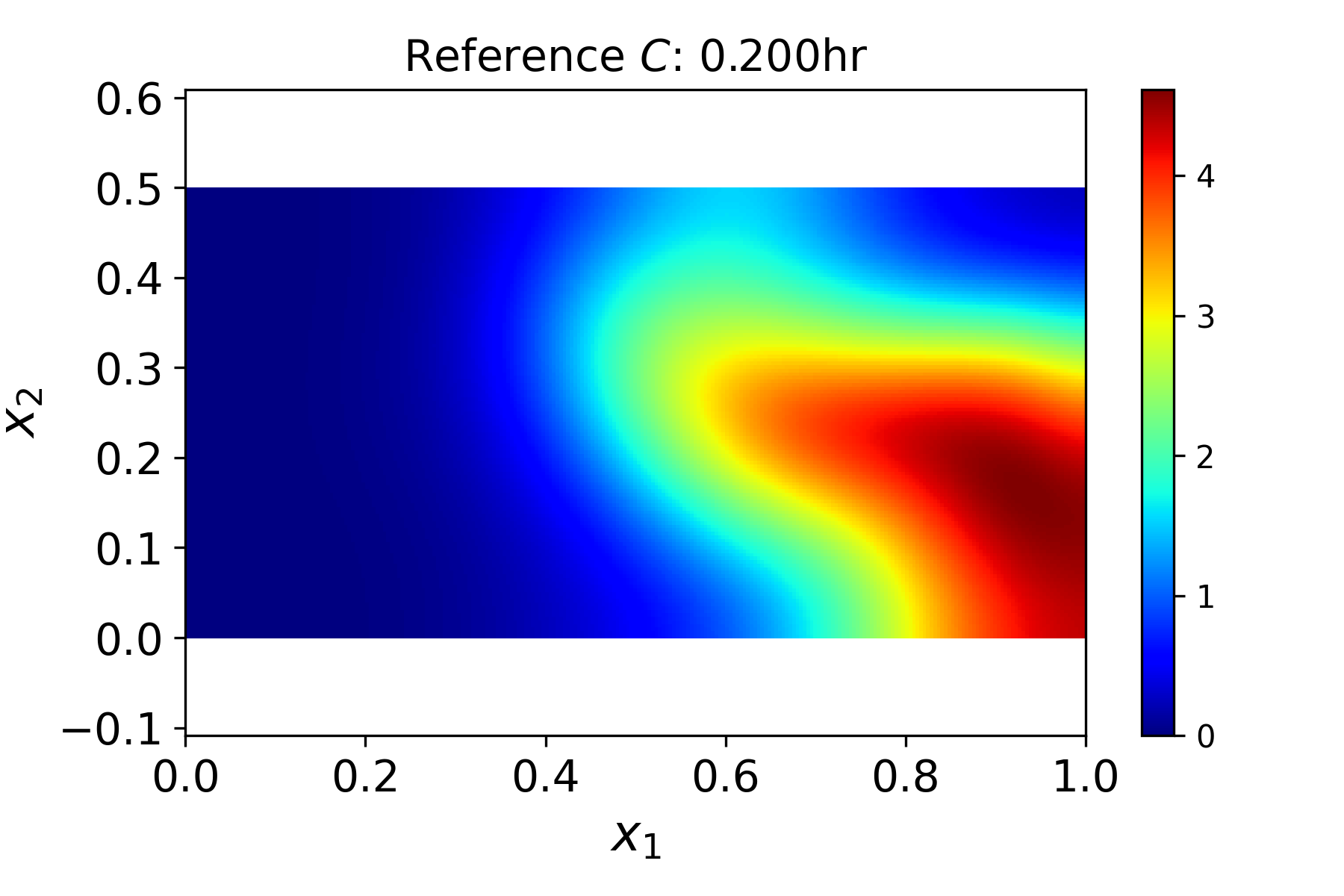}
\end{subfigure}
\begin{subfigure}[h]{0.32\textwidth}
\includegraphics[width=\linewidth]{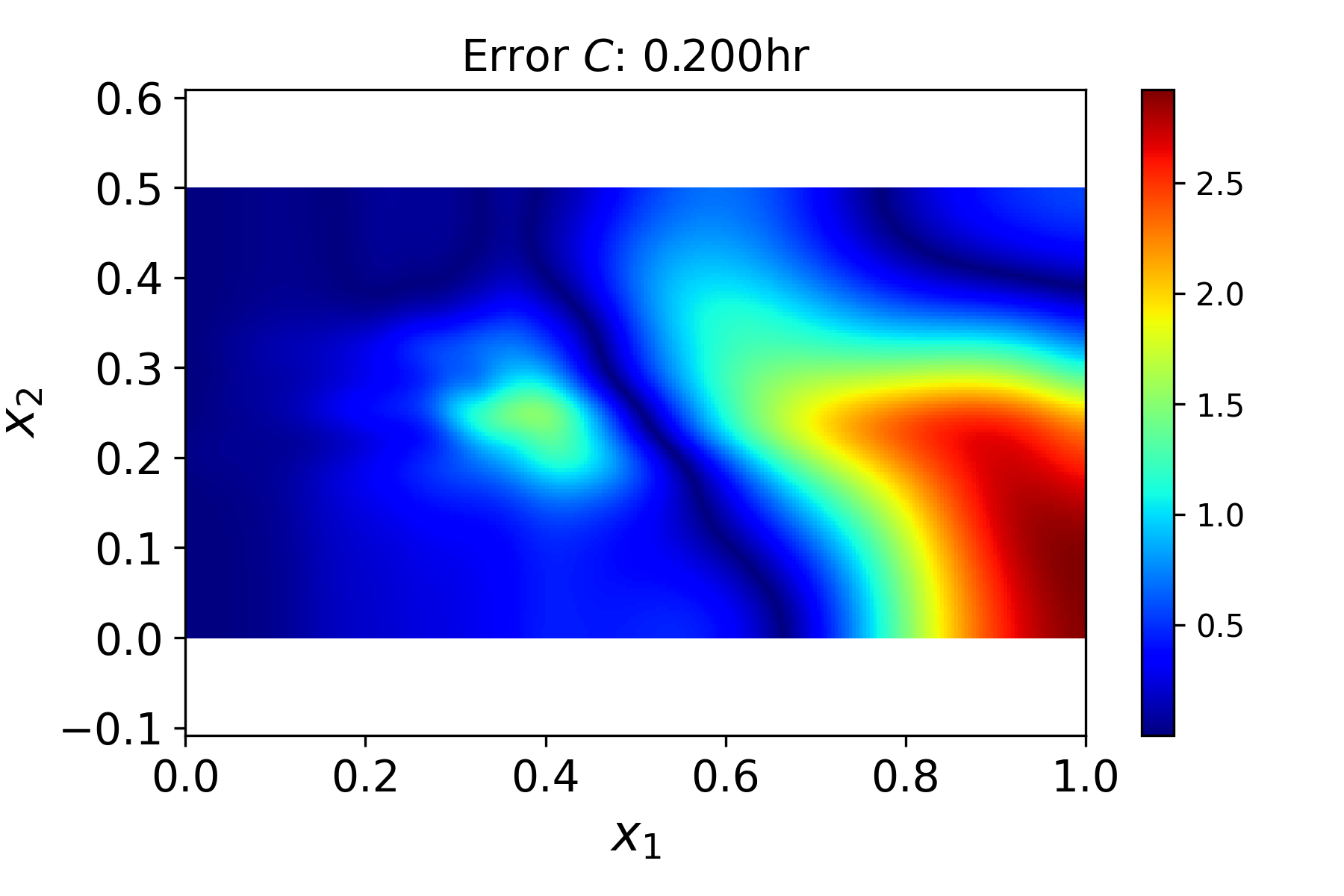}
\end{subfigure}

\begin{subfigure}[h]{0.32\textwidth}
\includegraphics[width=\linewidth]{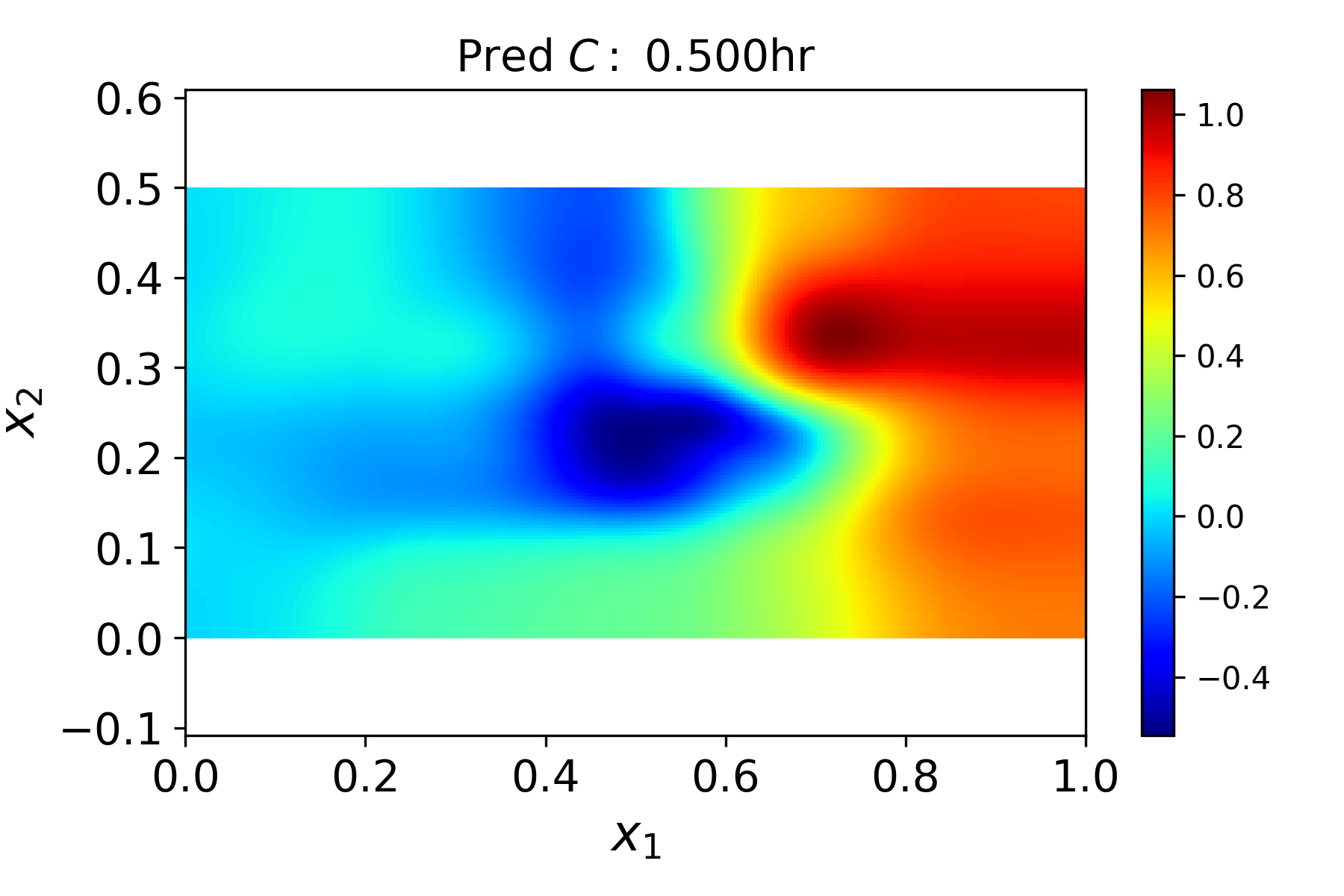}
\end{subfigure}
\begin{subfigure}[h]{0.32\textwidth}
\includegraphics[width=\linewidth]{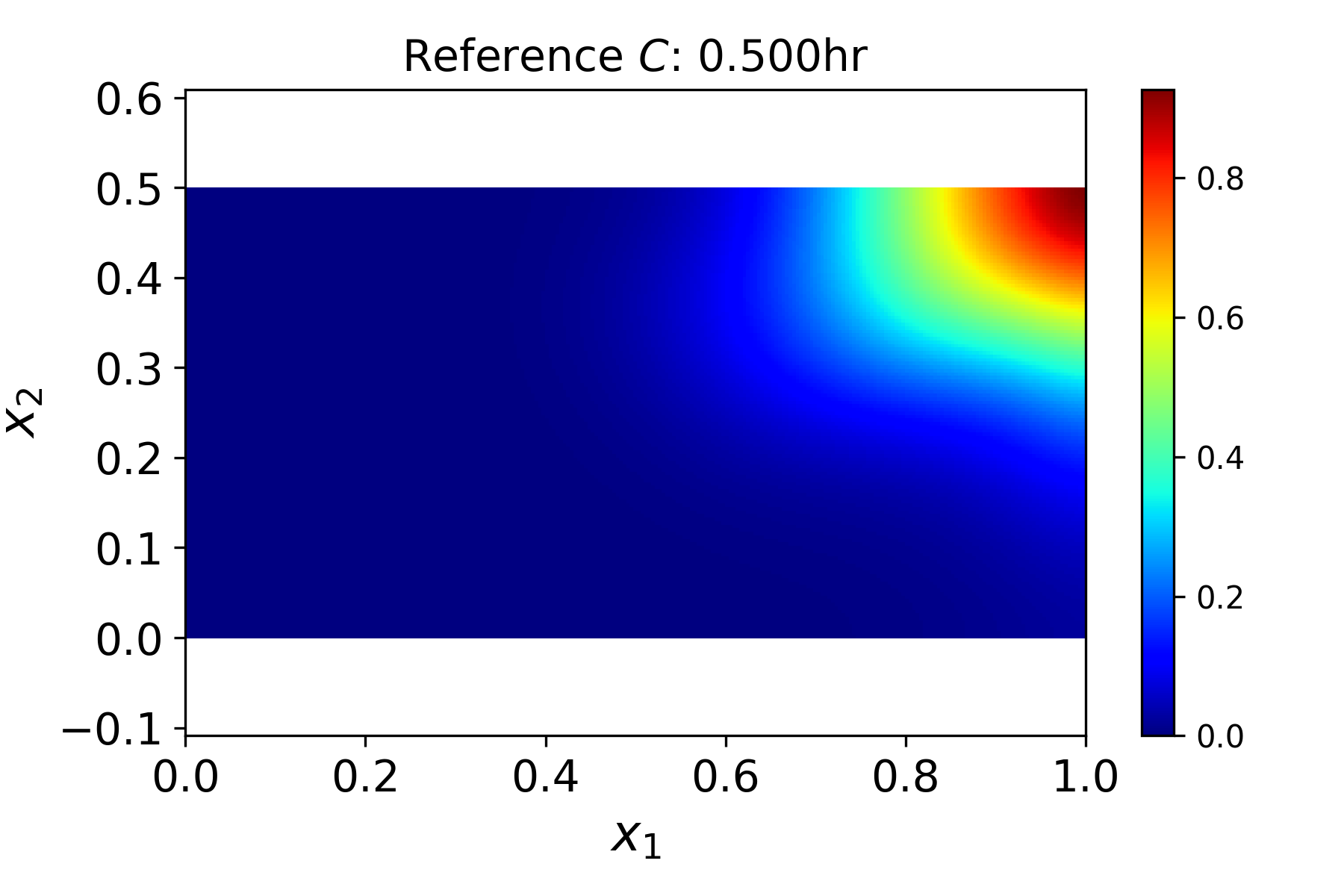}
\end{subfigure}
\begin{subfigure}[h]{0.32\textwidth}
  \includegraphics[width=\linewidth]{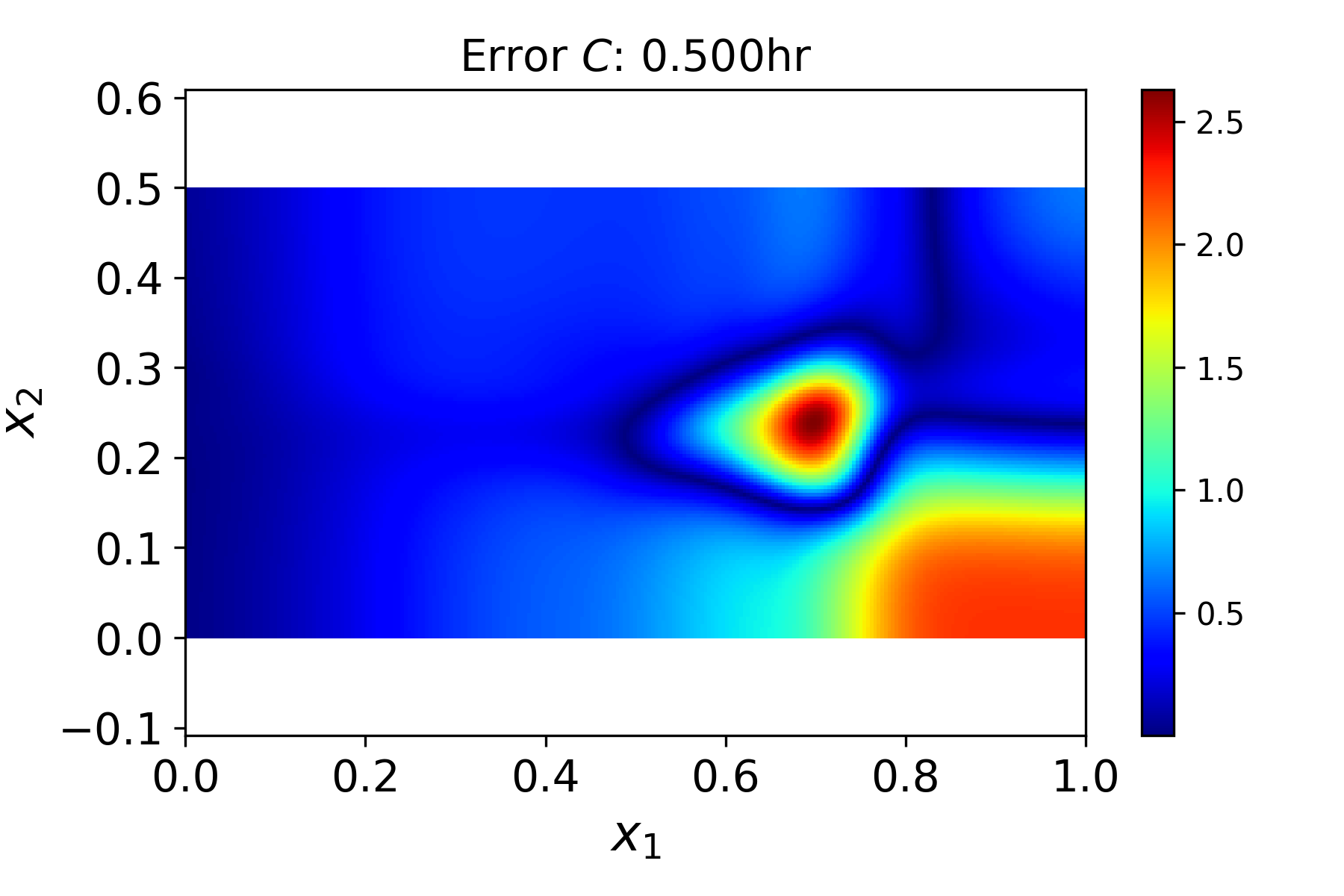}
\end{subfigure}
\caption{PINN solution $\hat{u}$ of the forward ADE \eqref{eq:ADE} (left column), reference $u$ field (central column), and the absolute point errors $\hat{u}-u$ (right column) at times $t=0.02$, 0.1, 0.2, and 0.5. The weights are set to $\lambda_{ic} = \lambda_{bcd}  = \lambda_{bcn} = 10$, $\lambda_{res} = 1$. This PINN solution has the smallest L2 error among all considered weight combinations in Section \ref{sec:PINN_limitations}. 
 }
\label{fig:PINN_forward_original}
\end{figure}
\begin{figure}[h]
\begin{subfigure}[h]{0.32\textwidth}
\includegraphics[width=\linewidth]{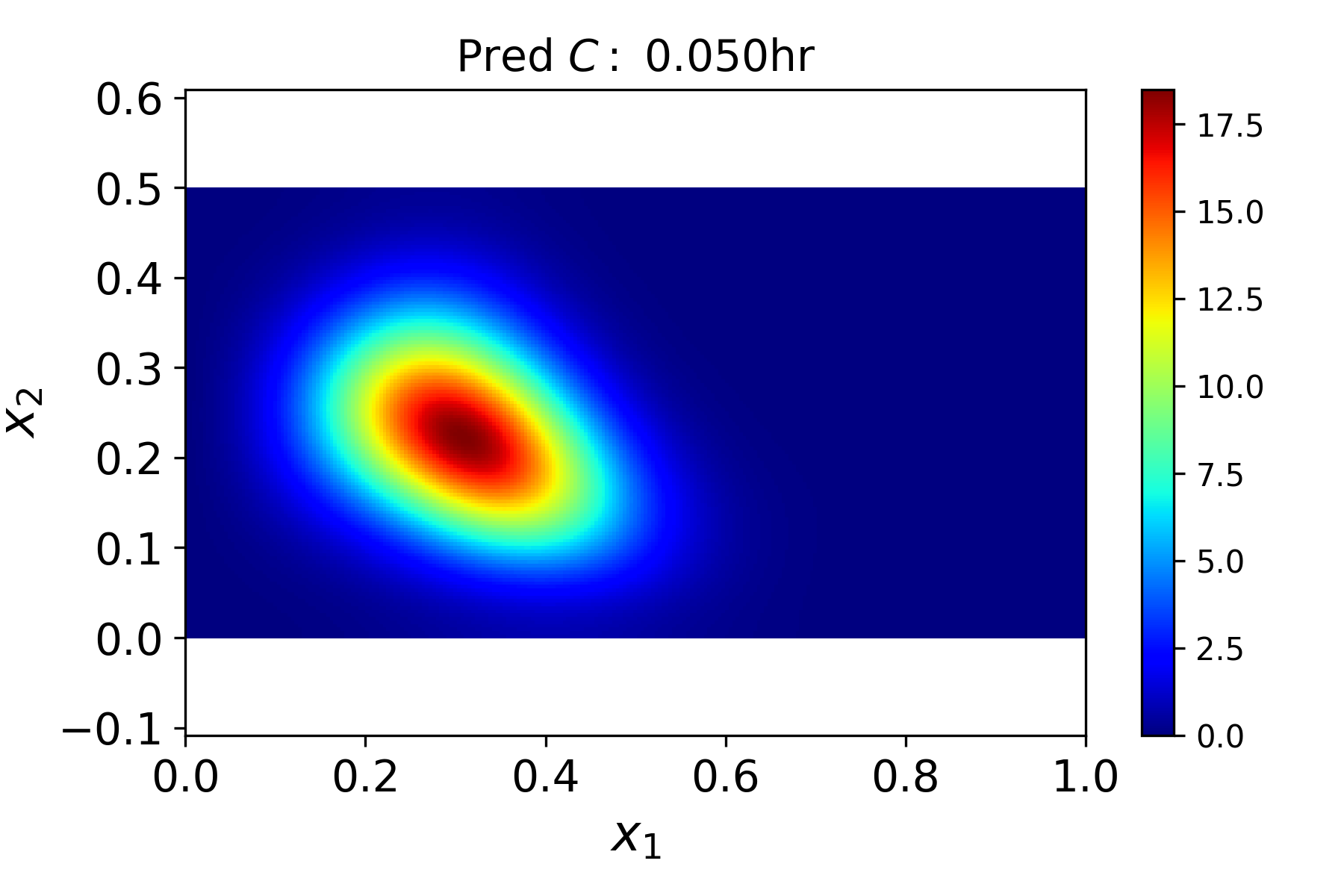}
\end{subfigure}
\begin{subfigure}[h]{0.32\textwidth}
\includegraphics[width=\linewidth]{draft_PINN_AD/figures/backward_plot/reference_t_0_050hr.png}
\end{subfigure}
\begin{subfigure}[h]{0.32\textwidth}
  \includegraphics[width=\linewidth]{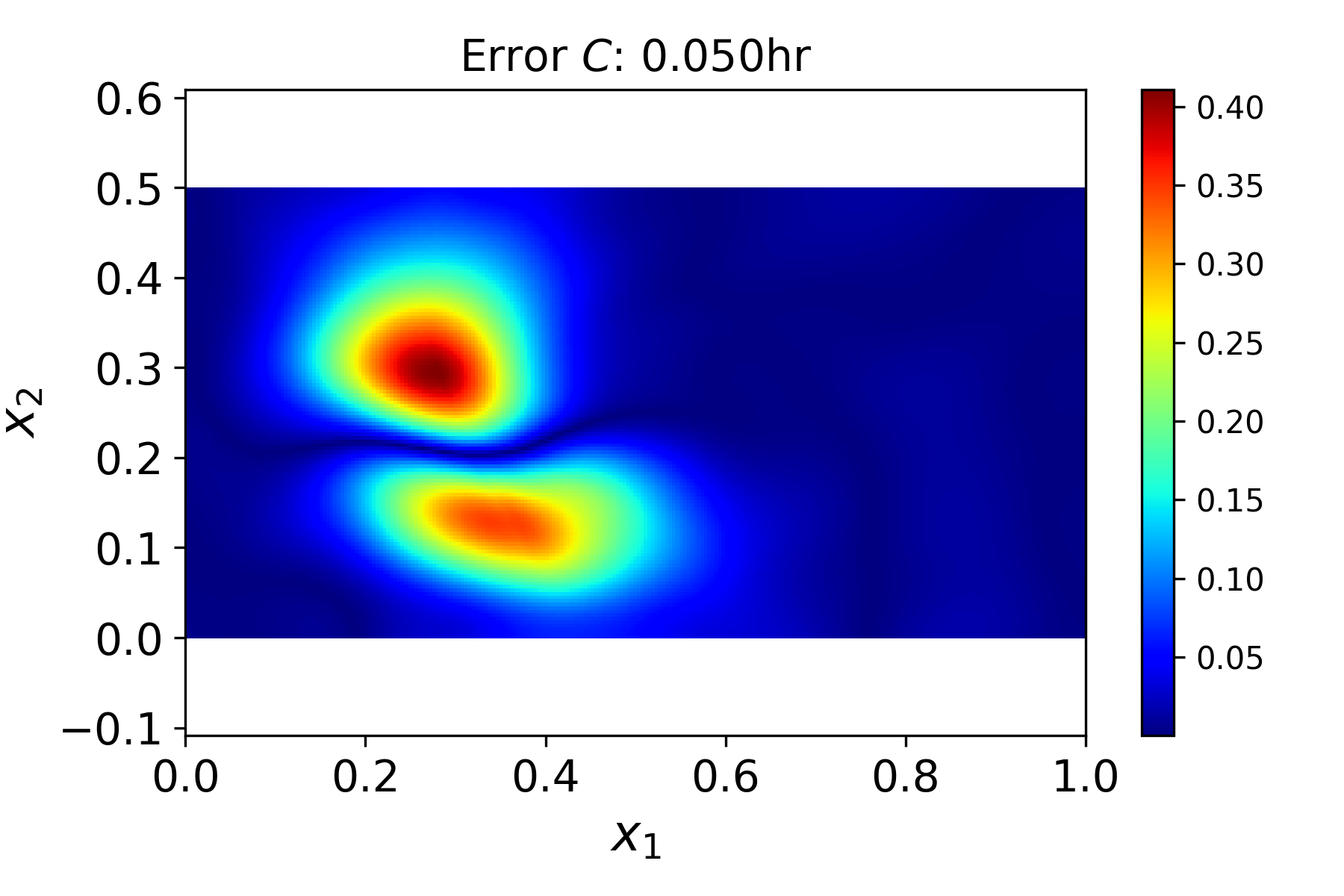}
\end{subfigure}

\begin{subfigure}[h]{0.32\textwidth}
\includegraphics[width=\linewidth]{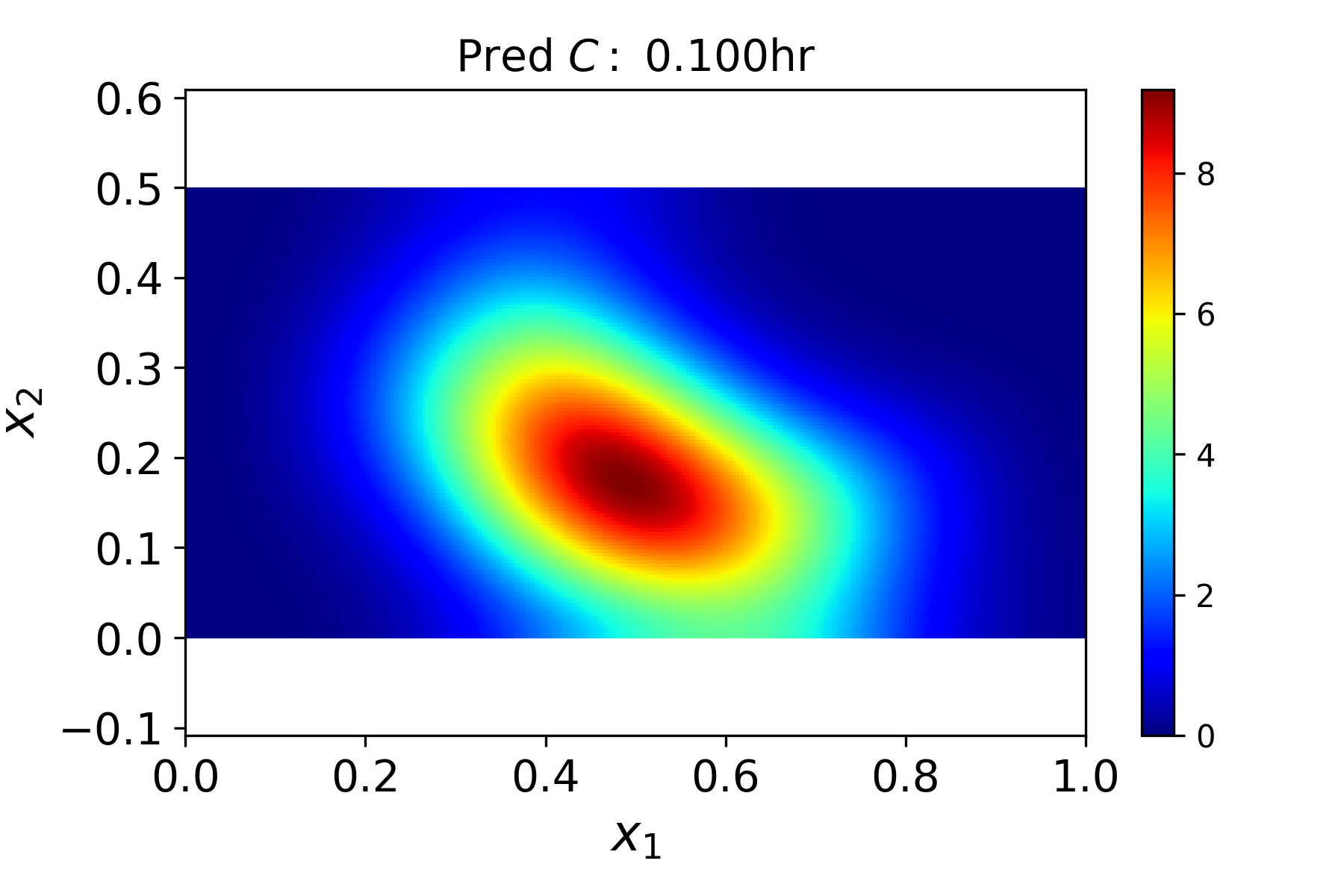}
\end{subfigure}
\begin{subfigure}[h]{0.32\textwidth}
\includegraphics[width=\linewidth]{draft_PINN_AD/figures/backward_plot/reference_t_0_100hr.png}
\end{subfigure}
\begin{subfigure}[h]{0.32\textwidth}
  \includegraphics[width=\linewidth]{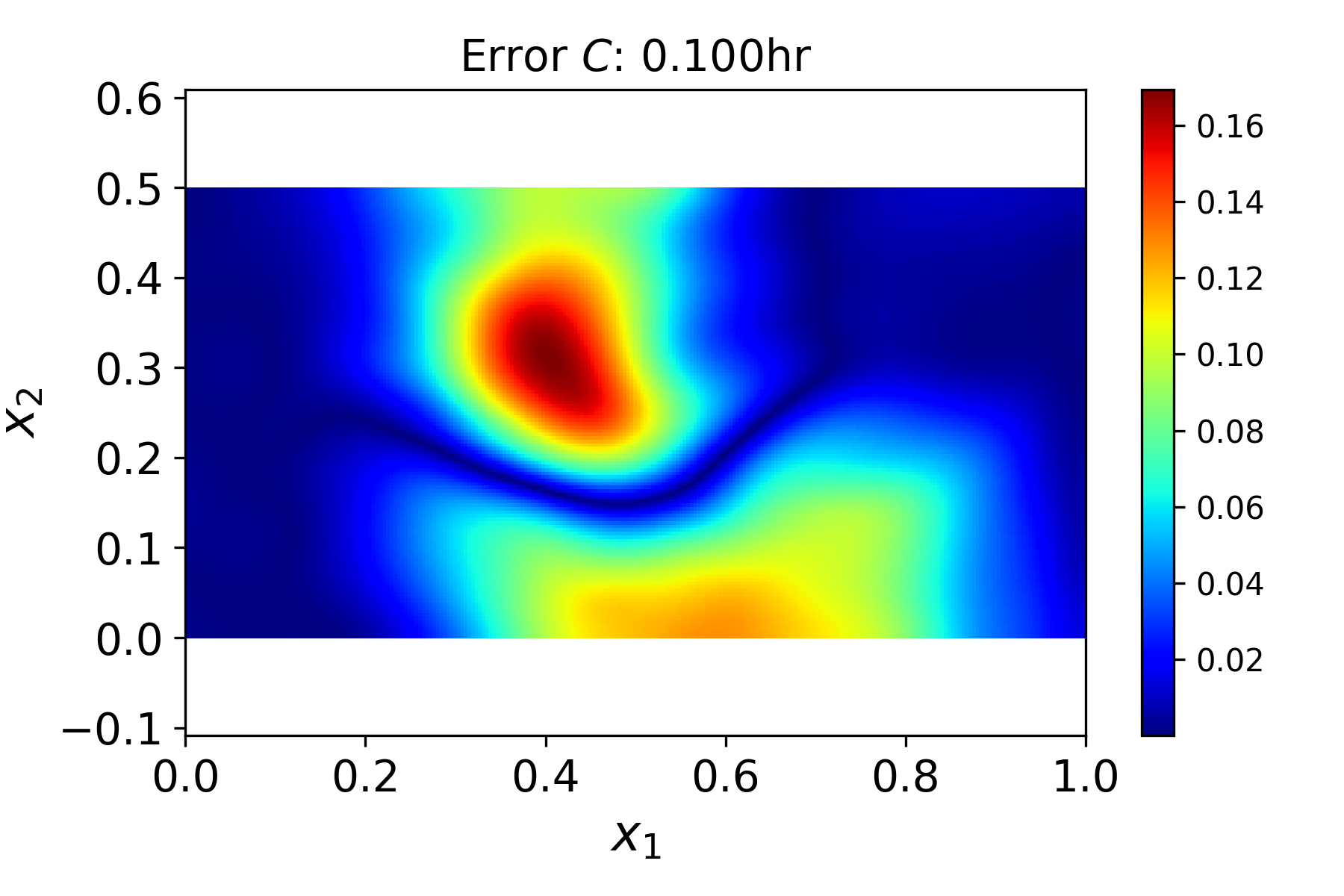}
\end{subfigure}

\begin{subfigure}[h]{0.32\textwidth}
\includegraphics[width=\linewidth]{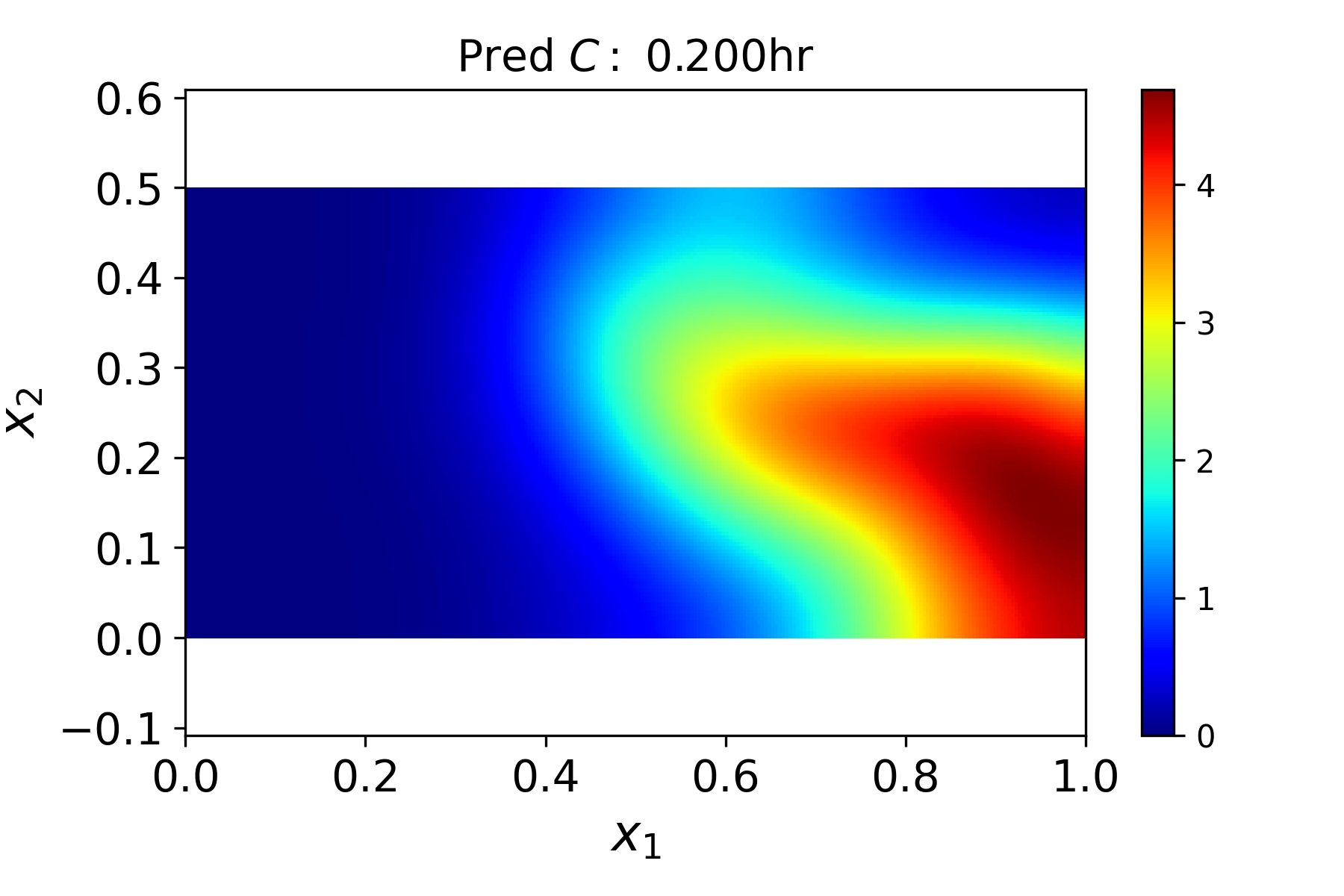}
\end{subfigure}
\begin{subfigure}[h]{0.32\textwidth}
\includegraphics[width=\linewidth]{draft_PINN_AD/figures/backward_plot/reference_t_0_200hr.png}
\end{subfigure}
\begin{subfigure}[h]{0.32\textwidth}
  \includegraphics[width=\linewidth]{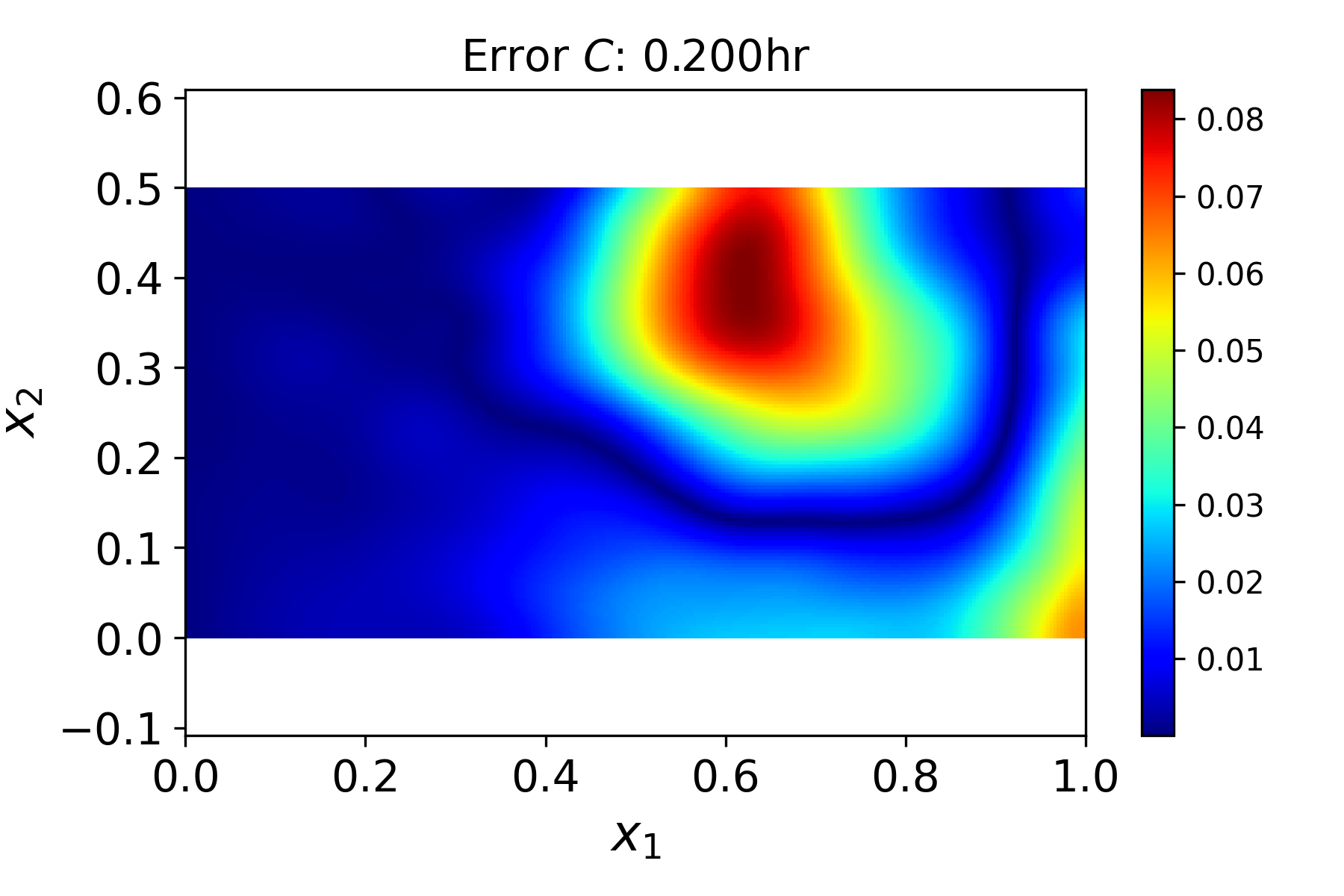}
\end{subfigure}

\begin{subfigure}[h]{0.32\textwidth}
\includegraphics[width=\linewidth]{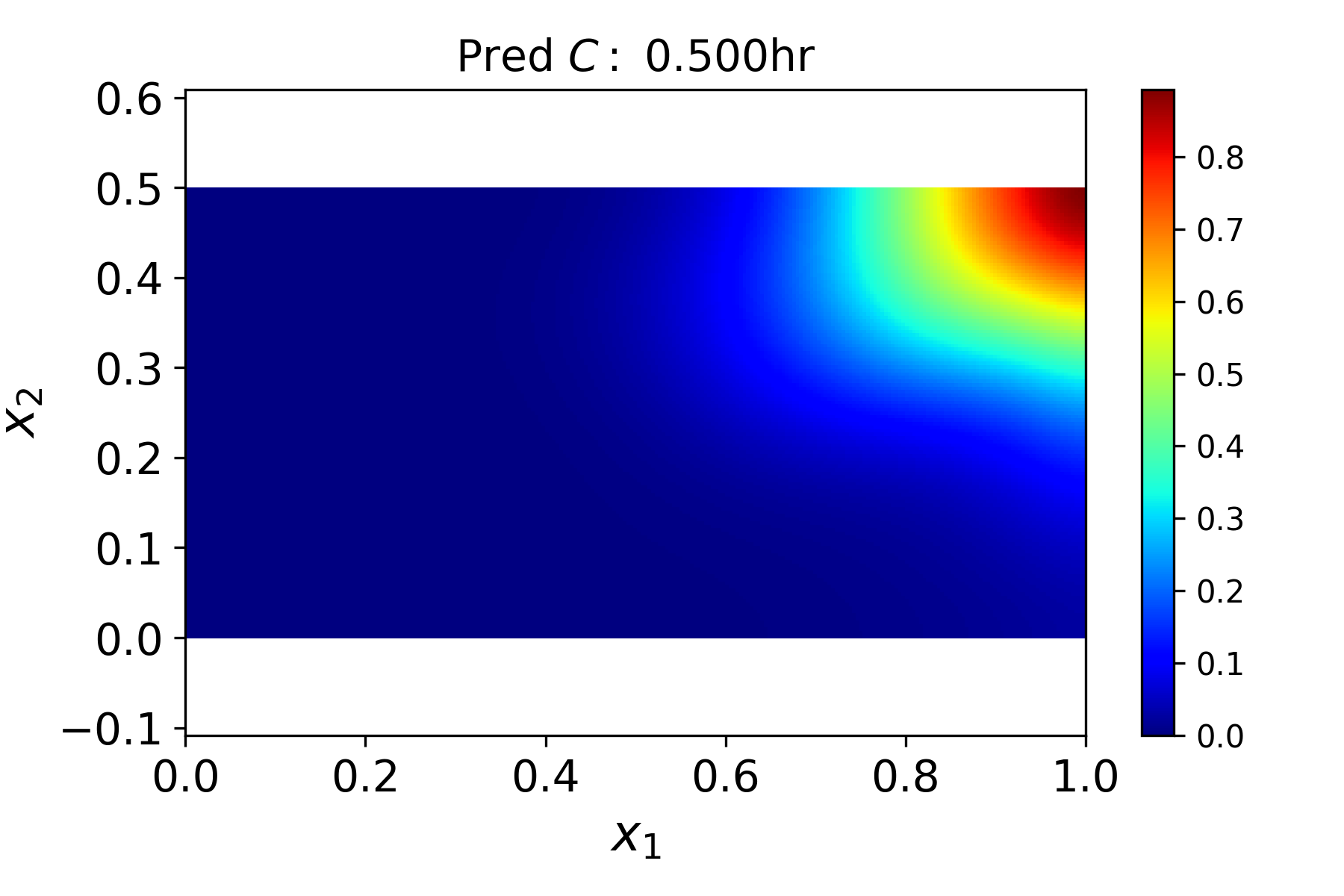}
\end{subfigure}
\begin{subfigure}[h]{0.32\textwidth}
\includegraphics[width=\linewidth]{draft_PINN_AD/figures/backward_plot/reference_t_0_500hr.png}
\end{subfigure}
\begin{subfigure}[h]{0.32\textwidth}
  \includegraphics[width=\linewidth]{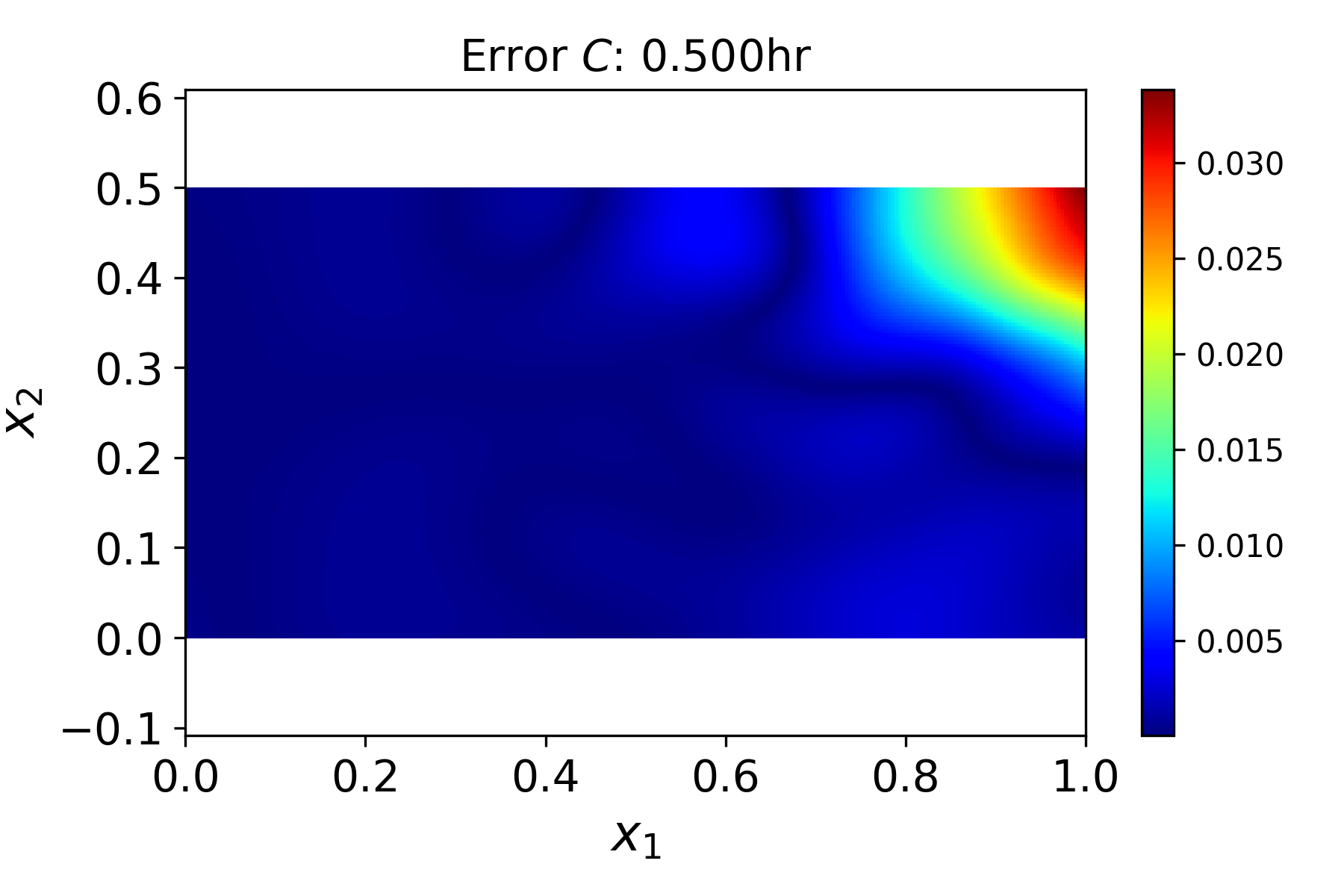}
\end{subfigure}
\caption{
PINN solution $\hat{u}$ of the normalized forward ADE \eqref{eq:ADE} (left column), reference normalized $\tilde{u}$ field (central column), and the absolute point errors $\hat{u}-\tilde{u}u$ (right column) at times $t=0.02$, 0.1, 0.2, and 0.5. The weights are set to $\lambda_{ic} = \lambda_{bcd} = 126$, $\lambda_{bcn} = \lambda_{res} = 1$.}
\label{fig:PINN_forward_normalized}
\end{figure}

\begin{figure}[h]
\begin{subfigure}[h]{0.33\textwidth}
\includegraphics[width=\linewidth]{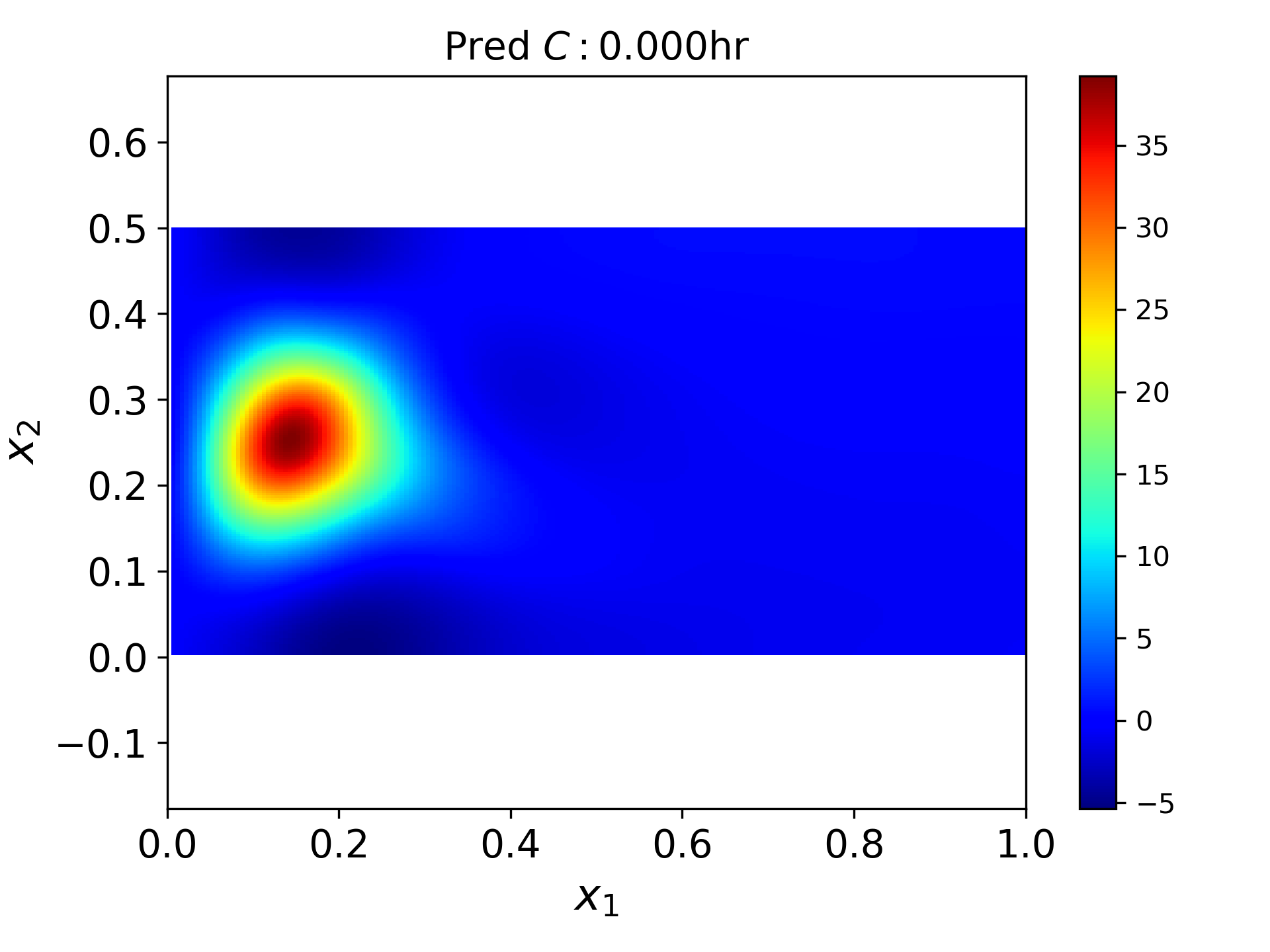}
\end{subfigure}\hfill
\begin{subfigure}[h]{0.33\textwidth}
\includegraphics[width=\linewidth]{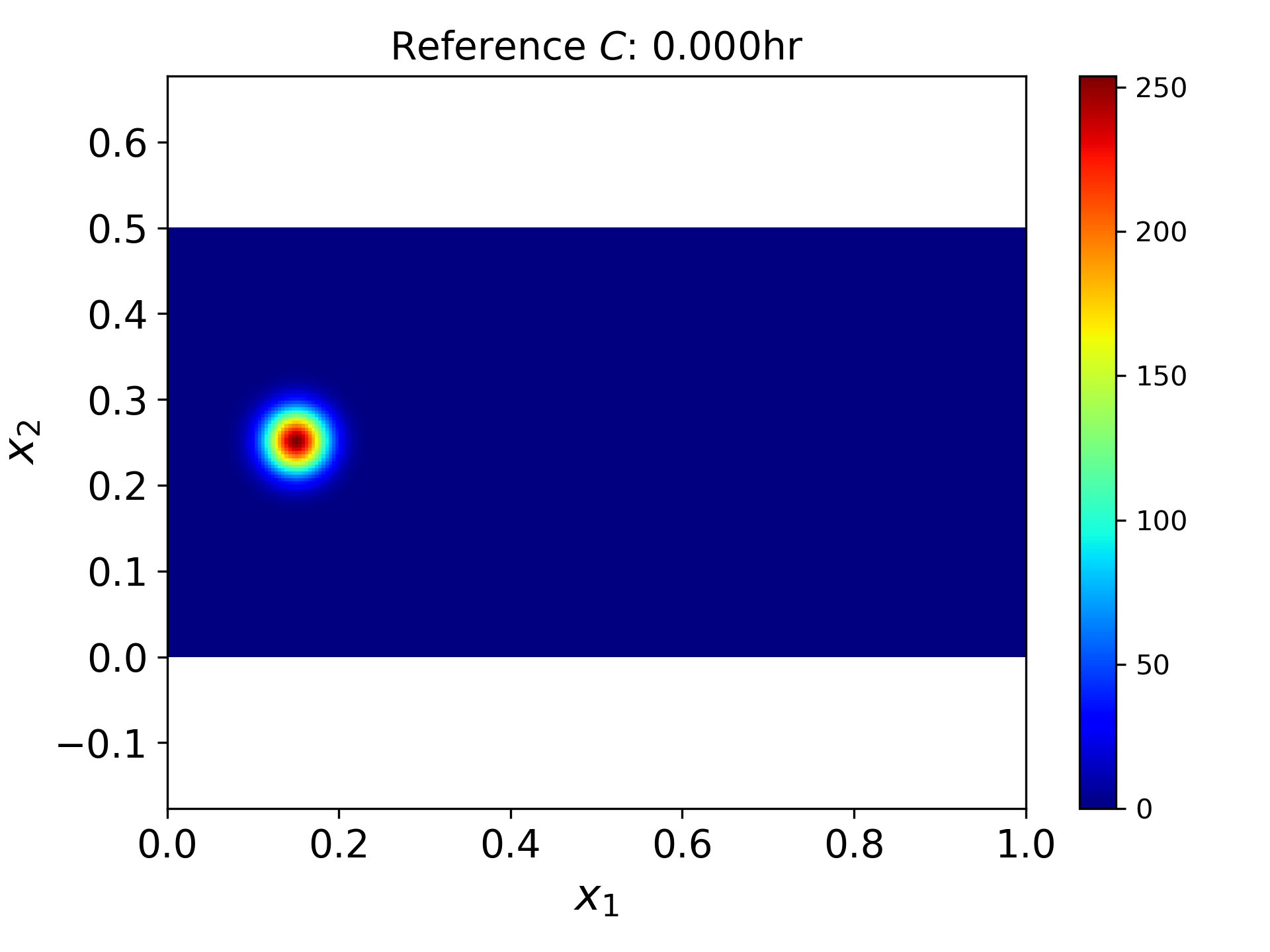}
\end{subfigure}\hfill
\begin{subfigure}[h]{0.33\textwidth}
\includegraphics[width=\linewidth]{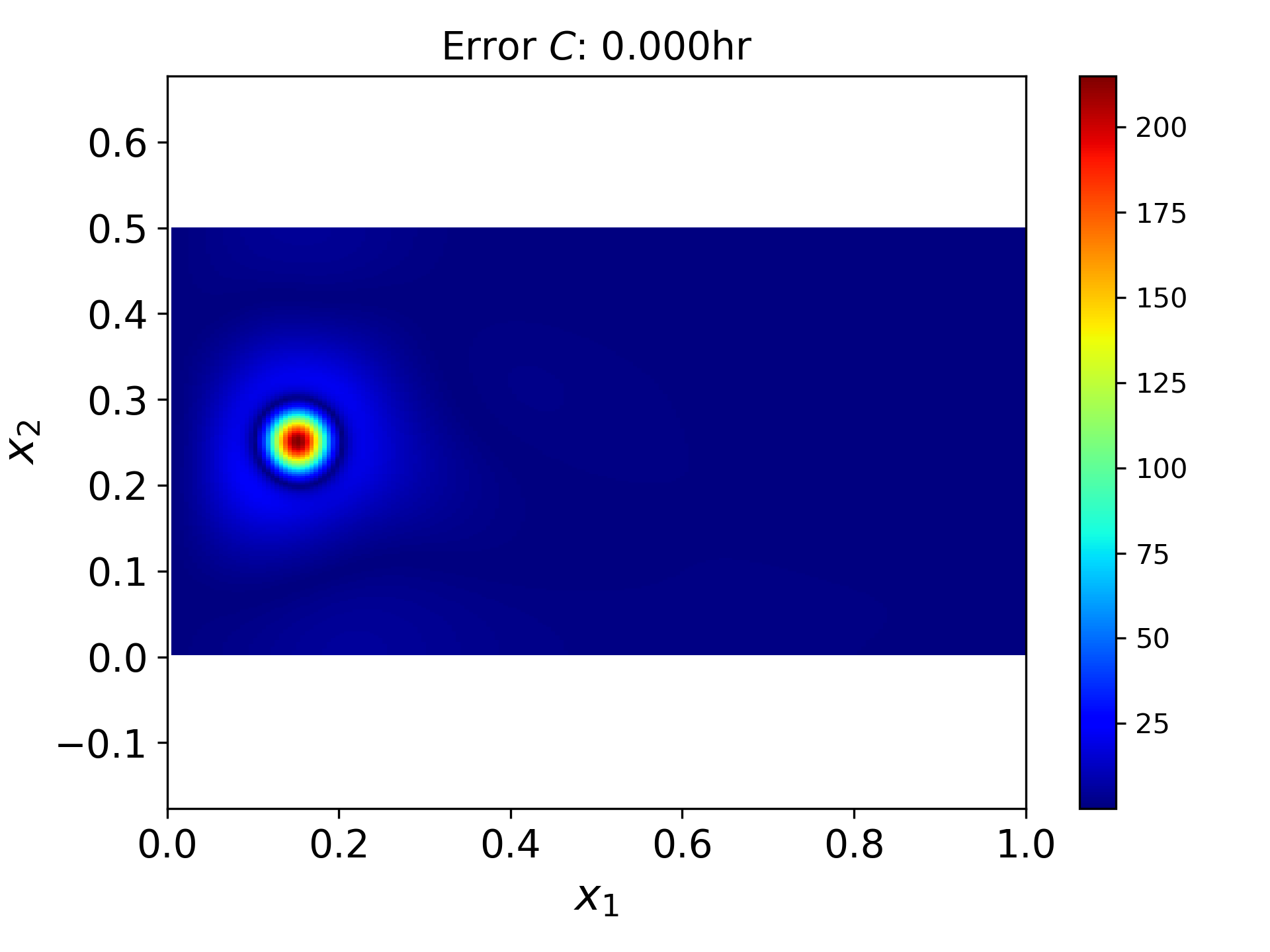}
\end{subfigure}

\begin{subfigure}[h]{0.33\textwidth}
\includegraphics[width=\linewidth]{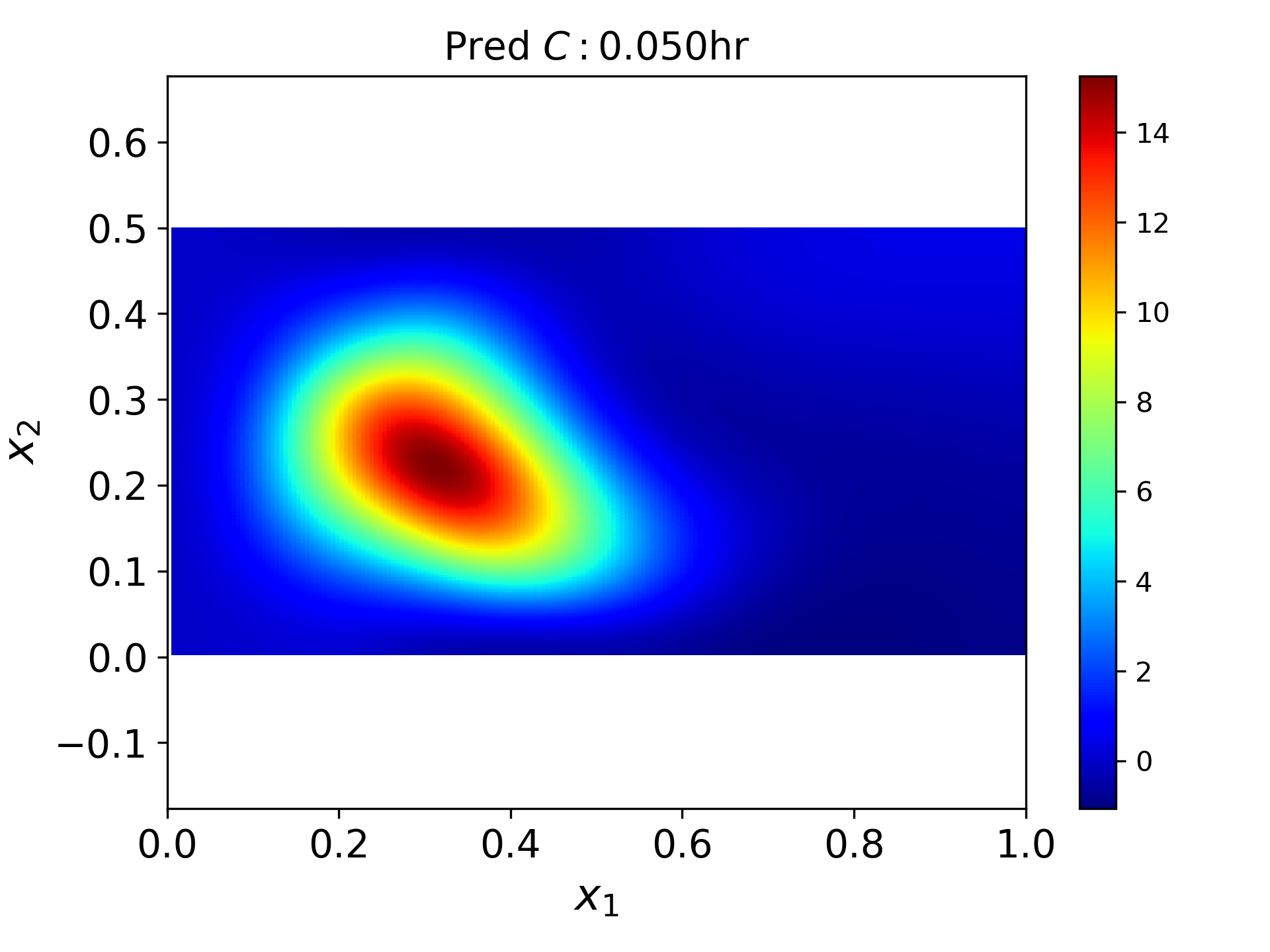}
\end{subfigure}\hfill
\begin{subfigure}[h]{0.33\textwidth}
\includegraphics[width=\linewidth]{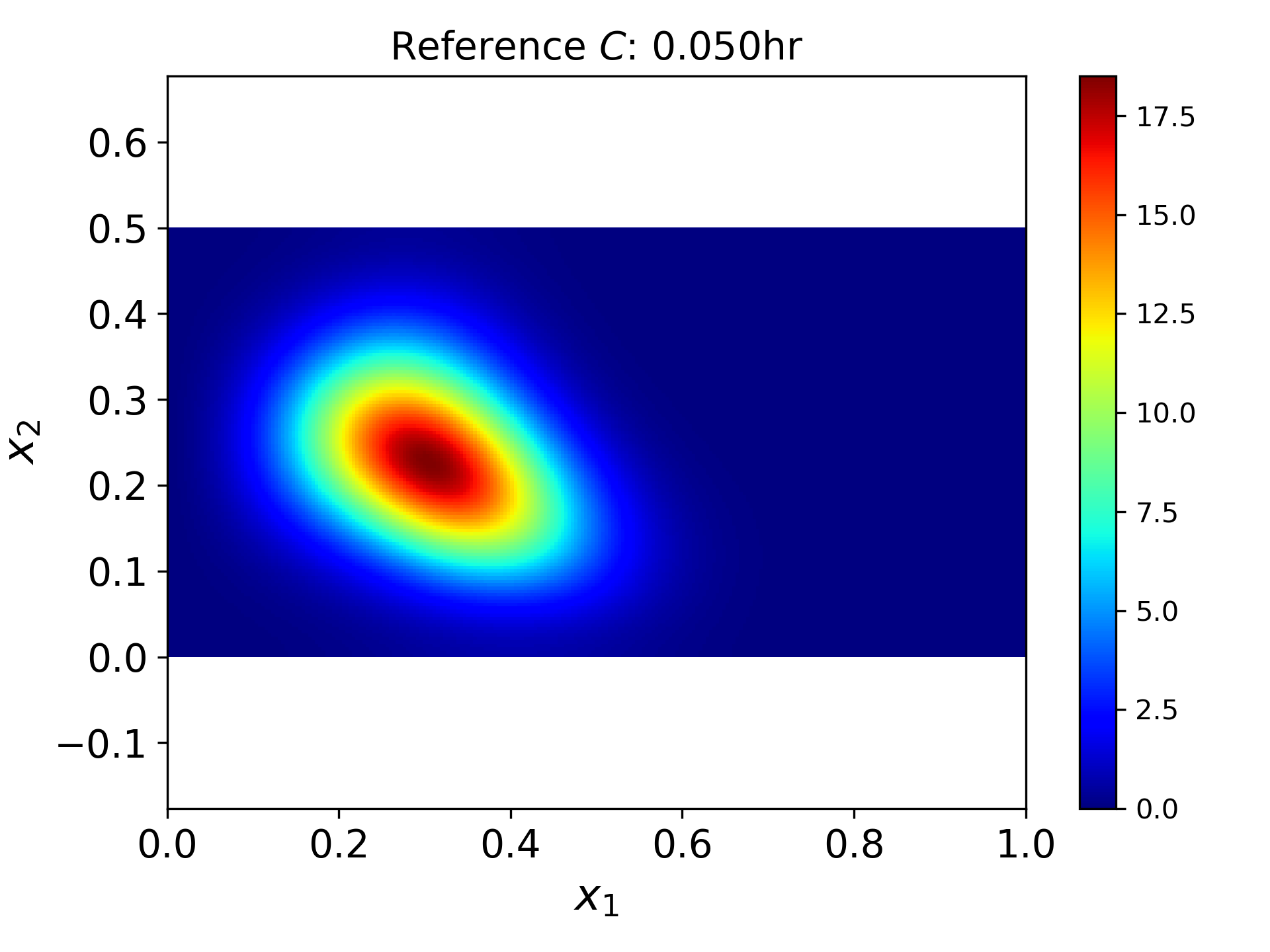}
\end{subfigure}\hfill
\begin{subfigure}[h]{0.33\textwidth}
\includegraphics[width=\linewidth]{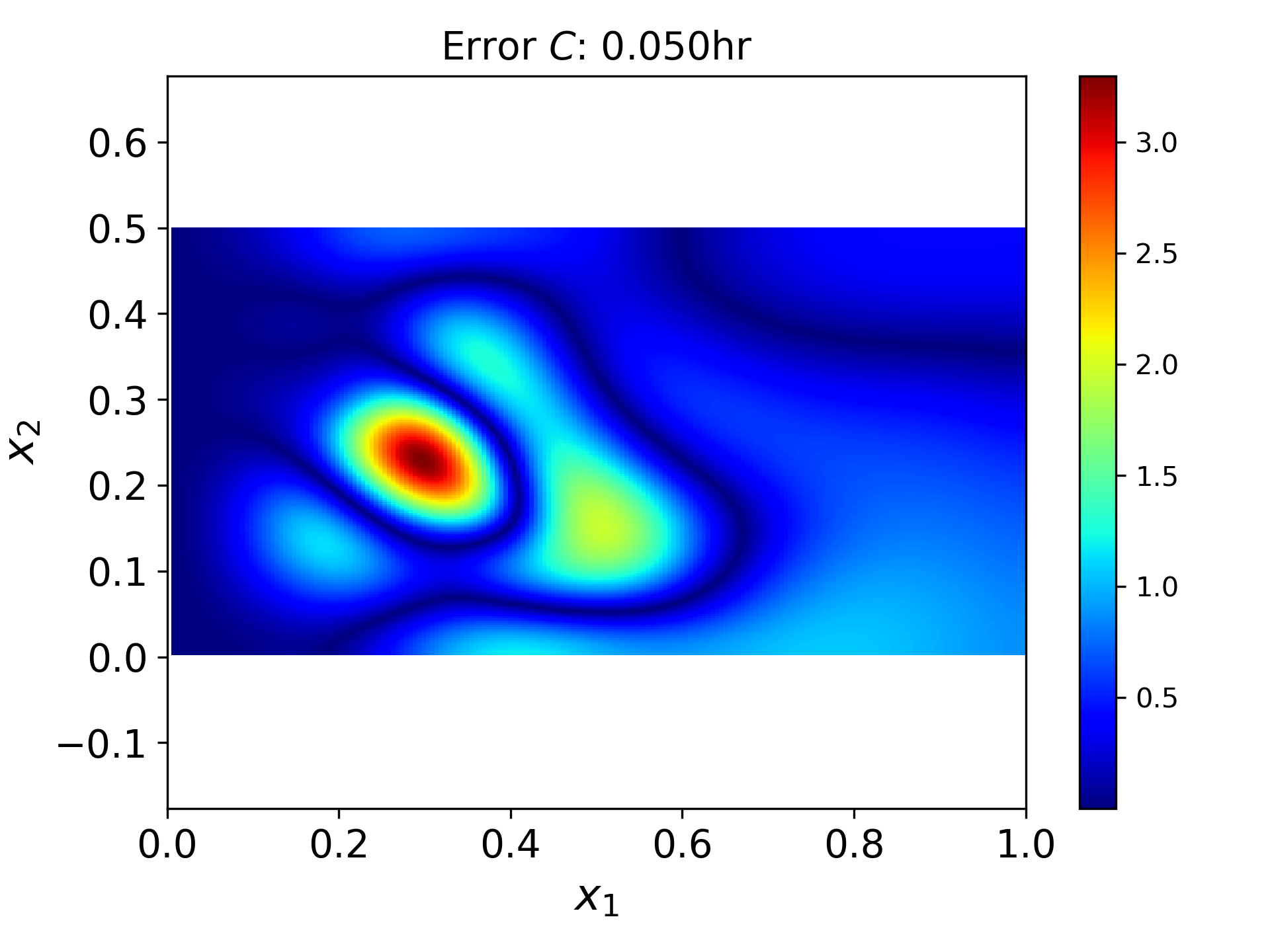}
\end{subfigure}

\begin{subfigure}[h]{0.33\textwidth}
\includegraphics[width=\linewidth]{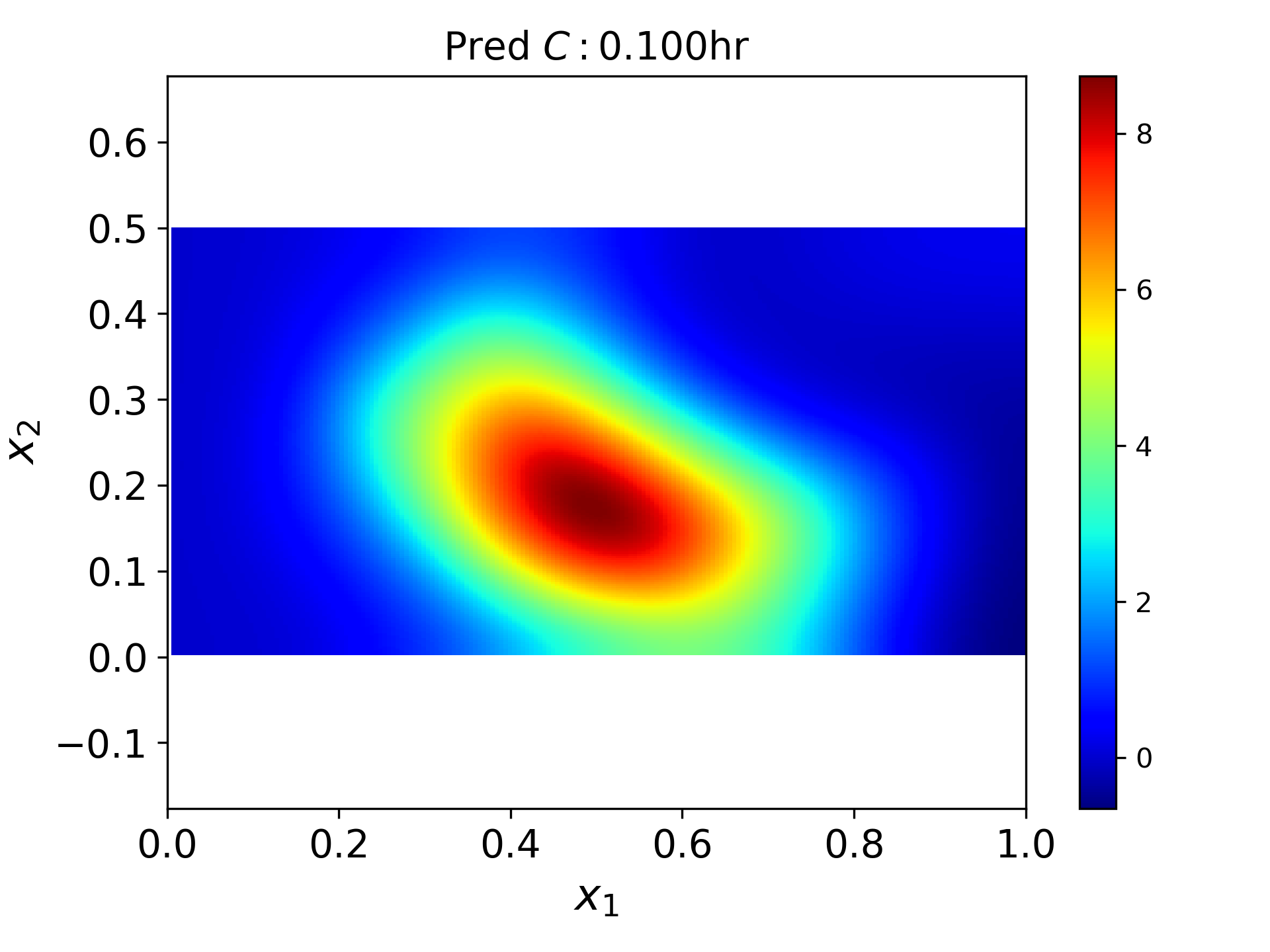}
\end{subfigure}\hfill
\begin{subfigure}[h]{0.33\textwidth}
\includegraphics[width=\linewidth]{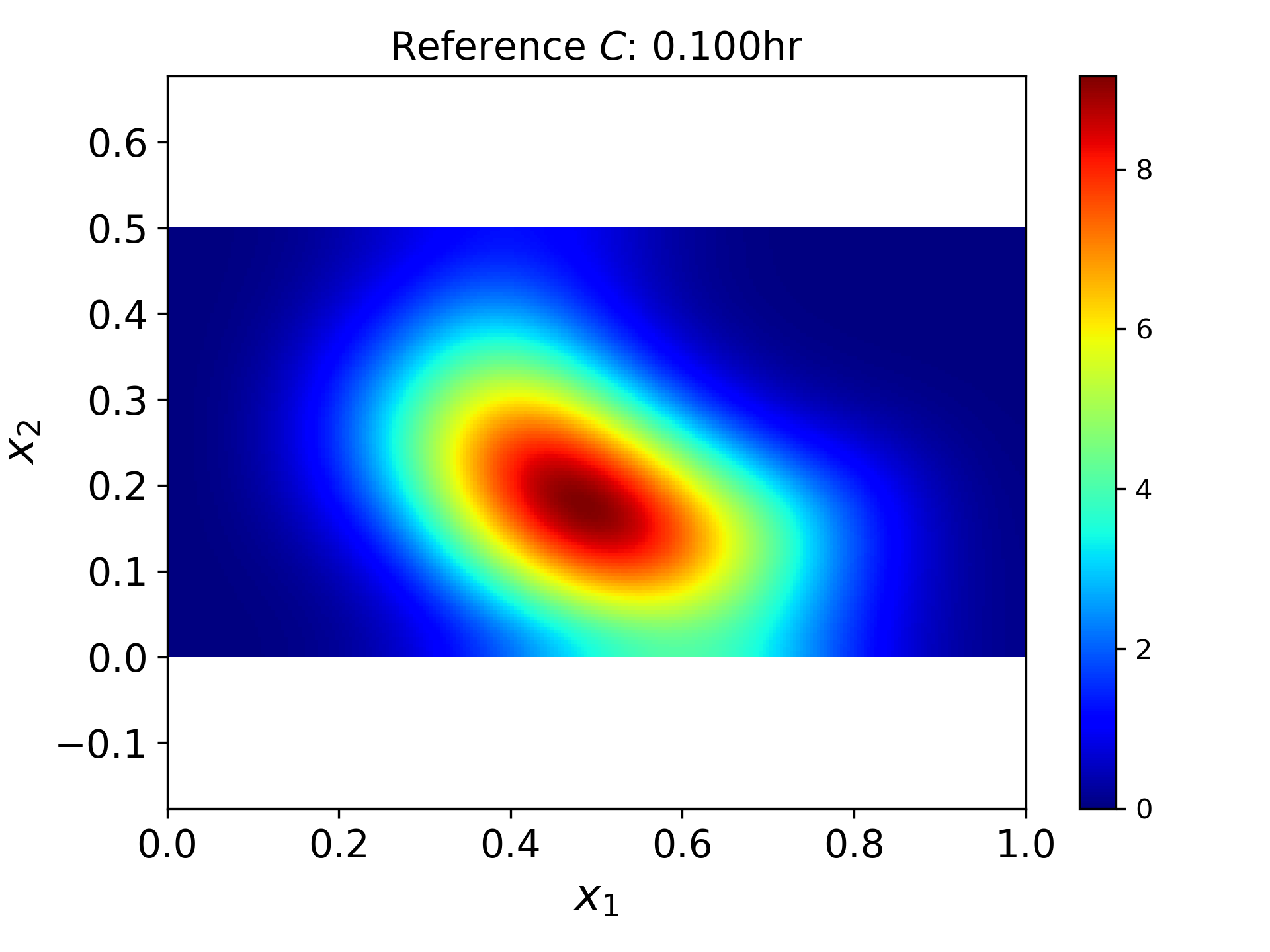}
\end{subfigure}\hfill
\begin{subfigure}[h]{0.33\textwidth}
\includegraphics[width=\linewidth]{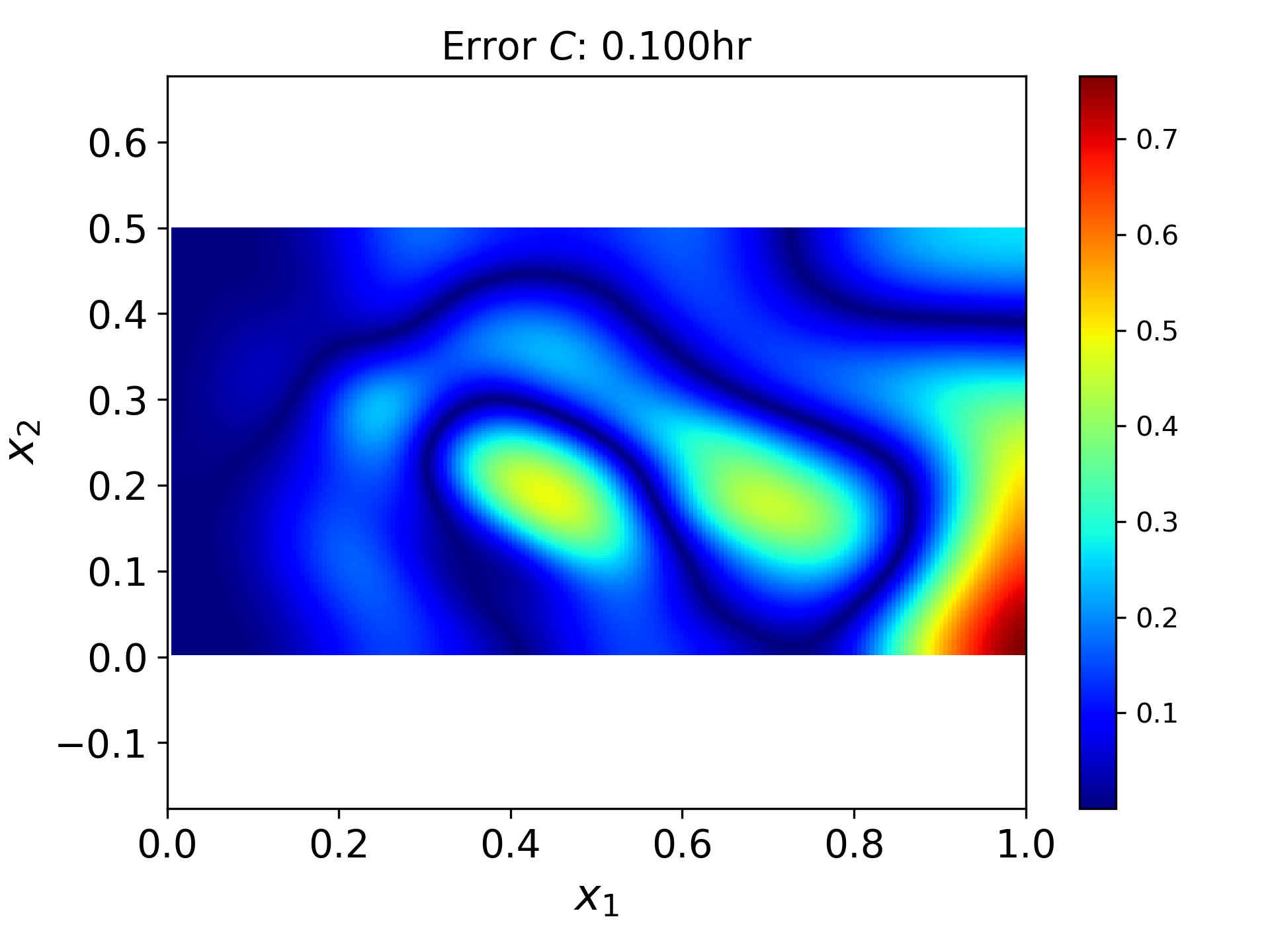}
\end{subfigure}

\begin{subfigure}[h]{0.33\textwidth}
\includegraphics[width=\linewidth]{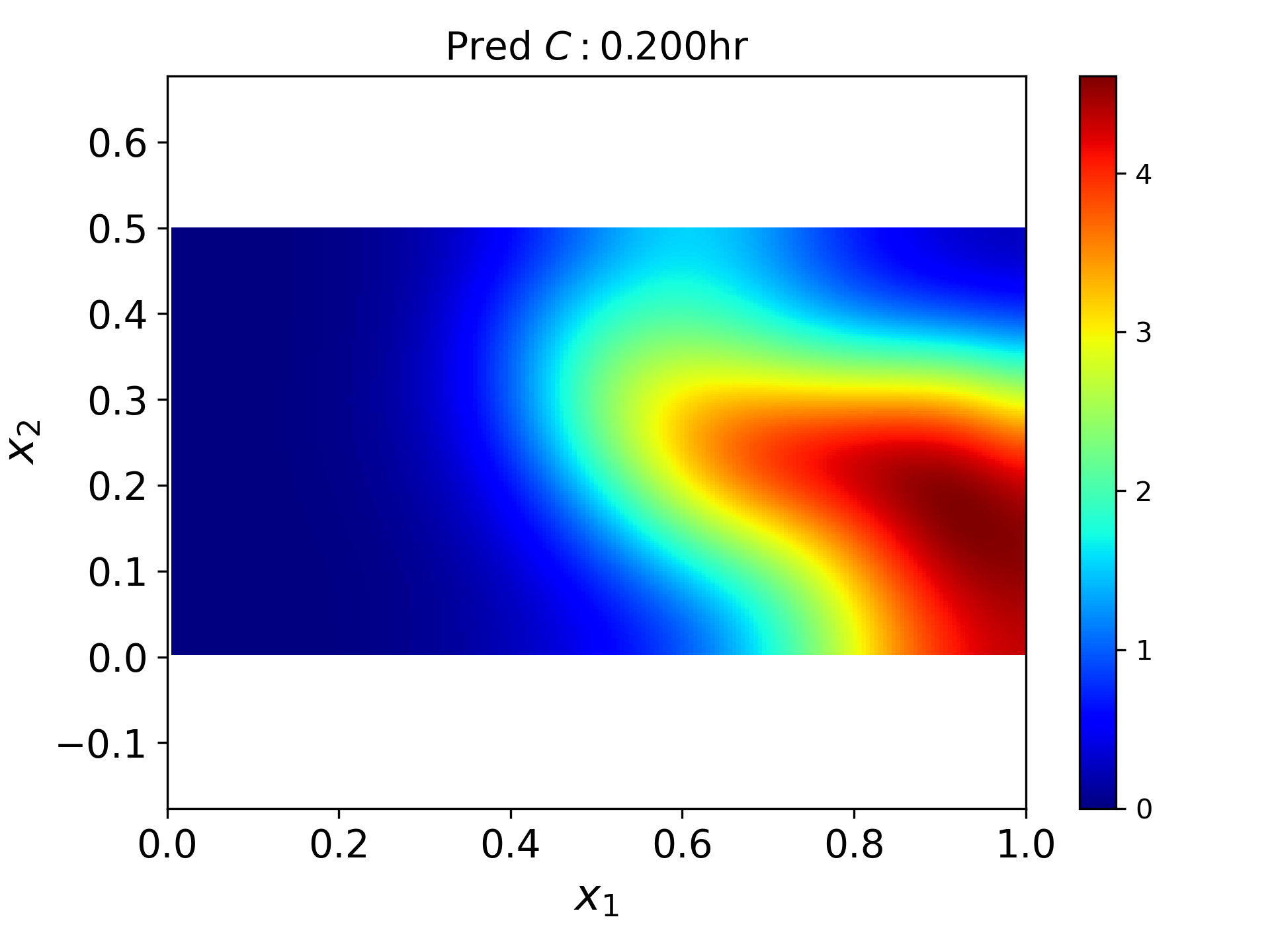}
\end{subfigure}\hfill
\begin{subfigure}[h]{0.33\textwidth}
\includegraphics[width=\linewidth]{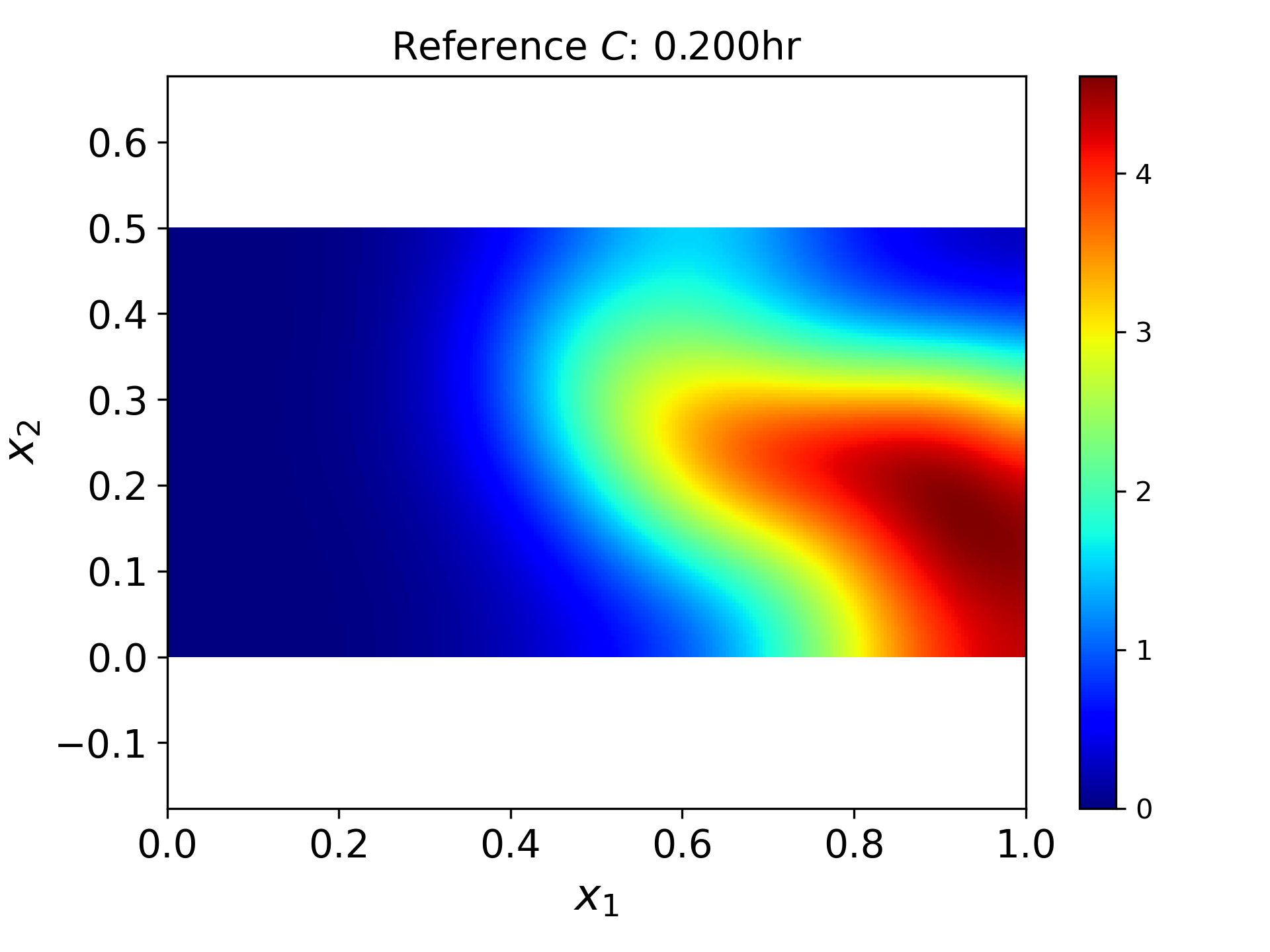}
\end{subfigure}\hfill
\begin{subfigure}[h]{0.33\textwidth}
\includegraphics[width=\linewidth]{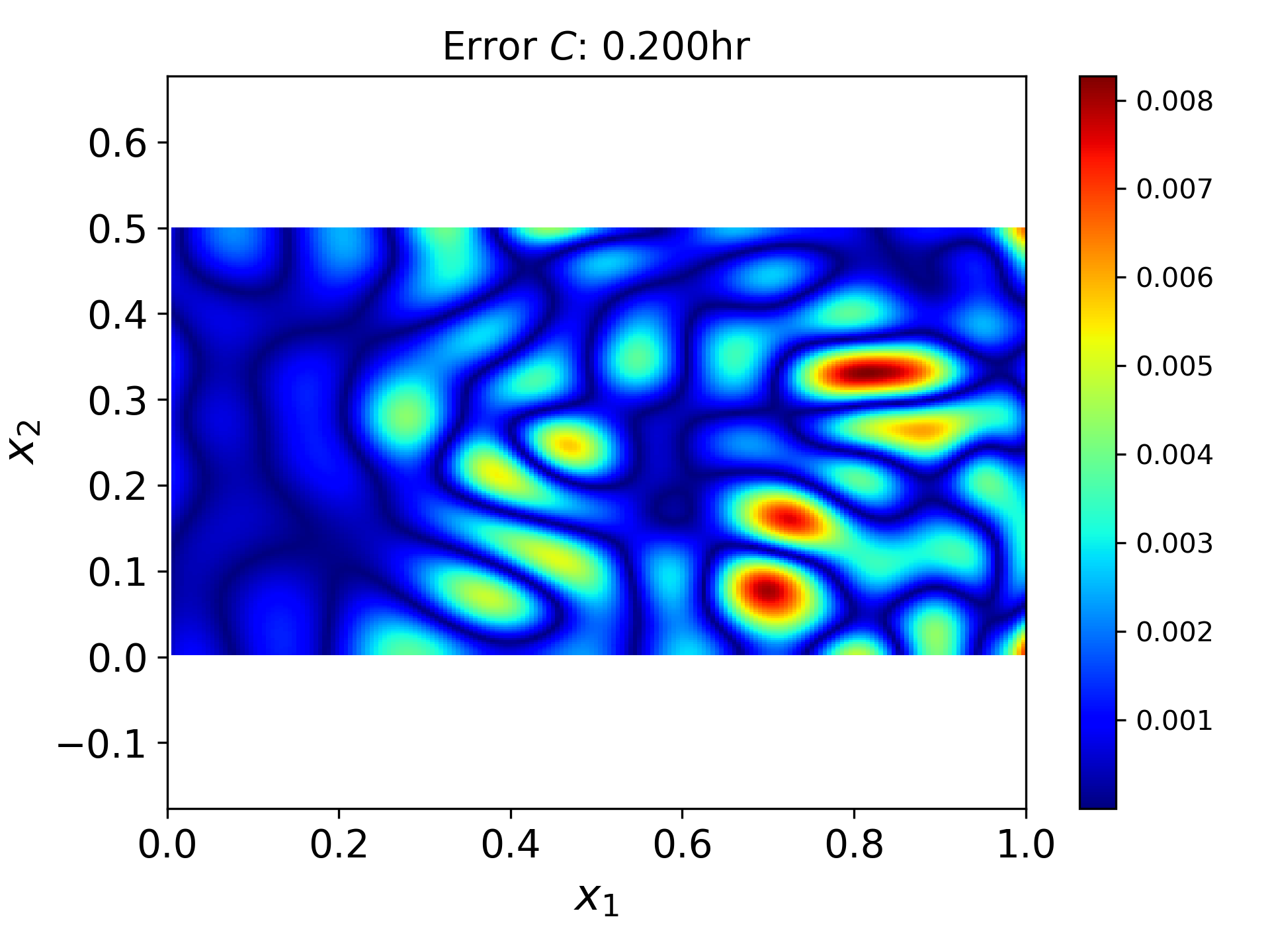}
\end{subfigure}
\caption{The backward PINN solution $\hat{u}$ (the first column), reference $u$ (the second column), and the absolute point error $u-\hat{u}$ (the third column) at $t=0$, 0.05, and 0.1, and at terminal time $T=0.2$. The identity activation function is used in the last layer. No measurements of $u$ are used to obtain the backward solution.}
\label{fig:PINN_backward_identity_0.2}
\end{figure}
\begin{figure}[h]
\begin{subfigure}[h]{0.33\textwidth}
\includegraphics[width=\linewidth]{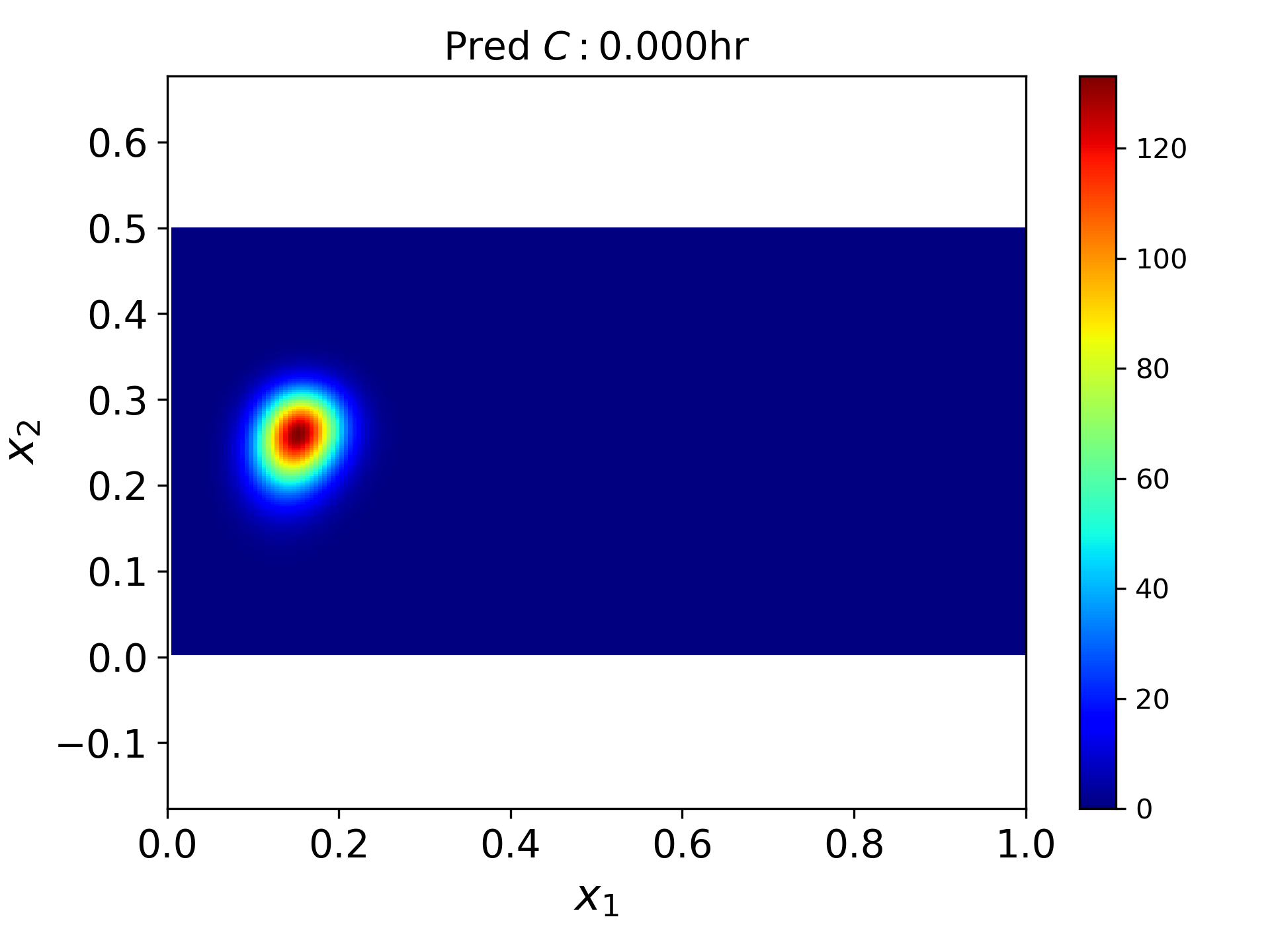}
\end{subfigure}\hfill
\begin{subfigure}[h]{0.33\textwidth}
\includegraphics[width=\linewidth]{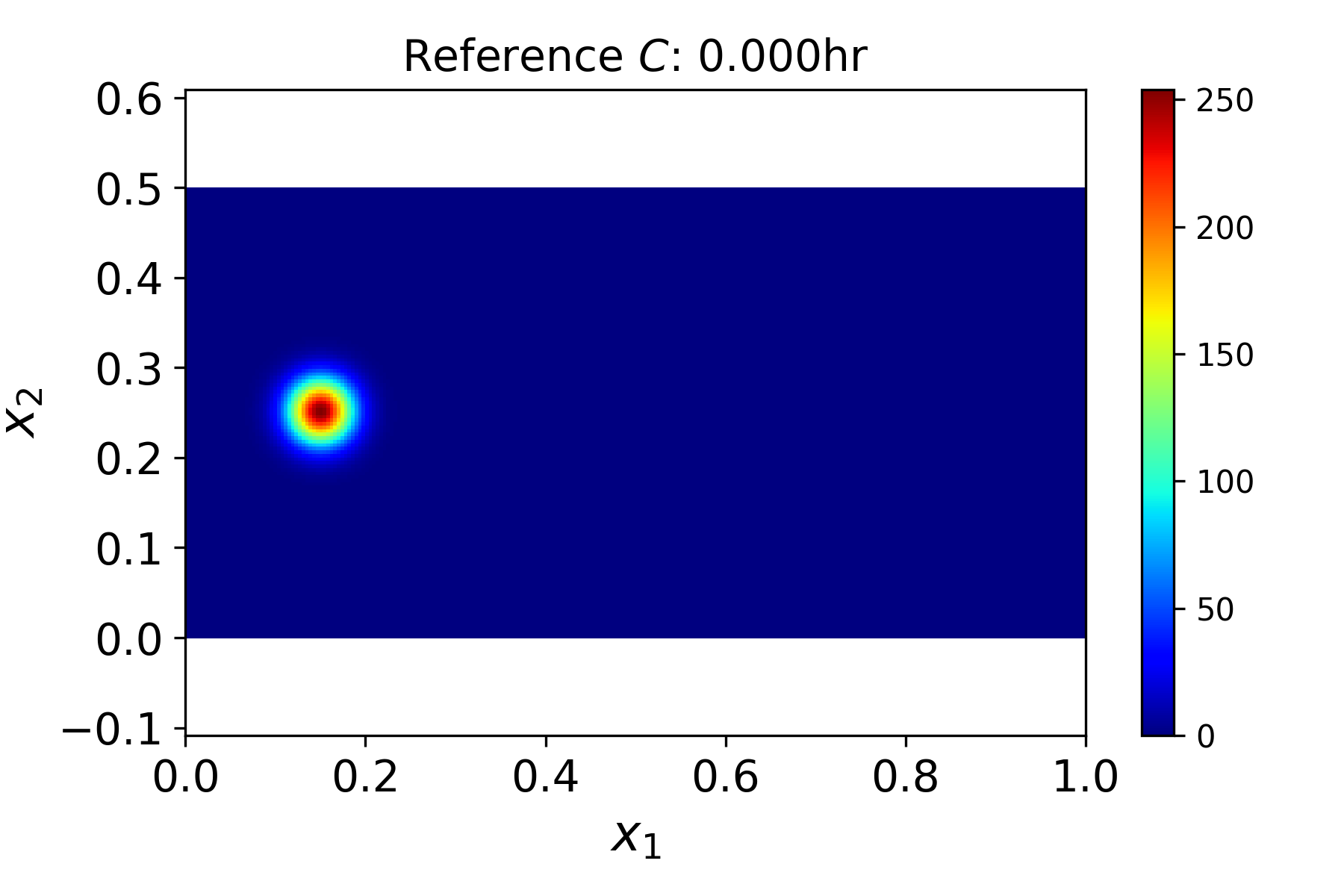}
\end{subfigure}\hfill
\begin{subfigure}[h]{0.33\textwidth}
\includegraphics[width=\linewidth]{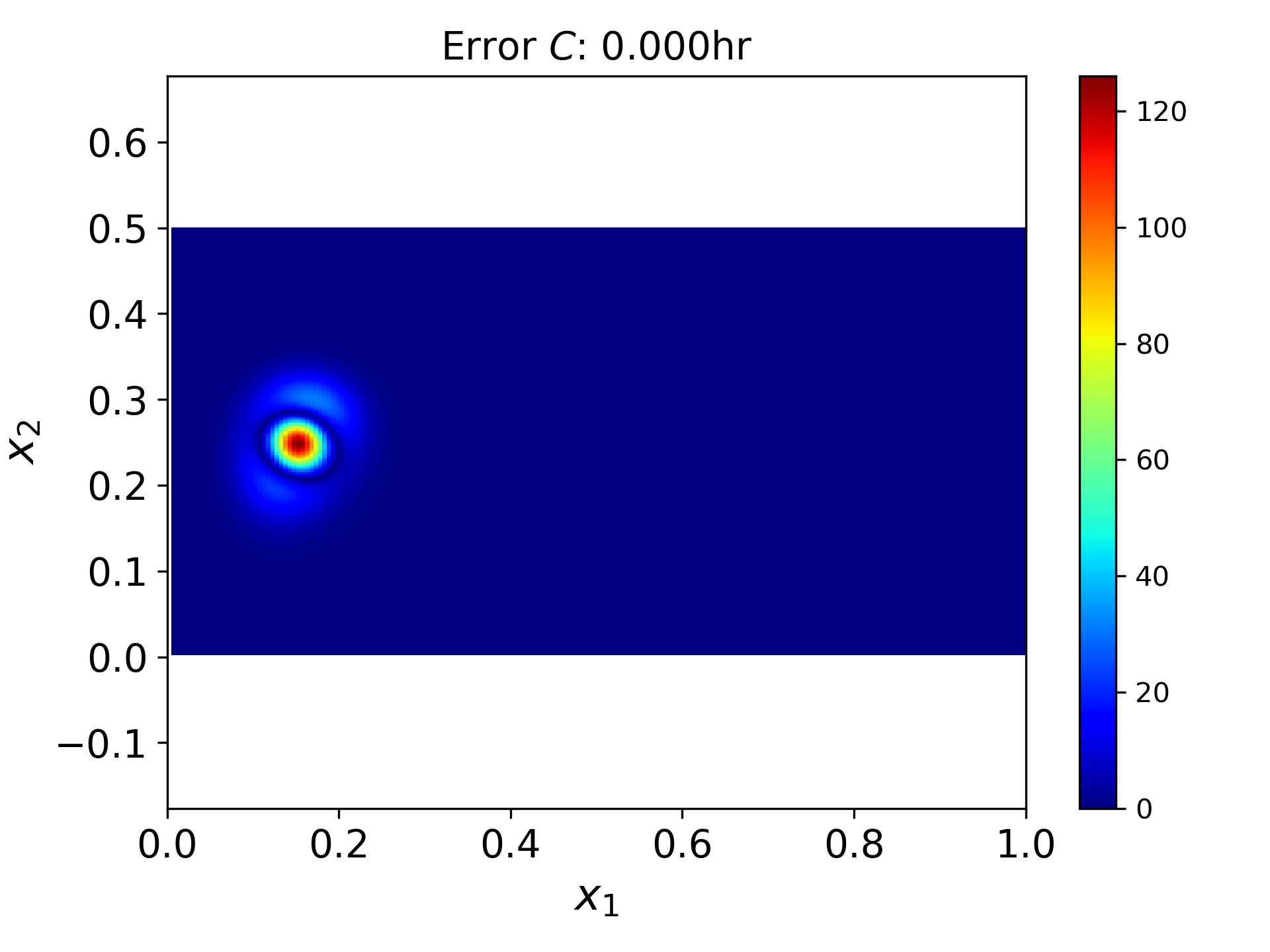}
\end{subfigure}

\begin{subfigure}[h]{0.33\textwidth}
\includegraphics[width=\linewidth]{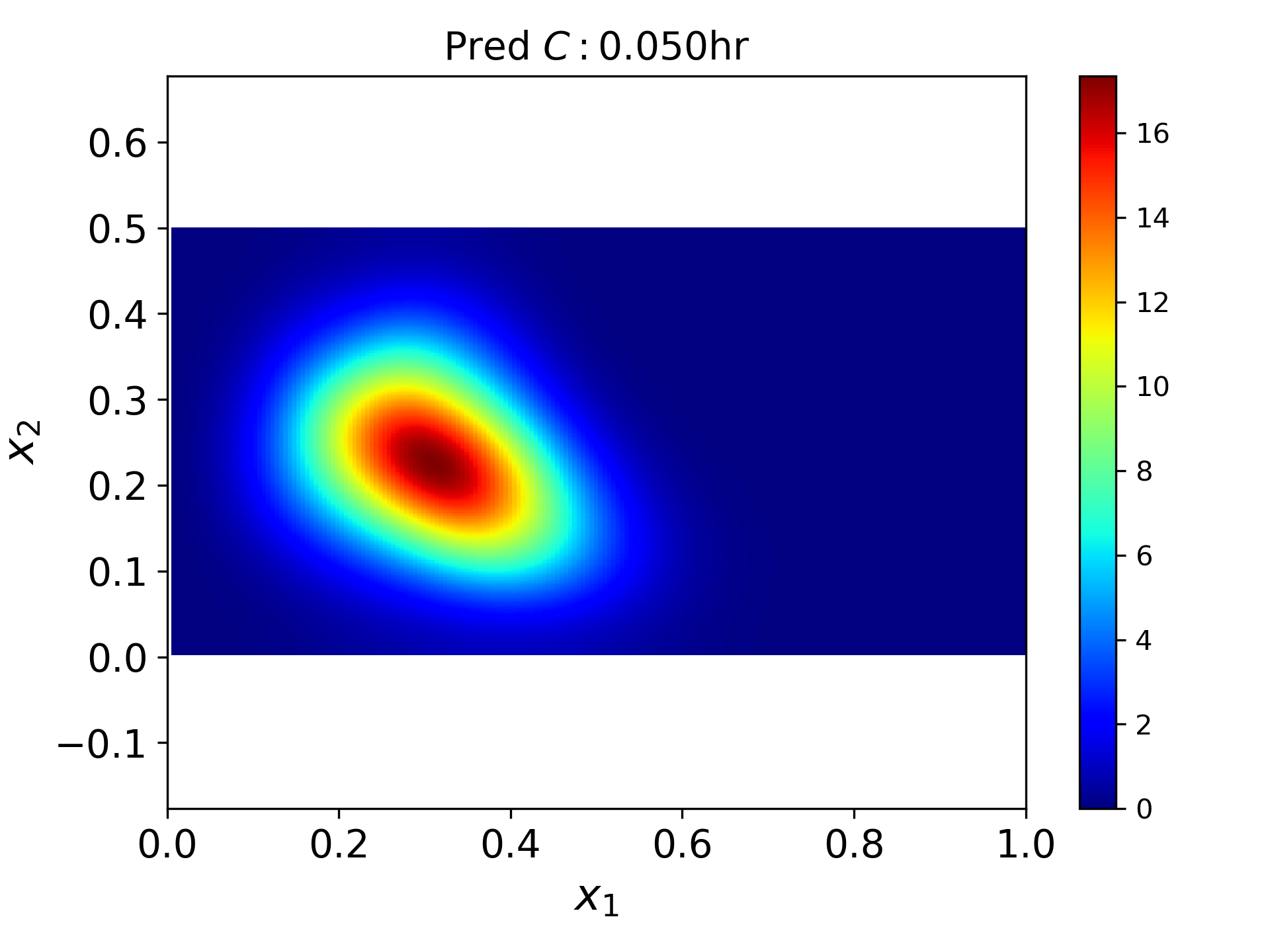}
\end{subfigure}\hfill
\begin{subfigure}[h]{0.33\textwidth}
\includegraphics[width=\linewidth]{draft_PINN_AD/figures/backward_plot/reference_t_0_050hr.png}
\end{subfigure}\hfill
\begin{subfigure}[h]{0.33\textwidth}
\includegraphics[width=\linewidth]{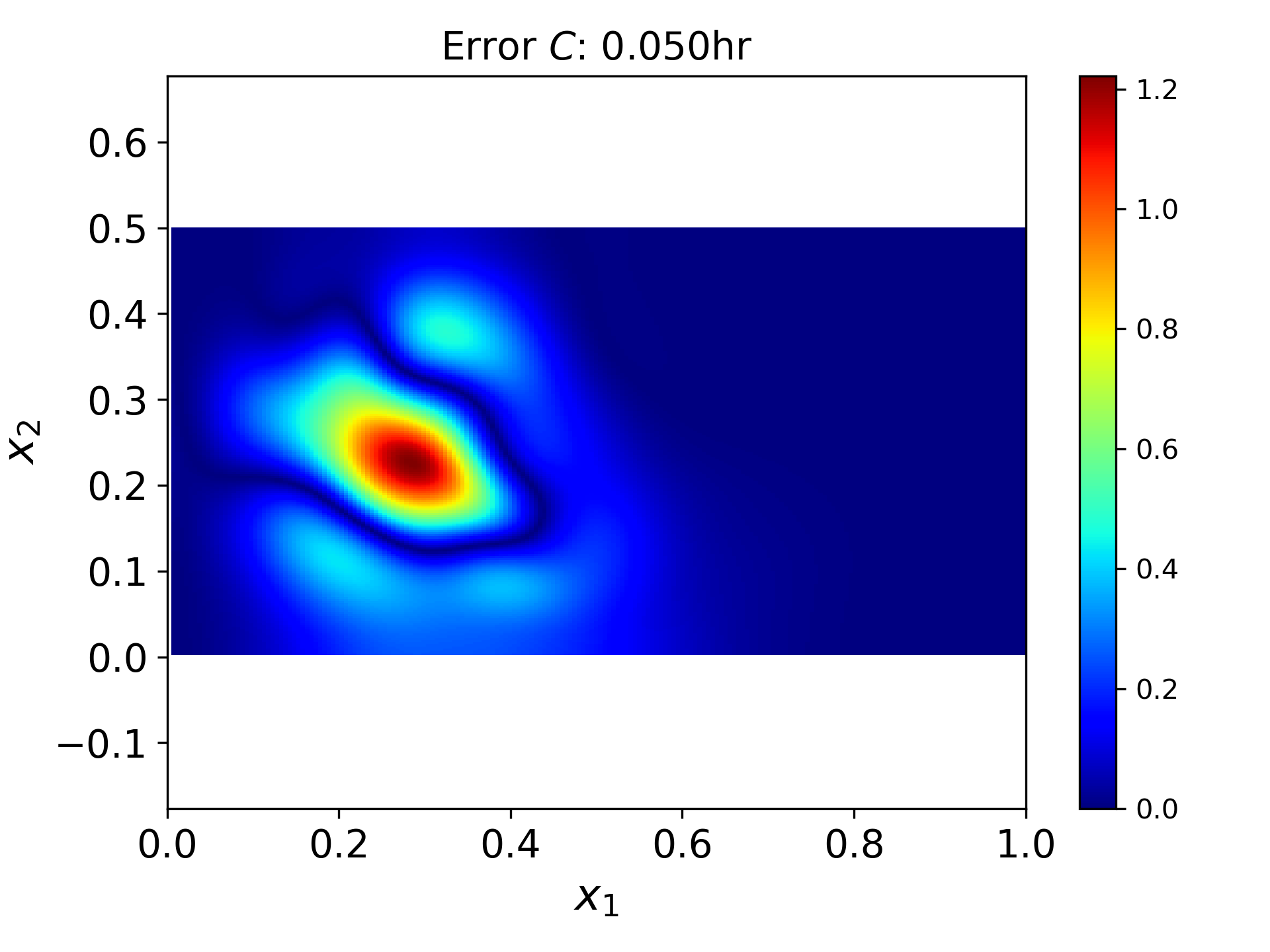}
\end{subfigure}

\begin{subfigure}[h]{0.33\textwidth}
\includegraphics[width=\linewidth]{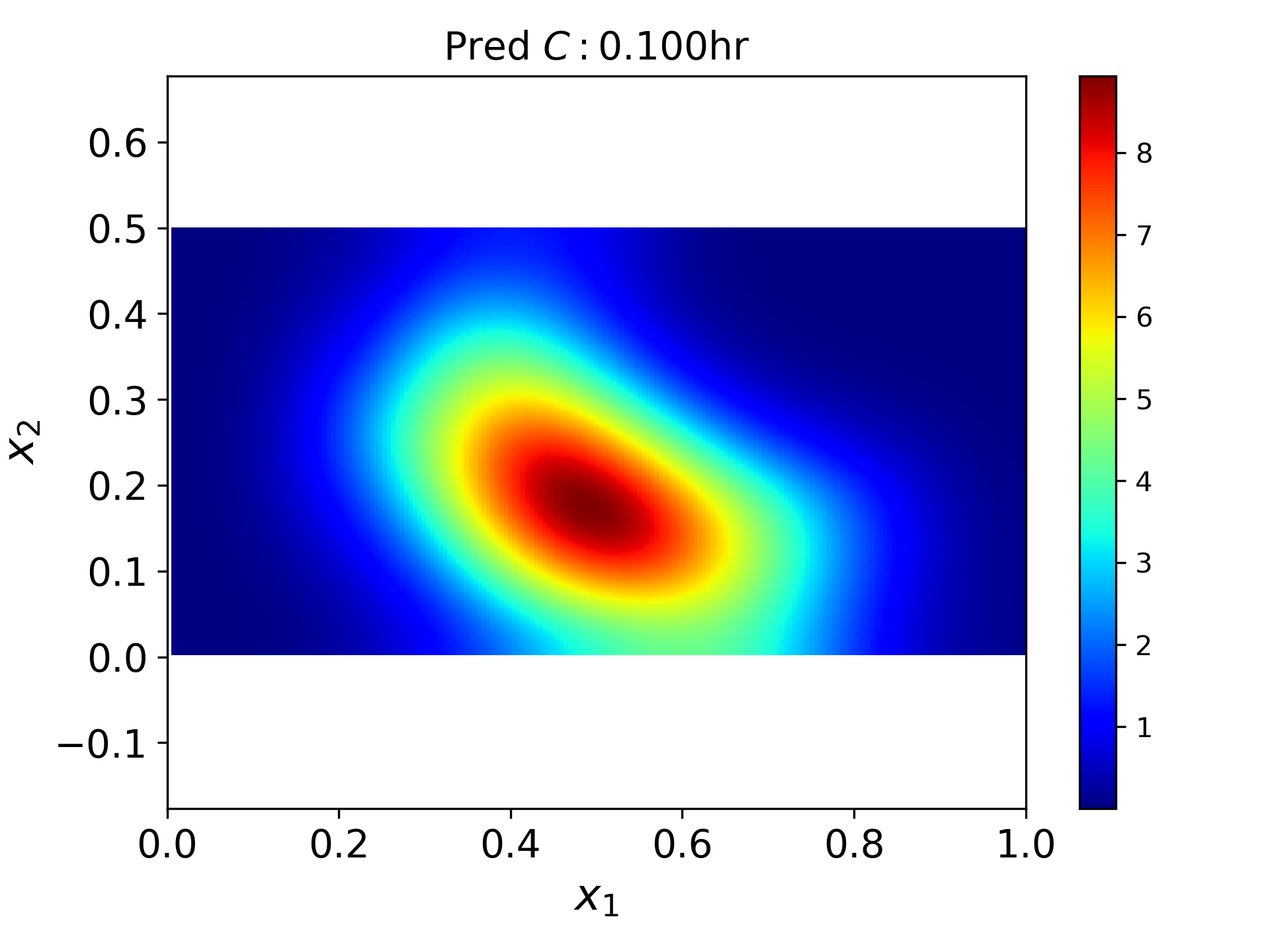}
\end{subfigure}\hfill
\begin{subfigure}[h]{0.33\textwidth}
\includegraphics[width=\linewidth]{draft_PINN_AD/figures/backward_plot/reference_t_0_100hr.png}
\end{subfigure}\hfill
\begin{subfigure}[h]{0.33\textwidth}
\includegraphics[width=\linewidth]{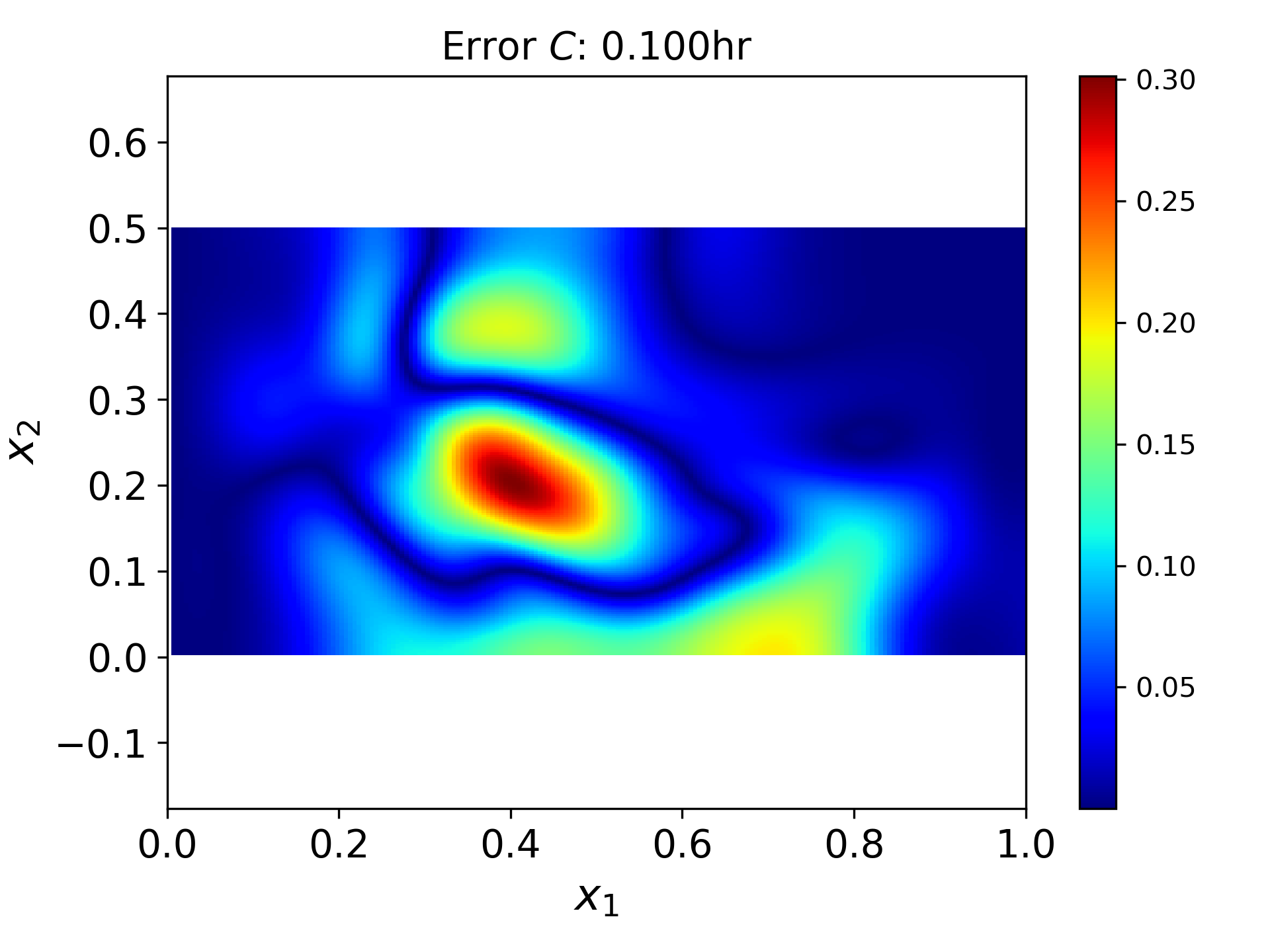}
\end{subfigure}

\begin{subfigure}[h]{0.33\textwidth}
\includegraphics[width=\linewidth]{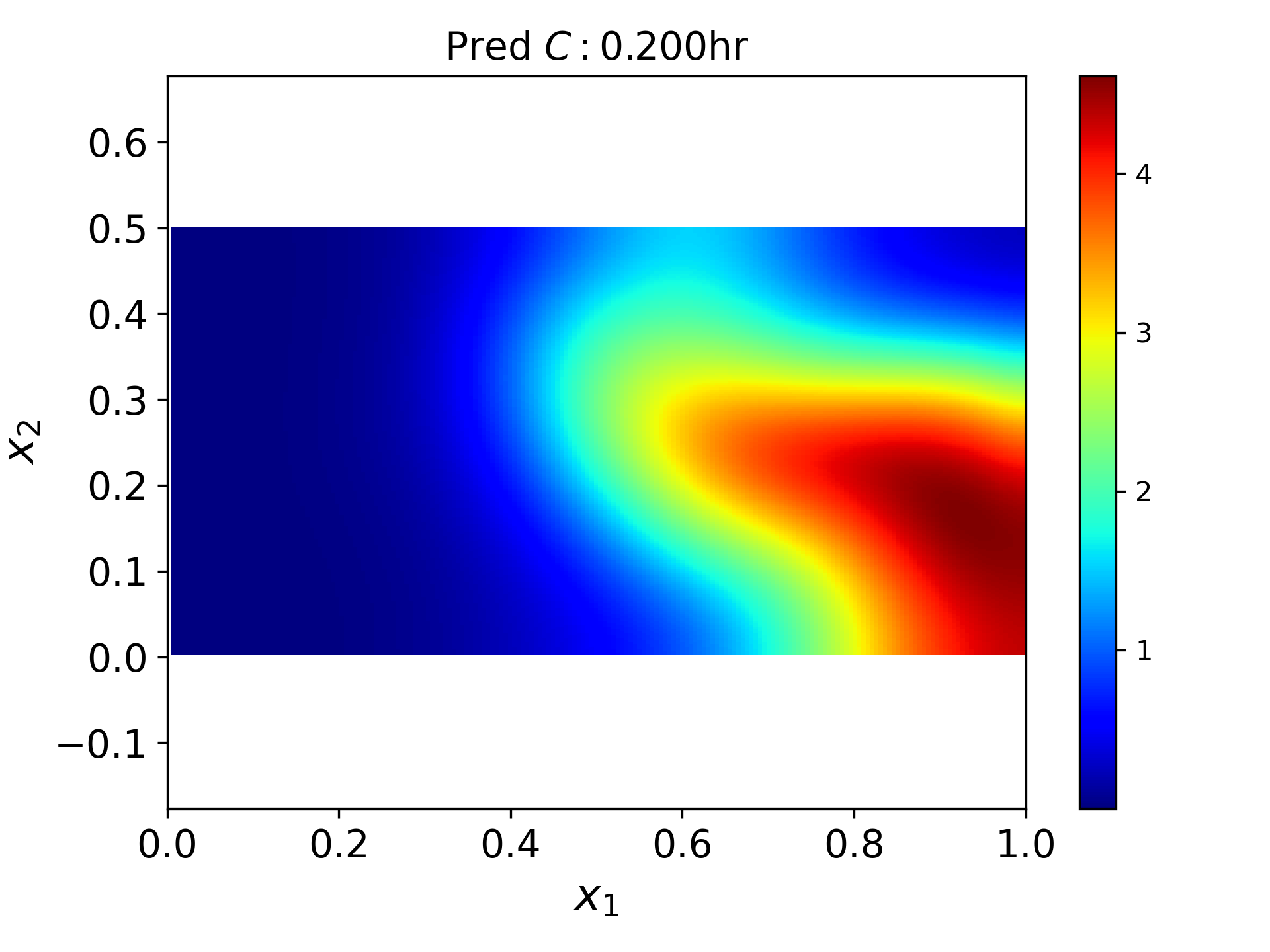}
\end{subfigure}\hfill
\begin{subfigure}[h]{0.33\textwidth}
\includegraphics[width=\linewidth]{draft_PINN_AD/figures/backward_plot/reference_t_0_200hr.png}
\end{subfigure}\hfill
\begin{subfigure}[h]{0.33\textwidth}
\includegraphics[width=\linewidth]{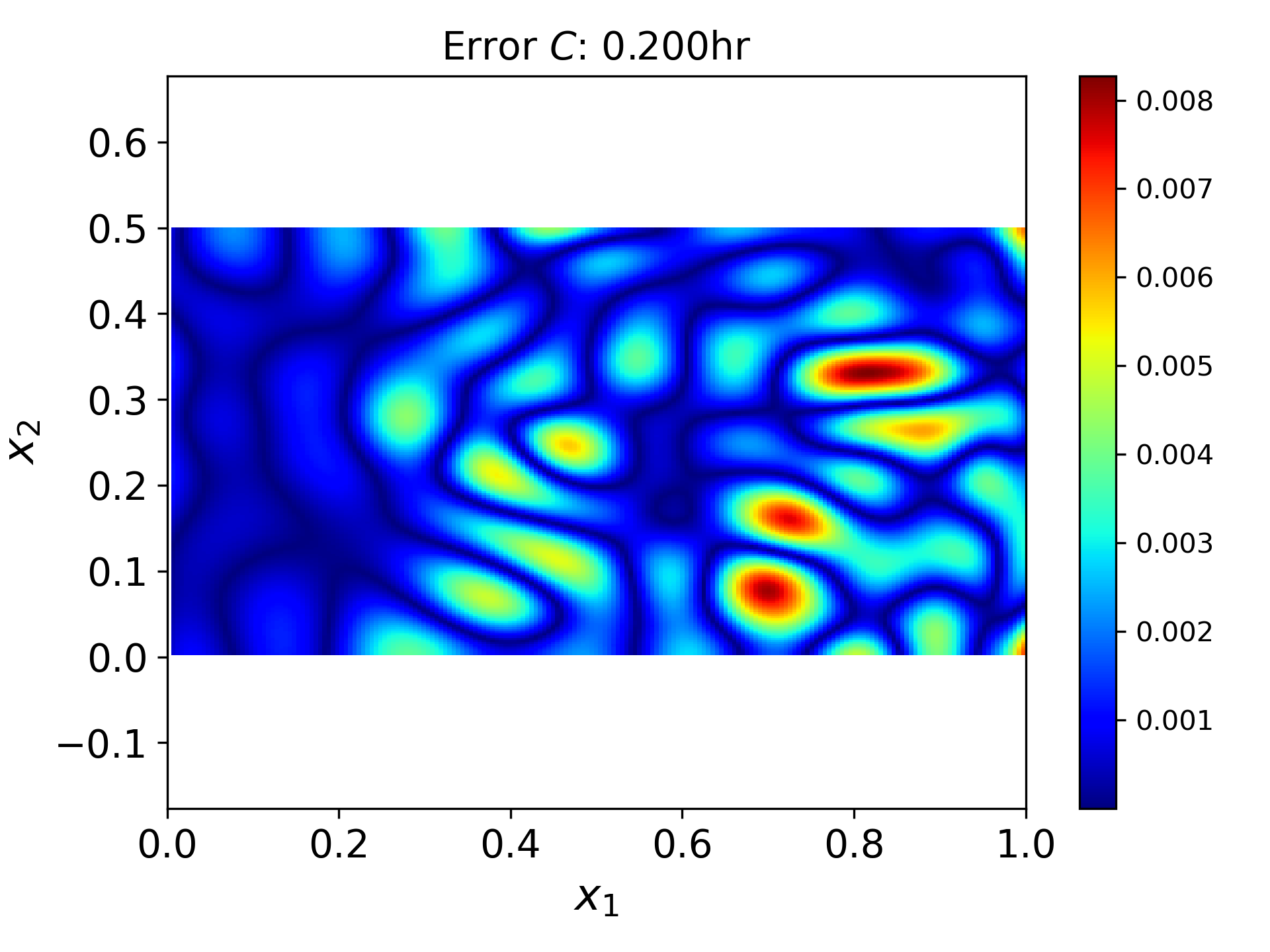}
\end{subfigure}
\caption{
The backward PINN solution $\hat{u}$ (the first column), reference $u$ (the second column), and the absolute point error $u-\hat{u}$ (the third column) at $t=0$, 0.05, and 0.1, and at terminal time $T=0.2$. The sigmoid activation function is used in the last layer. No measurements of $u$ are used to obtain the backward solution.}
\label{fig:PINN_backward_sigmoid_0.2}
\end{figure}
\begin{figure}[h]
\begin{subfigure}[h]{0.33\textwidth}
\includegraphics[width=\linewidth]{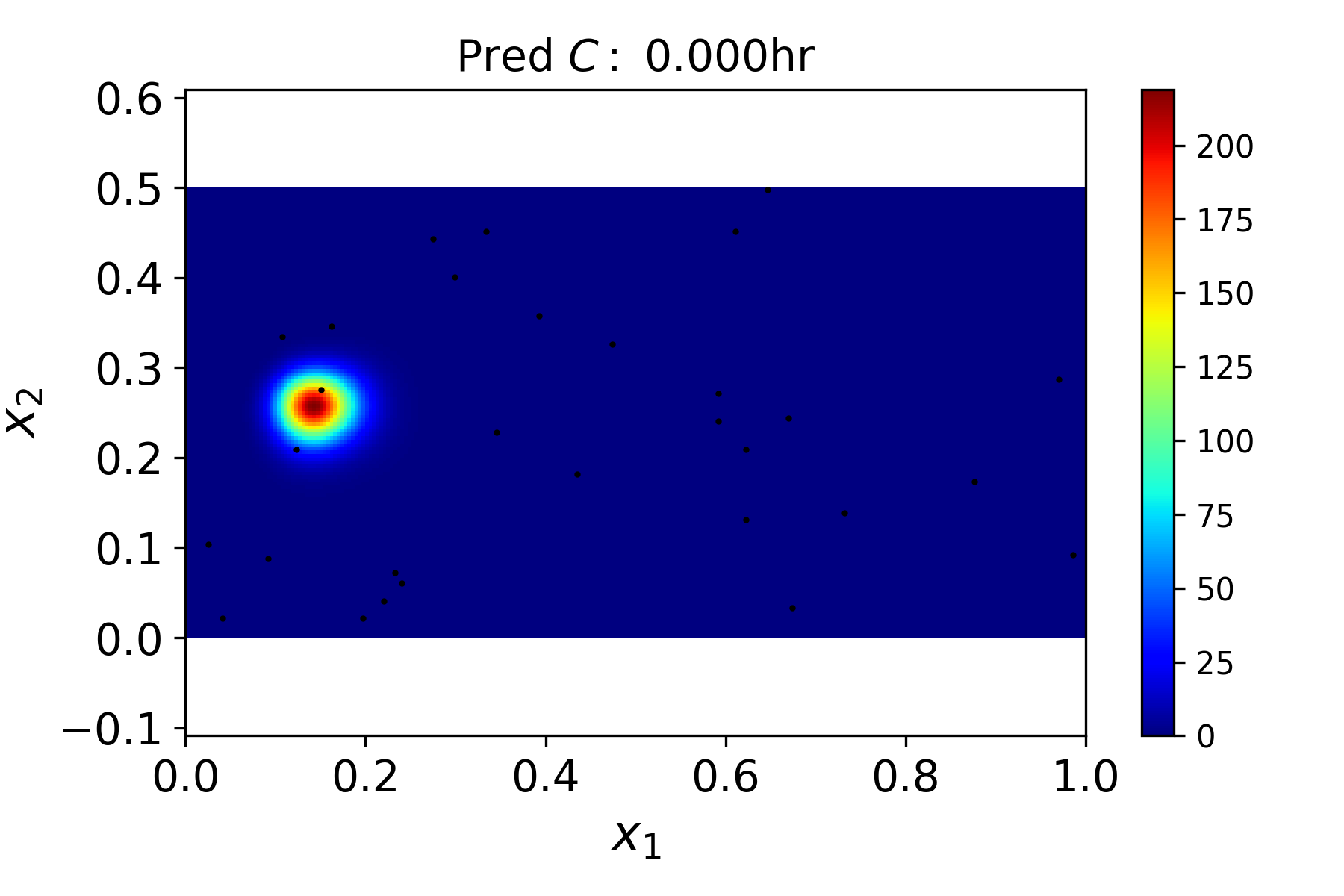}
\end{subfigure}\hfill
\begin{subfigure}[h]{0.33\textwidth}
\includegraphics[width=\linewidth]{draft_PINN_AD/figures/backward_plot/reference_t_0_000hr.png}
\end{subfigure}\hfill
\begin{subfigure}[h]{0.33\textwidth}
\includegraphics[width=\linewidth]{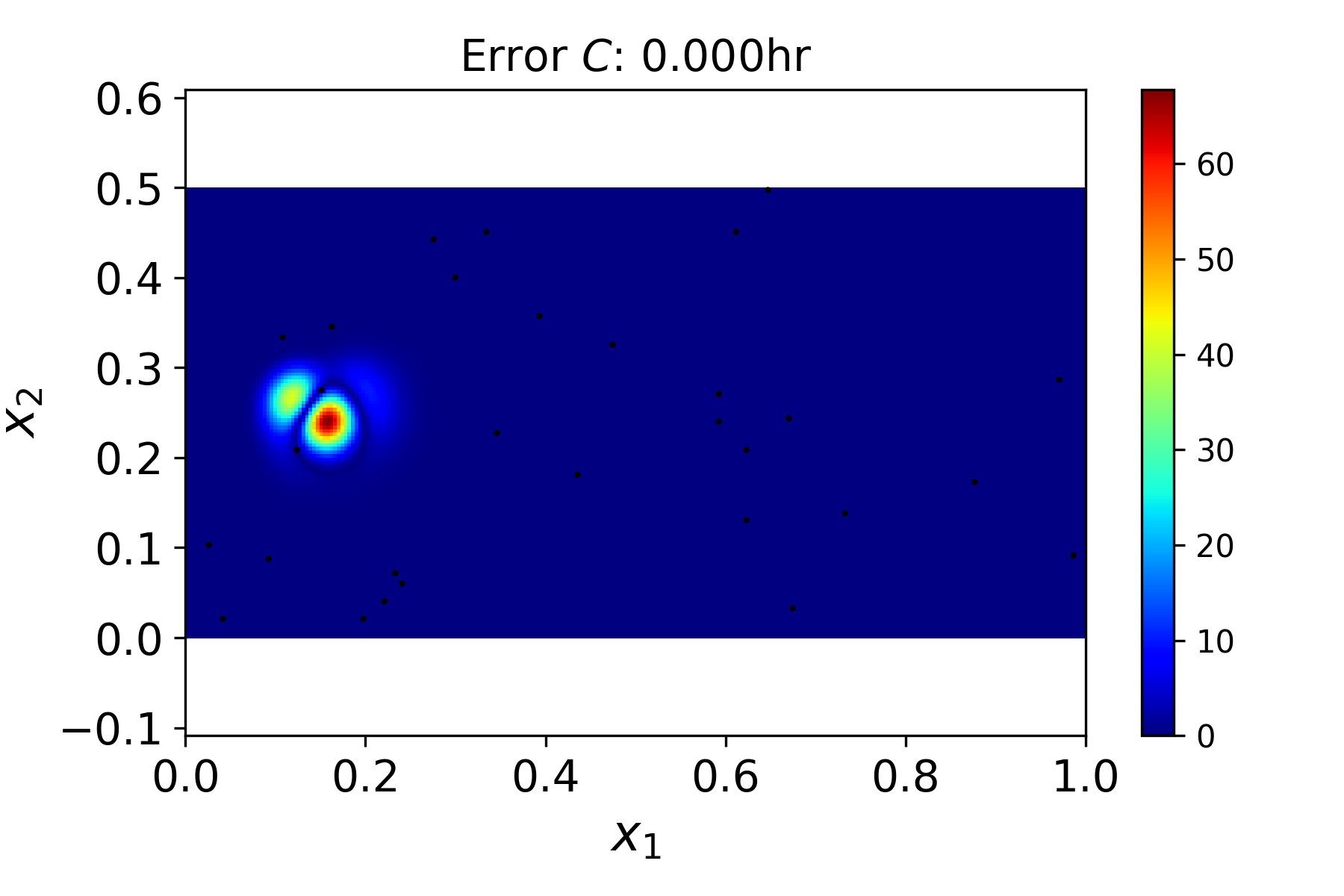}
\end{subfigure}

\begin{subfigure}[h]{0.33\textwidth}
\includegraphics[width=\linewidth]{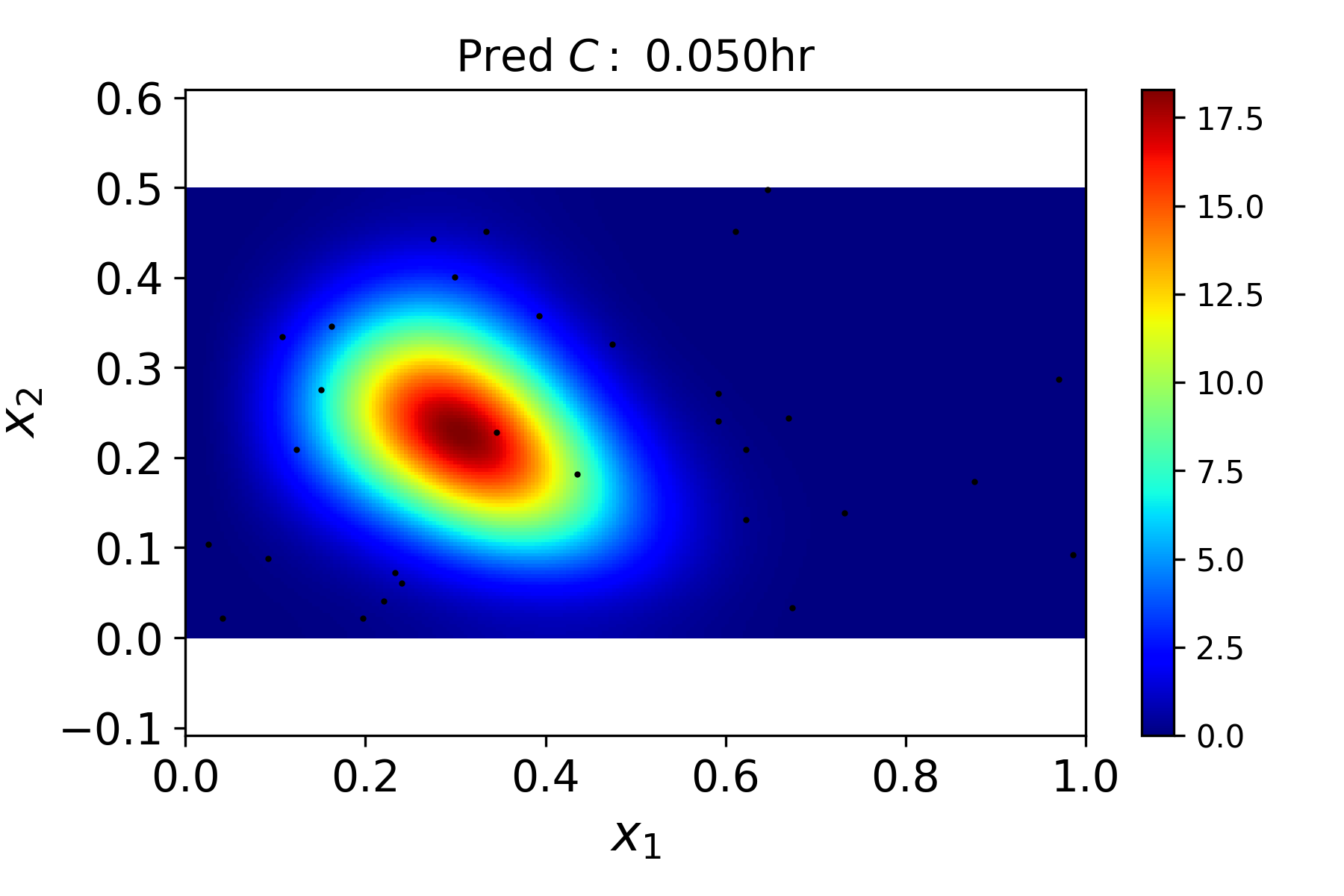}
\end{subfigure}\hfill
\begin{subfigure}[h]{0.33\textwidth}
\includegraphics[width=\linewidth]{draft_PINN_AD/figures/backward_plot/reference_t_0_050hr.png}
\end{subfigure}\hfill
\begin{subfigure}[h]{0.33\textwidth}
\includegraphics[width=\linewidth]{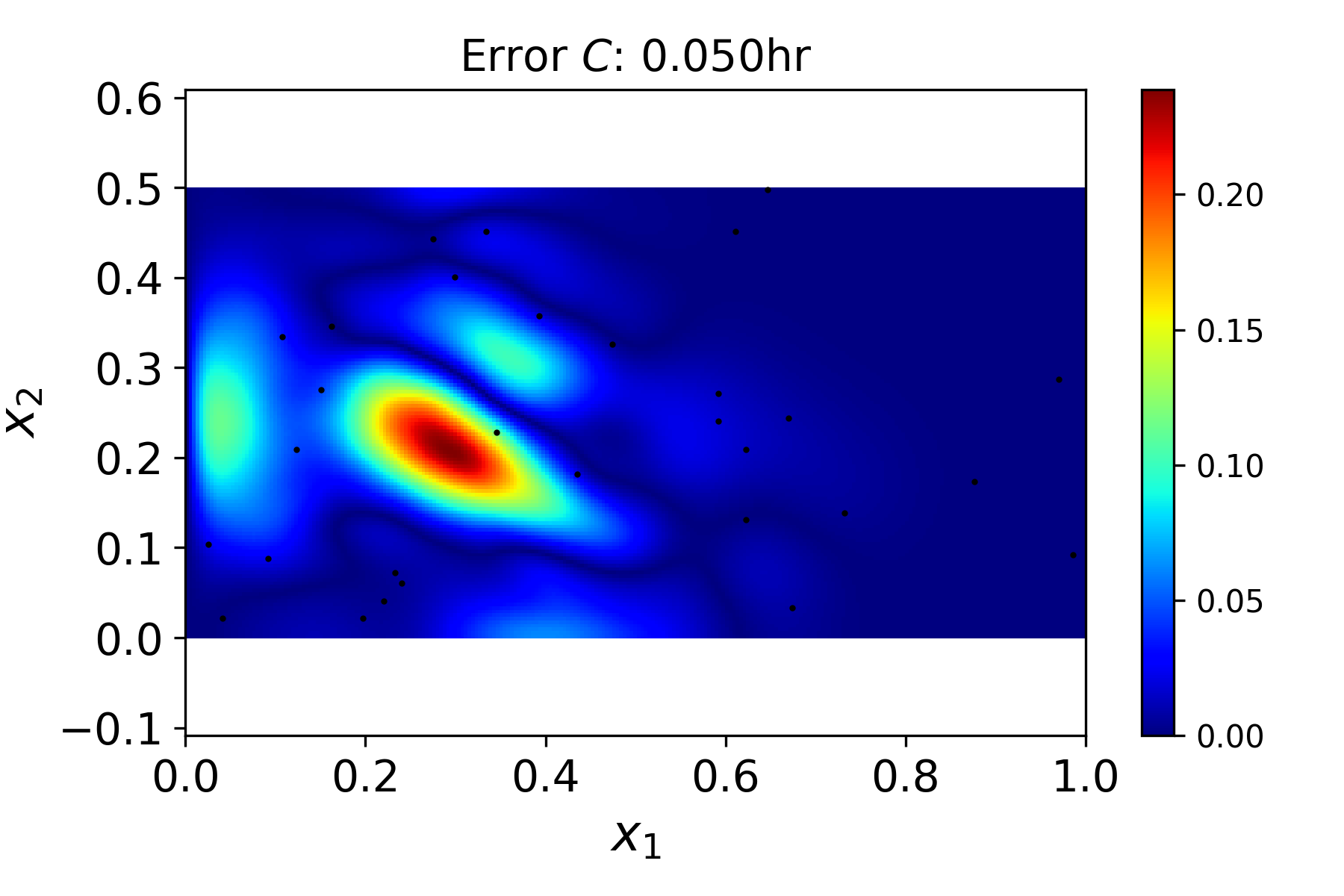}
\end{subfigure}

\begin{subfigure}[h]{0.33\textwidth}
\includegraphics[width=\linewidth]{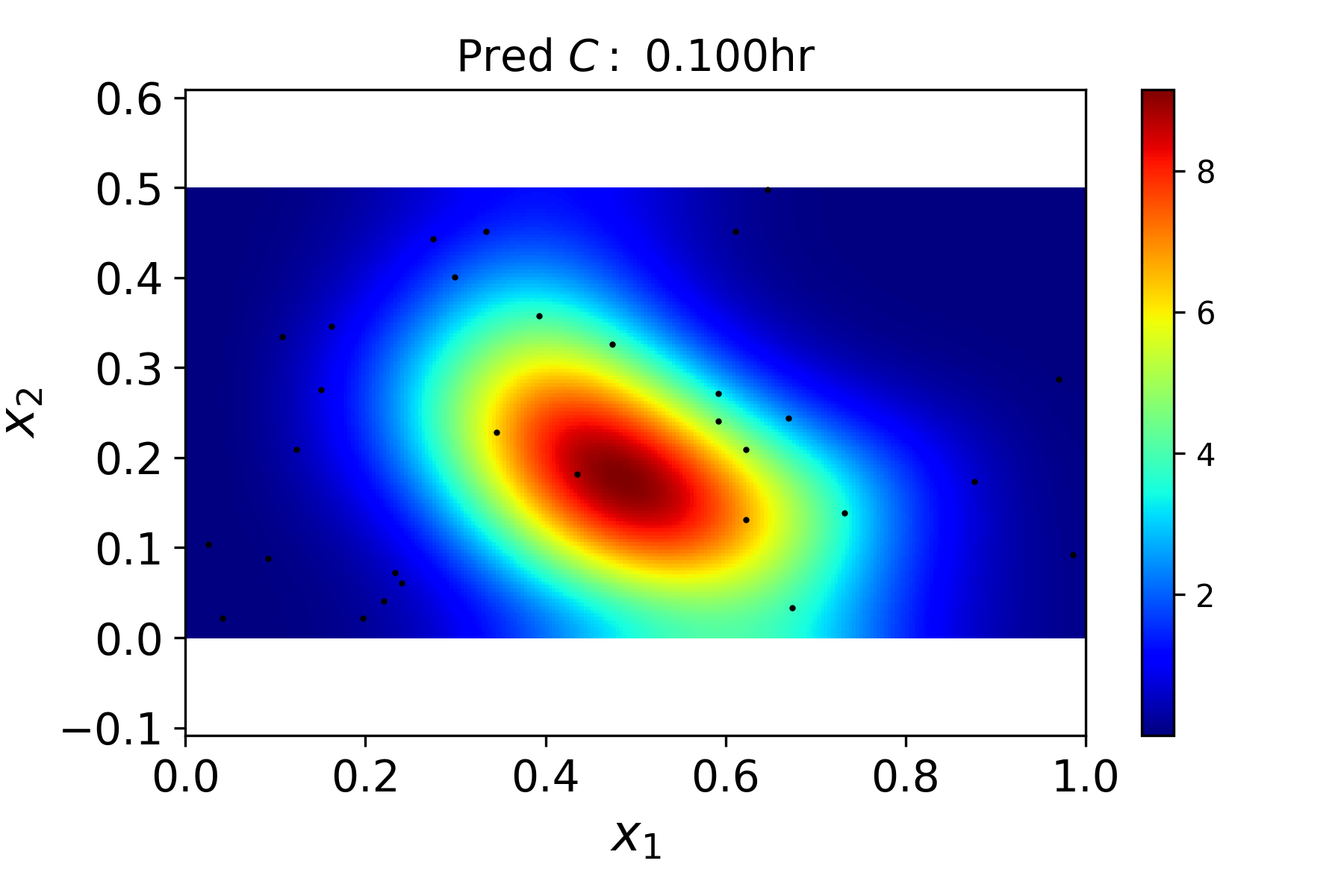}
\end{subfigure}\hfill
\begin{subfigure}[h]{0.33\textwidth}
\includegraphics[width=\linewidth]{draft_PINN_AD/figures/backward_plot/reference_t_0_100hr.png}
\end{subfigure}\hfill
\begin{subfigure}[h]{0.33\textwidth}
\includegraphics[width=\linewidth]{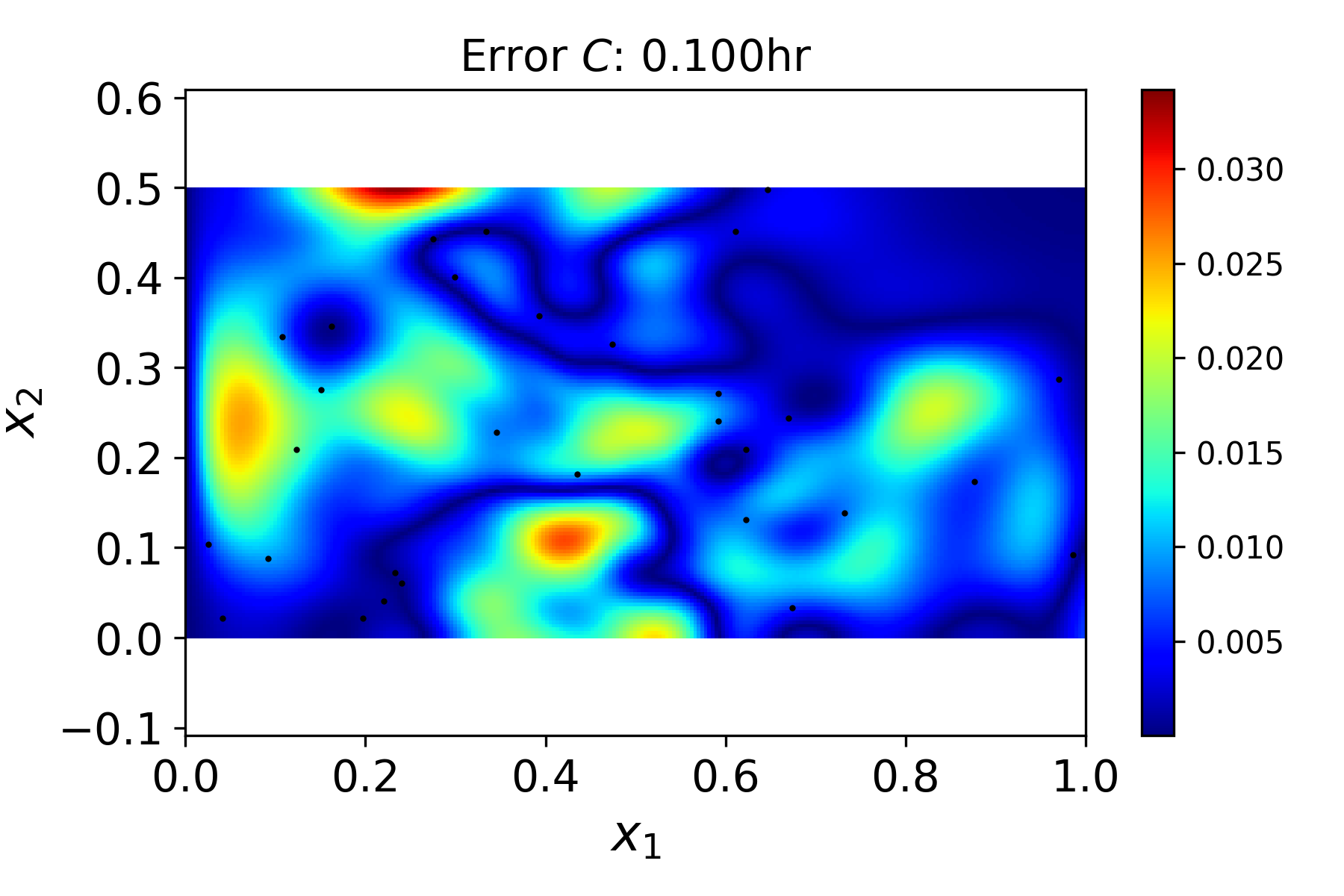}
\end{subfigure}

\begin{subfigure}[h]{0.33\textwidth}
\includegraphics[width=\linewidth]{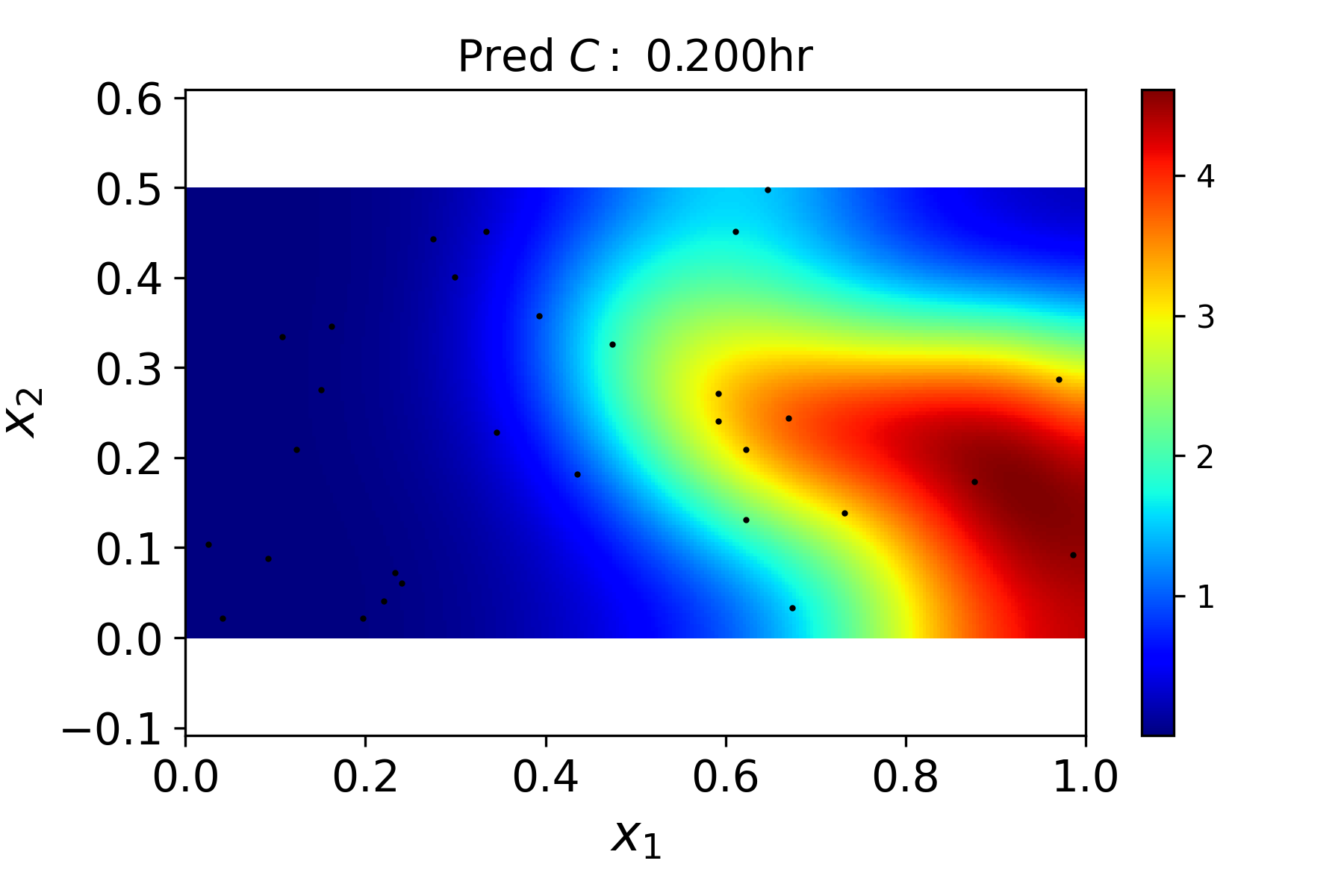}
\end{subfigure}\hfill
\begin{subfigure}[h]{0.33\textwidth}
\includegraphics[width=\linewidth]{draft_PINN_AD/figures/backward_plot/reference_t_0_200hr.png}
\end{subfigure}\hfill
\begin{subfigure}[h]{0.33\textwidth}
\includegraphics[width=\linewidth]{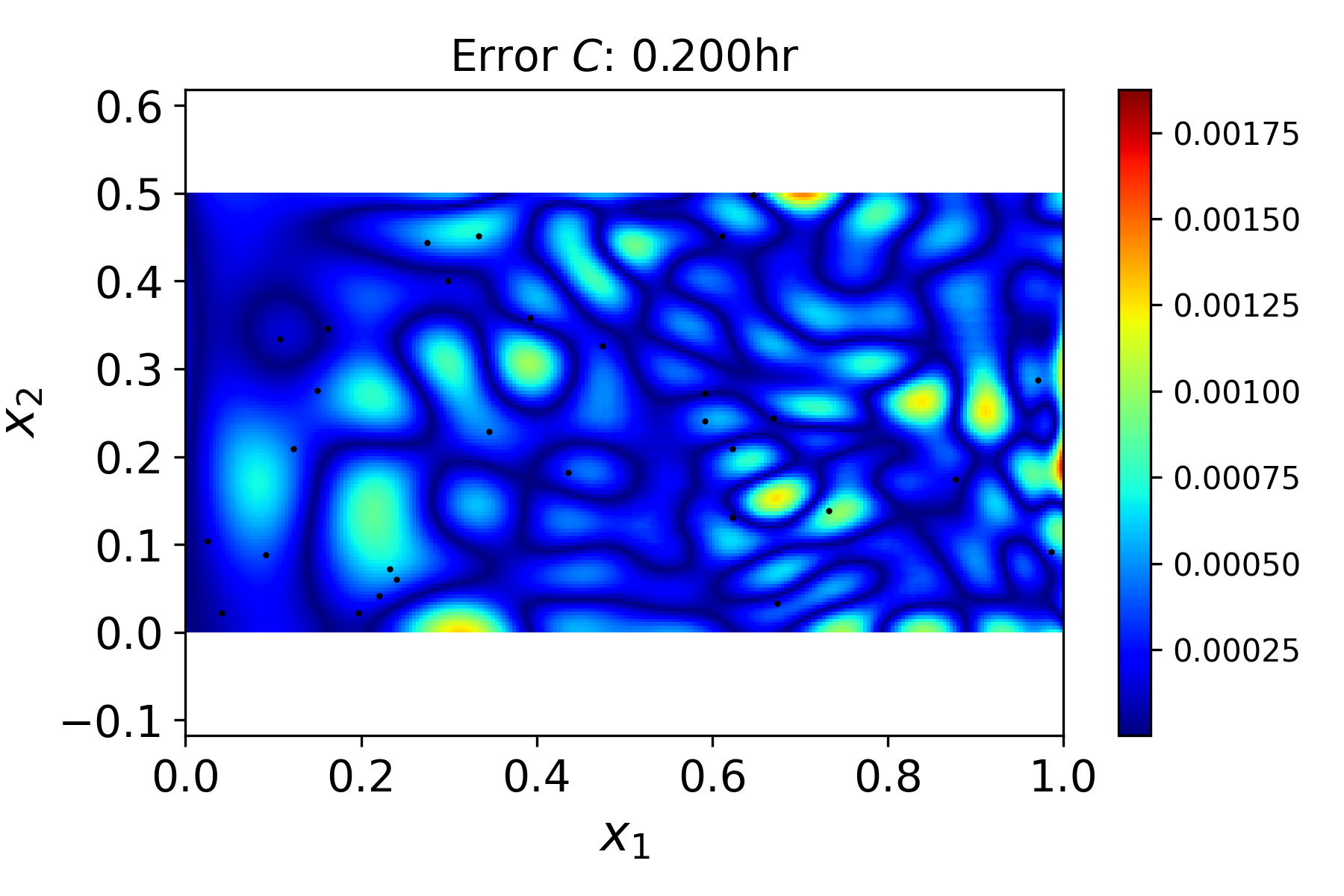}
\end{subfigure}
\caption{The backward PINN solution $\hat{u}$ (the first column), reference $u$ (the second column), and the absolute point error $u-\hat{u}$ (the third column) at $t=0$, 0.05, and 0.1, and at terminal time $T=0.2$. The sigmoid activation function is used in the last layer. $N_x = 30$ measurements of $u$ are used to obtain the PINN solution. Black dots denote the spatial locations of $u$ measurements.}
\label{fig:PINN_backward_sigmoid_0.2_Nx_30}
\end{figure}

\end{document}